\documentclass[10.5pt]{article}

\def\sqr#1#2{{\vcenter{\vbox{\hrule height.#2pt
              \hbox{\vrule width.#2pt height#1pt \kern#1pt \vrule width.#2pt}
              \hrule height.#2pt}}}}
\def\signed #1{{\unskip\nobreak\hfil\penalty50
              \hskip2em\hbox{}\nobreak\hfil#1
              \parfillskip=0pt \finalhyphendemerits=0 \par}}
\def\endpf{\signed {$\sqr69$}}

\def\dbR{{\mathop{\rm l\negthinspace R}}}

\def\3n{\negthinspace \negthinspace \negthinspace }
\def\2n{\negthinspace \negthinspace }
\def\1n{\negthinspace }

\def\dbE{{\mathop{\rm l\negthinspace E}}}
\def\dbF{{\mathop{\rm l\negthinspace F}}}

\def\dbH{{\mathop{\rm l\negthinspace H}}}

\def\dbP{{\mathop{\rm l\negthinspace P}}}
\def\dbR{{\mathop{\rm l\negthinspace R}}}
\def\={\buildrel \triangle \over =}

\def\ds{\displaystyle}

\def\ns{\noalign{\ss}}

\def\a{\alpha}
\def\b{\beta}
\def\g{\gamma}
\def\d{\delta}
\def\e{\varepsilon}
\def\z{\zeta}
\def\k{\kappa}
\def\l{\lambda}
\def\m{\mu}
\def\n{\nu}
\def\si{\sigma}
\def\t{\tau}
\def\f{\varphi}
\def\th{\theta}

\def\G{\Gamma}
\def\D{\Delta}
\def\Th{\Theta}

\def\F{\Phi}
\def\O{\Omega}

\def\cC{{\cal C}}
\def\cD{{\cal D}}
\def\cE{{\cal E}}
\def\cF{{\cal F}}
\def\cG{{\cal G}}

\def\cL{{\cal L}}

\def\cP{{\cal P}}

\def\cS{{\cal S}}

\def\cU{{\cal U}}

\def\cX{{\cal X}}

\def\cl{{\cal l}}
\def\ss{\smallskip}
\def\ms{\medskip}

\def\q{\quad}
\def\qq{\qquad}
\def\hb{\hbox}

\def\ua{\mathop{\uparrow}}
\def\da{\mathop{\downarrow}}

\def\lan{\mathop{\langle}}
\def\ran{\mathop{\rangle}}

\def\h{\widehat}
\def\wt{\widetilde}
\def\cd{\cdot}
\def\cds{\cdots}

\def\as{\hbox{\rm a.s.{ }}}
\def\sgn{\hbox{\rm sgn$\,$}}

\def\tr{\hbox{\rm tr$\,$}}
\def\cl{\overline}

\def\deq{\mathop{\buildrel\D\over=}}

\def\({\Big (}
\def\){\Big )}
\def\[{\Big[}
\def\]{\Big]}
\def\bde{\begin{definition}}
\def\ede{\end{definition}}
\def\be{\begin{equation}}
\def\bel{\begin{equation}\label}
\def\ee{\end{equation}}
\def\bt{\begin{theorem}}
\def\et{\end{theorem}}
\def\bc{\begin{corollary}}
\def\ec{\end{corollary}}
\def\bl{\begin{lemma}}
\def\el{\end{lemma}}
\def\bp{\begin{proposition}}
\def\ep{\end{proposition}}
\def\bas{\begin{assumption}}
\def\eas{\end{assumption}}
\def\br{\begin{remark}}
\def\er{\end{remark}}
\def\ba{\begin{array}}
\def\ea{\end{array}}
\def\ed{\end{document}}

\def\square#1{\vbox{\hrule\hbox{\vrule height#1%
     \kern#1\vrule}\hrule}}
\def\rectangle#1#2{\vbox{\hrule\hbox{\vrule height#1%
     \kern#2\vrule}\hrule}}

\font\tenbb=msbm10 \font\sevenbb=msbm7 \font\fivebb=msbm5

\newfam\bbfam
\scriptscriptfont\bbfam=\fivebb \textfont\bbfam=\tenbb
\scriptfont\bbfam=\sevenbb

\newtheorem{lemma}{Lemma}[section]
\newtheorem{remark}{Remark}[section]

\newtheorem{theorem}{Theorem}[section]
\newtheorem{corollary}{Corollary}[section]

\newtheorem{definition}{Definition}[section]
\newtheorem{proposition}{Proposition}[section]
\newtheorem{assumption}{Assumption}[section]

\makeatletter
   
   \@addtoreset{equation}{section}
\makeatother

\begin{document}
\title{\bf Time-Inconsistent
Optimal Control Problems\\ and the Equilibrium HJB
Equation\footnote{This work is supported in part by NSF Grant
DMS-1007514. Part of the work was done while the author was visiting
Department of Applied Mathematics, Polytechnic University, Hong
Kong. The author would like to thank Professors X.~Li and J.~Huang
for their hospitality.}}

\author{Jiongmin Yong\\ Department of
Mathematics, University of Central Florida, Orlando, FL 32816, USA}

\maketitle

\begin{abstract} A general time-inconsistent optimal control problem
is considered for stochastic differential equations with
deterministic coefficients. Under suitable conditions, a
Hamilton-Jacobi-Bellman type equation is derived for the equilibrium
value function of the problem. Well-posedness and some properties of
such an equation is studied, and time-consistent equilibrium
strategies are constructed. As special cases, the linear-quadratic
problem and a generalized Merton's portfolio problem are
investigated.

\end{abstract}


\ms


\section{Introduction}

Let $(\O,\cF,\dbF,\dbP)$ be a complete filtered probability space on
which a $d$-dimensional standard Brownian motion $W(\cd)$ is
defined, whose natural filtration is $\dbF=\{\cF_t\}_{t\ge0}$
(augmented by all the $\dbP$-null sets). Let $T>0$. For any
$t\in[0,T)$, we consider the following controlled stochastic
different equation (SDE, for short):
\bel{1.1}\left\{\ba{ll}
\ns\ds dX(s)=b(s,X(s),u(s))ds+\si(s,X(s),u(s))dW(s),\qq s\in[t,T],\\
\ns\ds X(t)=x,\ea\right.\ee
where $b:[0,T]\times\dbR^n\times U\to\dbR^n$ and
$\si:[0,T]\times\dbR^n\times U\to\dbR^{n\times d}$ are suitable
deterministic maps with $U$ being a metric space. In the above,
$(t,x)\in[0,T)\times\dbR^n$ is called an {\it initial pair},
$u:[t,T]\times\O\to U$ is called a {\it control process} and
$X:[0,T]\times\O\to\dbR^n$ is called a {\it state process}. We
define the set of all {\it admissible controls} by the following:
\bel{cU}\cU[t,T]=\Big\{u:[t,T]\times\O\to U\bigm|u(\cd)\hb{ is
$\dbF$-progressively measurable}\Big\}.\ee
Under some mild conditions, for any $(t,x)\in[0,T)\times\dbR^n$ and
$u(\cd)\in\cU[t,T]$, (\ref{1.1}) admits a unique solution
$X(\cd)\equiv X(\cd\,;t,x,u(\cd))$. To measure the performance of
the control process $u(\cd)\in\cU[t,T]$, we introduce the following
cost functional
\bel{1.2}J^0(t,x;u(\cd))=\dbE_t\[\int_t^Te^{-\d(s-t)}g^0(s,X(s),u(s))ds
+e^{-\d(T-t)}h^0(X(T))\],\ee
with some constant $\d\ge0$ (called the {\it discount rate}), some
deterministic maps $g^0:[0,T]\times\dbR^n\times U\to\dbR$ and
$h^0:\dbR^n\to\dbR$, and $\dbE_t[\,\cd\,]=\dbE[\,\cd\,|\,\cF_t]$. On
the right hand side of (\ref{1.2}), the first term is referred to as
a {\it running cost} and the second term is referred to as a {\it
terminal cost}. We can pose the following optimal control problem.

\ms

\bf Problem (C). \rm For any $(t,x)\in[0,T)\times\dbR^n$, find a
$\bar u(\cd)\in\cU[t,T]$ such that
\bel{V0}J^0(t,x;\bar
u(\cd))=\inf_{u(\cd)\in\cU[t,T]}J^0(t,x;u(\cd))\equiv V^0(t,x).\ee

\ms

Any $\bar u(\cd)\in\cU[t,T]$ satisfying (\ref{V0}) is called an {\it
optimal control} of Problem (C) for the initial pair $(t,x)$, the
corresponding state process $\bar X(\cd)\equiv X(\cd\,;t,x,\bar
u(\cd))$ and the pair $(\bar X(\cd),\bar u(\cd))$ are called the
corresponding {\it optimal state process} and {\it optimal pair},
respectively. The function $V^0(\cd\,,\cd)$ defined by (\ref{V0}) is
called the {\it value function} of Problem (C).

\ms

For Problem (C), we have the following Bellman's principle of
optimality (\cite{Fleming-Soner 2006, Yong-Zhou 1999}): For any
$\t\in[t,T]$,
\bel{DP}\ba{ll}
\ns\ds V^0(t,x)=\inf_{u(\cd)\in\cU[t,\t]}\dbE_t\[\int_t^\t
e^{-\d(s-t)}g^0(s,X(s),u(s))ds+e^{-\d(\t-t)}
V^0\big(\t,X(\t;t,x,u(\cd))\big)\],\ea\ee
where $\cU[t,\t]$ is defined similar to that of $\cU[t,T]$ replacing
$[t,T]$ by $[t,\t]$ (see (\ref{cU})). Now, if $(\bar X(\cd),\bar
u(\cd))$ is an optimal pair of Problem (C) for the initial pair
$(t,x)\in[0,T)\times\dbR^n$, then from (\ref{DP}), we obtain
\bel{}\ba{ll}
\ns\ds V^0(t,x)=J^0(t,x;\bar u(\cd))\\
\ns\ds=\dbE_t\[\int_t^\t e^{-\d(s-t)}g^0(s,\bar X(s),\bar
u(s))ds+e^{-\d(\t-t)}J^0\big(\t,\bar X(\t;t,x,u(\cd));\bar u(\cd)
\big|_{[\t,T]}\big)\]\\
\ns\ds\ge\dbE_t\[\int_t^\t e^{-\d(s-t)}g^0(s,\bar X(s),\bar
u(s))ds+e^{-\d(\t-t)}V^0\big(\t,\bar X(\t;t,x,\bar
u(\cd))\big)\]\\
\ns\ds\ge\inf_{u(\cd)\in\cU[t,\t]}\dbE_t\[\int_t^\t
e^{-\d(s-t)}g^0(s,X(s),u(s))ds+e^{-\d(\t-t)}V^0\big(\t,X(\t;t,x,u(\cd))\big)\]\1n
=\1n V^0(t,x).\ea\ee
Thus, one must have
$$\dbE_t\[J^0\big(\t,\bar X(\t);u(\cd)\big|_{[\t,T]}\big)-V^0(\t,\bar
X(\t))\]=0,\qq\as$$
Since
$$J^0\big(\t,\bar X(\t);u(\cd)\big|_{[\t,T]}\big)-V^0(\t,\bar
X(\t))\ge0,\qq\as,$$
it follows that
\bel{}J^0\big(\t,\bar X(\t);u(\cd)\big|_{[\t,T]}\big)=V^0(\t,\bar
X(\t))=\inf_{u(\cd)\in\cU[\t,T]}J^0\big(\t,\bar
X(\t);u(\cd)\big),\qq\as\ee
This means that the restriction $\bar
u(\cd)\big|_{[\t,T]}\in\cU[\t,T]$ of an optimal control $\bar
u(\cd)\in\cU[t,T]$ for the initial pair $(t,x)$ on any later time
interval $[\t,T]$ is optimal for the initial pair $\big(\t,\bar
X(\t;t,x,\bar u(\cd))\big)$. Such a phenomenon is called the {\it
time-consistency} of Problem (C).

\ms

The advantage of the time-consistency is that for any given initial
pair $(t,x)$, if an optimal control $\bar u(\cd)$ can be constructed
for that (initial pair), then it will stay optimal hereafter (for
the later initial pair along the optimal trajectory). This is very
ideal. However, in reality, the time-consistency could be lost. An
interesting situation that we are going to discuss in this paper is
what we call the {\it general discounting} situation which includes
the so-called {\it non-exponential discounting}, or {\it hyperbolic
discounting} situations. This amounts to saying the following: Due
to the possible subjectivity of people's preferences, the discount
factors $e^{-\d(s-t)}$ and $e^{-\d(T-t)}$ appeared in (\ref{1.2})
might be replaced by some general functions $\l(s,t)$ and $\n(T,t)$,
or more generally, we will consider the following cost functional:
\bel{1.8}J(t,x;u(\cd))=\dbE_t\[\int_t^Tg(t,s,X(s),u(s))ds+h(t,X(T))\],\ee
where the maps $g(\cd)$ and $h(\cd)$ explicitly depend on the
initial time $t$ in some general way. The optimal control problem
associated with (\ref{1.1}) and (\ref{1.8}), called Problem (N),
will not be time-consistent, or {\it time-inconsistent}, in general,
meaning that a restriction of an optimal control for a specific
initial pair on a later time interval might not be optimal for that
corresponding initial pair. Some concrete examples will be presented
in the next section.

\ms

In recent years, time-inconsistent optimal control problems have
attracted a number of researchers. See \cite{Ekeland-Lazrak 2010,
Ekeland-Privu 2008a, Bjork-Murgoci 2008, Yong 2011,
Marin-Solano-Shevkoplyas, Yong 2012, Ekeland-Mbodji-Pirvu 2012,
Bjork-Murgoci-Zhou, Hu-Jin-Zhou 2011} and references cited therein.

\ms

The purpose of this paper is to obtain time-consistent optimal
controls (which should be more properly called {\it equilibrium
control}) for Problem (N) mentioned above. Let us now briefly
describe our approach. Inspired by \cite{Yong 2011, Yong 2012}, we
introduce a sequence of multi-person hierarchical differential games
as follows. For any $N>1$, let $\Pi$ be a partition of the time
interval $[0,T]$ defined by
$$\Pi:0=t_1<t_2<\cds<t_N=T,$$
with the mesh size $\|\Pi\|$ given by
$$\|\Pi\|=\max_{1\le k\le N}(t_k-t_{k-1}).$$
The differential game, denoted by Problem (G$^\Pi$), associated with
partition $\Pi$ consists of $N$ players. The $k$-th player controls
the system on $[t_{k-1},t_k)$ by taking his/her control
$u^k(\cd)\in\cU[t_{k-1},t_k]$. The cost functional is constructed in
a sophisticated way, by using some techniques of forward-backward
stochastic differential equations (FBSDEs, for short) found in
\cite{MPY 1994,MY 1999}. The interaction among the players are as
follows: (i) The terminal pair $(t_k,X(t_k))$ of Player $k$ is the
initial pair of Player $(k+1)$; (ii) All the player know that each
player tries to find an optimal control for his/her own problem; and
(iii) Each player will discount the future costs in his/her own way,
regardless of the fact that the later players will control the
system. Under certain conditions, each player will have an optimal
control, denoted by $\bar u^k(\cd)\in\cU[t_{k-1},t_k]$, for his/her
own problem, as well as his/her own value function $V^k(\cd\,,\cd)$
defined on $[t_{k-1},t_k]\times\dbR^n$. Define
\bel{}\bar u^\Pi(t)=\sum_{k=1}^N\bar u^k(t)I_{[t_{k-1},t_k)}(t),\qq
t\in[0,T),\ee
and
\bel{}V^\Pi(t,x)=\sum_{k=1}^NV^k(t,x)I_{[t_{k-1},t_k)}(t),\qq(t,x)\in[0,T)
\times\dbR^n.\ee
We may call $\bar u^\Pi(\cd)$ and $V^\Pi(\cd\,,\cd)$ the {\it Nash
equilibrium control} and {\it Nash equilibrium value function} of
Problem (G$^\Pi$), respectively. When the following limits
\bel{}\left\{\ba{ll}
\ns\ds\lim_{\|\Pi\|\to0}\bar u^\Pi(t)=\bar u(t),\qq\qq t\in[0,T],\\
\ns\ds\lim_{\|\Pi\|\to0}V^\Pi(t,x)=V(t,x),\qq(t,x)\in[0,T]\times\dbR^n,\ea\right.\ee
exist for some $\bar u(\cd)\in\cU[0,T]$ and
$V:[0,T]\times\dbR^n\to\dbR$, we call them a {\it time-consistent
equilibrium control} and a {\it time-consistent equilibrium value
function} of Problem (N), respectively. As a major contribution of
this paper, we have derived the {\it equilibrium
Hamilton-Jacobi-Bellman equation} which can be used to characterize
the equilibrium value function $V(\cd\,,\cd)$, and it recovers the
result for time-inconsistent deterministic linear-quadratic problem
presented in \cite{Yong 2011}. The well-posedness of such an HJB
equation will be established for the case that the diffusion of the
state equation does not contain the control. The general case is
open at the moment, and we expect to present some more complete
results in our future publications. As important and interesting
special cases, we will construct equilibrium controls for stochastic
LQ problem with general discounting and for generalized Merton's
portfolio problem.

\ms

We refer the readers to \cite{Strotz,Pollak 1968,Goldman
1980,Miller-Salmon 1985,Tesfatsion 1986,L1997,Karp-Lee
2003,Palacios-Huerta,Caplin-Leahy 2006,Herings-Rohde,Grenadier-Wang
2007,Basak-Chabakarui 2008}, for some relevant results.

\section{Two Examples of Time-Inconsistent Optimal Control Problems}

In this section, we present two interesting examples of optimal
control problems which are time-inconsistent.

\ms

\bf Example 2.1. \rm Consider a one-dimensional controlled linear
SDE:
\bel{ex2.1a}\left\{\ba{ll}
\ns\ds dX(s)=u(s)ds+\si X(s)dW(s),\qq s\in[t,T],\\
\ns\ds X(t)=x,\ea\right.\ee
with cost functional
\bel{ex2.1b}J(t,x;u(\cd))=\dbE_t\[\int_t^T|u(s)|^2ds+g(t)|X(T)|^2\],\ee
where $\si>0$ is a constant and $g:[0,T]\to(0,\infty)$ is a
deterministic non-constant, continuous function. For a fixed initial
pair $(t,x)\in[0,T)\times\dbR$, the LQ problem associated with the
above (\ref{ex2.1a})--(\ref{ex2.1b}) is a standard one. We can show
that (see Appendix) the problem admits a unique optimal pair $(\bar
X(\cd),\bar u(\cd))$ given by
\bel{ex2.1c}\left\{\ba{ll}
\ns\ds\bar X(s)={\si^2+g(t)(e^{\si^2(T-s)}-1)\over\si^2+g(t)
(e^{\si^2(T-t)}-1)}\,e^{-{\si^2\over2}(s-t)+\si[W(s)-W(t)]}x,\\
\ns\ds\bar u(s)=-{\si^2g(t)e^{\si^2(T-s)}\over
\si^2+g(t)(e^{\si^2(T-t)}-1)}\,e^{-{\si^2\over2}(s-t)+\si[W(s)-W(t)]}x,\ea\right.\qq
s\in[t,T].\ee
Now, let $\t\in(t,T)$. We consider the corresponding LQ problem
starting from the initial pair $(\t,\bar X(\t))$. Then the optimal
pair, denoted by $(\h X(\cd),\h u(\cd))$, is given by
\bel{}\left\{\ba{ll}
\ns\ds\h X(s)={\si^2+g(\t)(e^{\si^2(T-s)}-1)\over\si^2+g(\t)
(e^{\si^2(T-\t)}-1)}\,e^{-{\si^2\over2}(s-\t)+\si[W(s)-W(\t)]}\bar
X(\t)\\
\ns\ds\qq~={\si^2+g(\t)(e^{\si^2(T-s)}-1)\over\si^2+g(t)
(e^{\si^2(T-s)}-1)}\cd{\si^2+g(t)(e^{\si^2(T-\t)}-1)\over
\si^2+g(\t)(e^{\si^2(T-\t)}-1)}\bar X(s),\\
\ns\ds\h u(s)=-{\si^2g(\t)e^{\si^2(T-s)}\over
\si^2+g(\t)(e^{\si^2(T-\t)}-1)}\,e^{-{\si^2\over2}(s-\t)+\si[W(s)-W(\t)]}\bar
X(\t)\\
\ns\ds\qq={g(\t)\over
g(t)}\cd{\si^2+g(t)(e^{\si^2(T-\t)}-1)\over\si^2+g(\t)(e^{\si^2(T-\t)}-1)}\,\bar
u(s),\ea\right.\qq s\in[\t,T].\ee
Hence, in general
$$\h X(s)=\bar X(s),\q\h u(s)=\bar u(s),\qq s\in[\t,T],$$
may fail. In fact, we have the following comparison of the cost
functional values corresponding to $\bar u(\cd)\big|_{[\t,T]}$ and
$\h u(\cd)$, with the same initial pair $(\t,\bar X(\t))$:
$$\ba{ll}
\ns\ds J\big(\t,\bar X(\t);\bar u(\cd)\big|_{[\t,T]}\big)-J(\t,\bar
X(\t);\h
u(\cd))={\si^4(e^{\si^2(T-\t)}-1)e^{\si^2(T-\t)}[g(t)-g(\t)]^2|\bar
X(\t)|^2\over[\si^2+g(t)(e^{\si^2(T-\t)}-1)]^2[\si^2+g(\t)(e^{\si^2(T-\t)}-1)]}\\
\ns\ds={\si^4(e^{\si^2(T-\t)}-1)[g(t)-g(\t)]^2
e^{\si^2(T-\t)-\si^2(\t-t)+2\si[W(\t)-W(t)]}x^2\over[\si^2+g(t)(e^{\si^2(T-t)}-1)]^2
[\si^2+g(\t)(e^{\si^2(T-\t)}-1)]} \,>0,\ea$$
unless $g(\t)=g(t)$ or $x=0$. Therefore, the LQ problem associated
with (\ref{ex2.1a})--(\ref{ex2.1b}) is time-inconsistent.

\ms

Note that by sending $\si^2\to0$, using l'H\^opital's rule, we
recover the example presented in \cite{Yong 2011, Yong 2012} for a
deterministic LQ problem.

\ms

\bf Example 2.2. (Merton's portfolio problem) \rm We consider a
one-dimensional controlled SDE:
\bel{ex2.2a}\left\{\ba{ll}
\ns\ds dX(s)=\big[rX(s)+(\m-r)u(s)-c(s)\big]ds+\si u(s)dW(s),\qq
s\in[t,T],\\
\ns\ds X(t)=x,\ea\right.\ee
where $(u(\cd),c(\cd))$ is the control process with $u(\cd)$ being
the dollar amount invested in the stock (which could be positive and
negative), $c(\cd)$ is the consumption rate process (which has to be
non-negative), and the solution $X(\cd)\equiv
X(\cd\,;t,x,u(\cd),c(\cd))$ of (\ref{ex2.2a}) corresponding to
$(t,x,u(\cd),c(\cd))$ is the wealth process, which is required to be
non-negative. In the above, $r>0$ is the interest rate for the bank
account, $\m>r$ and $\si$ are the appreciation rate and volatility
of the stock, respectively. Let
$$\ba{ll}
\ns\ds\cU[t,T]=\Big\{u:[t,T]\times\O\to\dbR\bigm|u(\cd)\hb{ is
$\dbF$-progressively
measurable, }\dbE\[\int_t^T|u(s)|^2ds\]^{1\over2}<\infty\Big\},\\
\ns\ds\cC[t,T]=\Big\{c:[t,T]\times\O\to[0,\infty)\bigm|c(\cd)\hb{ is
$\dbF$-progressively measurable,
}\dbE\int_t^Tc(s)ds<\infty\Big\}.\ea$$
The payoff functional is given by the following:
\bel{}J(t,x;u(\cd),c(\cd))=\dbE_t\[\int_t^T\n(t,s)c(s)^\b
ds+\rho(t)X(T)^\b\],\ee
where $\n(\cd\,,\cd)$ and $\rho(\cd)$ are given positive-valued
functions, and $\b\in(0,1)$. As a convention, we define
\bel{}x^\b=-\infty,\qq x<0.\ee
The following is referred to as a {\it generalized Merton's
portfolio problem}:

\ms

\bf Problem (M). \rm For any given $(t,x)\in[0,T)\times\dbR$, find a
pair $(\bar u(\cd),\bar c(\cd))\in\cU[t,T]\times\cC[t,T]$ such
that
\bel{}J(t,x;\bar u(\cd),\bar
c(\cd))=\sup_{(u(\cd),c(\cd))\in\cU[t,T]\times\cC[t,T]}
J(t,x;u(\cd),c(\cd)).\ee

\ms

When
\bel{classical}\n(t,s)=e^{-\d(s-t)},\q\rho(t)=e^{-\d(T-t)},\qq0\le
t\le s\le T,\ee
the above problem is called a {\it classical Merton's portfolio
problem} (which is time-consistent). We now look at the above
general case. For a fixed $t\in[0,T)$, Problem (M) can be regarded
as a standard optimal control problem. Therefore, we can treat it in
a standard way. We can show that (see Appendix) the problem admits a
unique optimal control $(\bar u^t(\cd),\bar c^t(\cd))$ for the
initial pair $(t,x)$, which is given by the following:
\bel{}\left\{\ba{ll}
\ns\ds\bar u^t(s)={(\m-r)\bar X^t(s)\over\si^2(1-\b)},\\
\ns\ds\bar c^t(s)={\n(t,s)^{1\over1-\b}\bar X^t(s)\over
e^{{\l\over1-\b}(T-s)}\rho(t)^{1\over1-\b}
+\int_s^Te^{{\l\over1-\b}(\t-s)}\n(t,\t)^{1\over1-\b}d\t},\ea\right.\qq
s\in[t,T],\ee
with
\bel{l2}\l={[2r\si^2(1-\b)+(\m-r)^2]\b\over2\si^2(1-\b)},\ee
and
\bel{}\ba{ll}
\ns\ds
V^t(t,x)\equiv\sup_{(u(\cd),c(\cd))\in\cU[t,T]\times\cC[t,T]}J(t,x;u(\cd),c(\cd))\\
\ns\ds\qq\q~=\[e^{{\l\over1-\b}(T-t)}\rho(t)^{1\over1-\b}
+\int_t^Te^{{\l\over1-\b}(\t-t)}\n(t,\t)^{1\over1-\b}
d\t\]^{1-\b}x^\b,\qq(t,x)\in[0,T]\times[0,\infty).\ea\ee
Now, let $\bar t\in(t,T)$, and consider Problem (M) starting from
the initial pair $(\bar t,\bar X^t(\bar t))$. Then
$$\ba{ll}
\ns\ds V^{\bar t}(\bar t,\bar X^t(\bar t))=\[e^{{\l\over1-\b}(T-\bar
t)}\rho(\bar t)^{1\over1-\b}+\int_{\bar t}^Te^{{\l\over1-\b}(\t-\bar
t)}\n(\bar t,\t)^{1\over1-\b} d\t\]^{1-\b}\bar X^t(\bar t)^\b.\ea$$
We can show that (see Appendix) if $\bar t\in(t,T)$ is such that
\bel{2.13b}\int_{\bar
t}^T\[{e^{\l\t}\n(t,\t)\over\rho(t)}\]^{1\over1-\b}d\t\ne\int_{\bar
t}^T\[{e^{\l\t}\n(\bar t,\t)\over\rho(\bar t)}\]^{1\over1-\b}d\t,\ee
then
\bel{}J\big(\bar t,\bar X^t(\bar t);\bar u(\cd)\big|_{[\bar
t,T]},\bar c(\cd)\big|_{[\bar t,T]}\big)<V^{\bar t}(\bar t,\bar
X^t(\bar t))=\sup_{(u(\cd),c(\cd))\in\cU[\bar t,T]\times\cC[\bar
t,T]}J(\bar t,\bar X^t(\bar t);u(\cd),c(\cd)),\ee
which means that the restriction $\big(\bar u(\cd)\big|_{[\bar
t,T]},\bar c(\cd)\big|_{[\bar t,T]}\big)$ of optimal control $(\bar
u(\cd),\bar c(\cd))$ for the initial pair $(t,x)$ on $[\bar t,T]$ is
not optimal for the initial pair $(\bar t,\bar X^t(\bar t))$ if
(\ref{2.13b}) holds.

\ms

Interestingly, in the case that (\ref{classical}) holds, one has
$${\n(t,\t)\over\rho(t)}={e^{-\d(\t-t)}\over e^{-\d(T-t)}}=e^{\d(T-\t)}.$$
Thus, in this case, (\ref{2.13b}) does not hold. As a matter of
fact, the problem is time-consistent and (\ref{2.13b}) should not be
true.

\ms

\section{Some Preliminaries}

\ms

For convenience, let us rewrite the state equation and the cost
functional below.
\bel{3.1}\left\{\ba{ll}
\ns\ds dX(s)=b(s,X(s),u(s))ds+\si(s,X(s),u(s))dW(s),\qq s\in[t,T],\\
\ns\ds X(t)=x,\ea\right.\ee
and
\bel{3.2}
J(t,x;u(\cd))=\dbE_t\[\int_t^Tg(t,s,X(s),u(s))ds+h(t,X(T))\].\ee
Clearly, our cost functional covers the non-exponential/hyperbolic
discounting situations. By comparing the above state equation and
cost functional, it seems that we may consider a little more general
state equation of the following form:
\bel{3.3}\left\{\ba{ll}
\ns\ds dX(s)=b(t,x,s,X(s),u(s))ds+\si(t,x,s,X(s),u(s))dW(s),\qq s\in[t,T],\\
\ns\ds X(t)=x.\ea\right.\ee
However, according to \cite{Yong 2012} (see also \cite{Yong 2011}),
we know that when an equilibrium pair (see below for definition) is
constructed, the eventual effective state equation will take the
following form:
\bel{}\left\{\ba{ll}
\ns\ds dX(s)=b(s,X(s),s,X(s),u(s))ds+\si(s,X(s),s,X(s),u(s))dW(s),\qq s\in[t,T],\\
\ns\ds X(t)=x,\ea\right.\ee
Therefore, it suffices to consider the state equation of form
(\ref{3.1}).

\ms

In what follows, we let $T>0$ be a fixed time horizon, and
$U\subseteq\dbR^m$ be a closed subset, which could be either bounded
or unbounded (it is allowed that $U=\dbR^m$). We will use $K>0$ as a
generic constant which can be different from line to line. Let
$\cS^n$ be the set of all $(n\times n)$ symmetric real matrices.
Denote
$$D[0,T]=\Big\{(t,s)\in[0,T]^2\bigm|0\le t\le s\le T\Big\}.$$
Recall from Section 1 that
$$\cU[t,T]=\Big\{u:[t,T]\times\O\to U\bigm|u(\cd)\hb{ is $\dbF$-progressively
measurable}\Big\}.$$
Further, for $q\ge1$, let
$$\cU^q[t,T]=\Big\{u:[t,T]\times\O\to U\bigm|u(\cd)\hb{ is $\dbF$-progressively measurable, }
\dbE\int_t^T|u(s)|^qds<\infty\Big\}.$$
Note that in the case $U$ is bounded, for different $q\ge1$, all the
$\cU^q[t,T]$ are the same as $\cU[t,T]$. We introduce the following
standing assumptions.

\ms

{\bf(H1)} The maps $b:[0,T]\times\dbR^n\times U\to\dbR^n$,
$\si:[0,T]\times\dbR^n\times U\to\dbR^{n\times d}$ are continuous
and there exist constants $L>0$ and $k\ge0$ such that
\bel{3.5}\left\{\ba{ll}
\ns\ds|b(t,x,u)-b(t,y,u)|\le
L\big(1+(|x|\vee|y|)^k+|u|\big)|x-y|,\\
\ns\ds\lan x-y,b(t,x,u)-b(t,y,u)\ran\le
L|x-y|^2,\\
\ns\ds|\si(t,x,u)-\si(t,y,u)|\le
L|x-y|,\ea\right.\qq\forall(t,u)\in[0,T]\times U,~x,y\in\dbR^n,\ee
where $|x|\vee|y|=\max\{|x|,|y|\}$, and
\bel{}|b(t,0,u)|+|\si(t,0,u)|\le
L(1+|u|),\qq\qq\forall(t,u)\in[0,T]\times U.\ee

\ms

{\bf(H2)} Maps $g:D[0,T]\times\dbR^n\times U\to\dbR$ and
$h:[0,T]\times\dbR^n\to\dbR$ are continuous, and there exist
constants $L>0$ and $q\ge0$ such that
\bel{3.7}\left\{\ba{ll}
\ns\ds0\le g(\t,t,x,u)\le L(1+|x|^q+|u|^q\big),\\
\ns\ds0\le h(\t,x)\le
L(1+|x|^q),\ea\right.\qq\qq\forall(\t,t,x,u)\in
D[0,T]\times\dbR^n\times U.\ee

\ms

Let us make a couple of remarks on (H1). First of all, if $x\mapsto
b(t,x,u)$ is uniformly Lipschitz, then the first two conditions of
(\ref{3.5}) hold. On the other hand, we point out that the first
condition in (\ref{3.5}) merely implies that $x\mapsto b(t,x,u)$ is
locally Lipschitz, and the second condition in (\ref{3.5}) alone
does not imply the global Lipschitz condition for the map $x\mapsto
b(t,x,u)$. A simple example that the first and the second conditions
in (\ref{3.5}) are satisfied but $x\mapsto b(t,x,u)$ is not uniform
Lipschitz is the following:
$$b(t,x,u)=-|x|^2x-|u|x,\qq(t,x,u)\in[0,T]\times\dbR^n\times U.$$
It is clear that the above map is not uniformly Lipschitz in $x$,
the first condition in (\ref{3.5}) holds with $k=2$, and we can
check that
$$\lan x-y,b(t,x,u)-b(t,y,u)\ran\le0,\qq\forall(t,u)\in[0,T]\times U,~x,y\in\dbR^n.$$
Note that under (\ref{3.5}), one has
$$\ba{ll}
\ns\ds\lan x,b(t,x,u)\ran=\lan x,b(t,x,u)-b(t,0,u)\ran+\lan
x,b(t,0,u)\ran\\
\ns\ds\qq\qq\qq\le
L|x|^2+L(1+|u|)|x|\le2L\big(1+|x|^2+|u|^2\big),\ea$$
and
$$|\si(t,x,u)|^2\le L^2\big(1+|x|+|u|\big)^2\le3L^2\big(1+|x|^2+|u|^2\big).$$
Now, for (H2), we note that the nonnegativity of $g(\cd)$ and
$h(\cd)$ can be replaced by the condition that both $g(\cd)$ and
$h(\cd)$ are bounded from below. The following result is concerning
the well-posedness of the state equation.

\ms

\bf Proposition 3.1. \sl Let {\rm(H1)} hold. Then for any
$(t,x)\in[0,T]\times\dbR^n$ and $u(\cd)\in\cU^q[t,T]$, with $q\ge2$,
state equation $(\ref{3.1})$ admits a unique solution $X(\cd)\equiv
X(\cd\,;t,x,u(\cd))$, and the following estimate hold:
\bel{3.8}\dbE|X(s;t,x,u(\cd))|^q\le
K\big(1+|x|^q+\int_t^s|u(r)|^qdr\big),\qq s\in[t,T].\ee
Moreover, if $\h x\in\dbR^n$ is another point, the following
holds:
\bel{3.9}\dbE|X(s;t,x,u(\cd)))-X(s;t,\h x,u(\cd))|^2\le K|x-\h
x|^2,\qq\forall t\in[0,T],~x,\h x\in\dbR^n,~u(\cd)\in\cU^2[t,T].\ee

\ms

\it Proof. \rm Let $(t,x)\in[0,T)\times\dbR^n$ be given and
$u(\cd)\in\cU^q[t,T]$. By the local Lipschitz condition, we see that
the solution $X(\cd)\equiv X(\cd\,;t,x,u(\cd))$ exists locally, say
on $[t,\t)$ for some stopping time $\t$. Next, applying It\^o's
formula to $|X(\cd)|^q$ on $[t,\t)$, we have the following:
$$\ba{ll}
\ns\ds\dbE|X(s)|^q=|x|^q+\dbE\int_t^s\(q|X(r)|^{q-2}\lan
X(r),b(r,X(r),u(r))\ran\\
\ns\ds+{q(q-2)\over2}|X(r)|^{q-4}\lan
X(r),\si(r,X(r),u(r))\ran{}^2+{q\over2}|X(r)|^{q-2}\tr\big[\si(s,X(s),u(s))
\si(s,X(s),u(s))^T\big]\)dr\\
\ns\ds\le|x|^q+K\dbE\int_t^s|X(r)|^{q-2}\big(1+|X(r)|^2+|u(r)|^2\big)dr
\le|x|^q+K\dbE\int_t^s\big(1+|X(r)|^q+|u(r)|^q\big)dr.\ea$$
By Gronwall's inequality, we obtain
$$\dbE|X(s)|^q\le K\big(1+|x|^q+\dbE\int_t^s|u(r)|^qdr\big),\qq s\in[t,\t).$$
This implies that $X(\cd)$ must be globally exists on $[t,T]$, and
estimate (\ref{3.8}) holds.

\ms

Next, for any $t\in[0,T)$, $x,\h x\in\dbR^n$, and
$u(\cd)\in\cU^q[t,T]$, let $X(\cd)=X(\cd\,;t,x,u(\cd))$ and $\h
X(\cd)=X(\cd\,;t,\h x,u(\cd))$. Then applying It\^o's formula to
$|X(\cd)-\h X(\cd)|^2$, we get
$$\ba{ll}
\ns\ds\dbE|X(s)-\h X(s)|^2=|x-\h x|^2+\dbE\int_t^s\(\lan X(s)-\h
X(s),b(s,X(s),u(s))-b(s,\h X(s),u(s))\ran\\
\ns\ds\qq\qq\qq\qq\qq\qq\qq\qq\qq+|\si(s,X(s),u(s))-\si(s,\h
X(s),u(s))|^2\)ds\\
\ns\ds\qq\qq\qq\qq\le|x-\h x|^2+K\dbE\int_t^s|X(s)-\h
X(s)|^2ds.\ea$$
Hence, applying Gronwall's inequality, we obtain (\ref{3.9}). \endpf

\ms

From the above result, we see that for any
$(t,x)\in[0,T)\times\dbR^n$ and $u(\cd)\in\cU^{q\vee2}[t,T]$, the
following holds:
$$\ba{ll}
\ns\ds\dbE|g(\t,s,X(s),u(s))|\le
L\big(1+\dbE|X(s)|^{q\vee2}+\dbE|u(s)|^{q\vee2}\big)\\
\ns\ds\qq\qq\qq\qq\q\le
K\big(1+|x|^{q\vee2}+\dbE|u(s)|^{q\vee2}+\int_t^s\dbE|u(r)|^{q\vee2}dr\big),\ea$$
and
$$\dbE|h(\t,X(T))|\le L\big(1+\dbE|X(T)|^{q\vee2}\big)\le
K\big(1+|x|^{q\vee2}+\dbE\int_t^T|u(s)|^{q\vee2}ds\big).$$
Hence, $J(t,x;u(\cd))$ is finite for any
$u(\cd)\in\cU^{q\vee2}[t,T]$. For simplicity, hereafter, we adopt
the following convention:
\bel{convention}J(t,x;u(\cd))\deq+\infty,\qq\hb{if $J(t,x,u(\cd))$
is not finite or not defined}.\ee
We now formally state our optimal control problem.

\ms

\bf Problem (N). \rm For any given initial pair
$(t,x)\in[0,T)\times\dbR^n$, find a $\bar u(\cd)\in\cU[t,T]$ such
that
\bel{}J(t,x;\bar u(\cd))=\inf_{u(\cd)\in\cU[t,T]}J(t,x;u(\cd)).\ee

\ms

From the examples presented in the previous section, we know that
the above Problem (N) is time-inconsistent, in general. Our goal is
to find time-consistent equilibrium controls and characterize the
equilibrium value function, which will be made precise below.

\ms

We denote
$$a(t,x,u)={1\over2}\si(t,x,u)\si(t,x,u)^T,\qq\forall(t,x,u)\in[0,T]\times\dbR^n\times
U.$$
Define
\bel{dbH}\ba{ll}
\ns\ds\dbH(\t,t,x,u,p,P)=\lan
b(t,x,u),p\ran+\tr\big[a(t,x,u)P\big]+g(\t,t,x,u),\\
\ns\ds\qq\qq\qq\qq\qq\qq\qq\forall(\t,t,x,u,p,P)\in
D[0,T]\times\dbR^n\times U\times\dbR^n\times\cS^n,\ea\ee
and let
\bel{H}H(\t,t,x,p,P)=\inf_{u\in U}\dbH(\t,t,x,u,p,P),\qq\forall
(\t,t,x,p,P)\in D[0,T]\times\dbR^n\times\dbR^n\times\cS^n.\ee
Note that $H(\t,t,x,p,P)$ is not necessarily finite on the whole
space $D[0,T]\times\dbR^n\times\dbR^n\times\cS^n$ when $U$ is
unbounded. Therefore, we denote
\bel{D(H)}\cD(H)=\Big\{(\t,t,x,p,P)\in
D[0,T]\times\dbR^n\times\dbR^n\times\cS^n
\bigm|H(\t,t,x,p,P)>-\infty\Big\},\ee
call it the {\it domain} of $H$. Next, let
\bel{}\ba{ll}
\ns\ds\arg\min\dbH(\t,t,x,\cd\,,p,P)=\Big\{\bar u\in
U\bigm|\dbH(\t,t,x,\bar u,p,P)=\min_{u\in
U}\dbH(\t,t,x,u,p,P)\Big\},\\
\ns\ds\qq\qq\qq\qq\qq\qq\qq\qq\qq\forall(\t,t,x,p,P)\in\cD(H).\ea\ee
This is a multi-valued map. Suppose we can define a map
$\psi:\cD(\psi)\subseteq\cD(H)\to U$ such that
\bel{psi1}\ba{ll}
\ns\ds
H(\t,t,x,p,P)\equiv\dbH(\t,t,x,\psi(\t,t,x,p,P),p,P)=\inf_{u\in
U}\dbH(\t,t,x,u,p,P)>-\infty,\\
\ns\ds\qq\qq\qq\qq\qq\qq\qq\qq\qq\qq\qq\forall(\t,t,x,p,P)\in\cD(\psi).\ea\ee
The set $\cD(\psi)$ is called the {\it domain} of $\psi$, which
consists of all points $(\t,t,x,p,P)\in\cD(H)$ such that the infimum
in (\ref{psi1}) is achieved at $\psi(\t,t,x,p,P)$. It is clear that
$\psi(\cd)$ is actually a selection of $\arg\min\dbH(\cd)$, i.e.,
\bel{psi2}\ba{ll}
\ns\ds\psi(\t,t,x,p,P)\in\arg\min\dbH(\t,t,x,\cd\,,p,P),\qq\forall(\t,t,x,p,P)\in\cD(\psi).\ea\ee
The map $\psi(\cd)$ will play an important role later. Therefore,
let us say a little bit more on it. Note that when $U$ is bounded
(since it is assumed to be closed, it is compact in this case), one
has
$$\cD(\psi)=\cD(H)=D[0,T]\times\dbR^n\times\dbR^n\times\cS^n.$$
However, when $U$ is unbounded, say, $U=\dbR^m$, one might have
\bel{3.16}\cD(\psi)\ne\cD(H)\ne
D[0,T]\times\dbR^n\times\dbR^n\times\cS^n.\ee
To say something about the case of (\ref{3.16}), let us present the
following simple lemma.

\ms

\bf Lemma 3.2. \sl Let $U\subseteq\dbR^m$ be closed and $f:U\to\dbR$
be a lower semi-continuous map such that
$$\inf_{u\in U}f(u)\equiv\bar f>-\infty.$$
For any $\e>0$, let
$$f^\e(u)=f(u)+\e|u|^2,\qq u\in U.$$
Then there exists a $u^\e\in U$ such that
$$f^\e(u^\e)=\inf_{u\in U}f^\e(u)\equiv\bar f^\e.$$
Moreover,
$$\lim_{\e\da0}\bar f^\e=\bar f,\qq\lim_{\e\da0}\e|u^\e|^2=0.$$

\ms

\it Proof. \rm First of all, fix a $u_0\in U$. For any minimizing
sequence $u_k\in U$ of $f^\e(\cd)$, we may assume that
$$f^\e(u_0)\ge f^\e(u_k)=f(u_k)+\e|u_k|^2\ge\bar f+\e|u_k|^2,\qq\forall k\ge1.$$
Thus, $u_k$ is bounded. Consequently, by the closeness of $U$, we
may assume that $u_k\to u^\e\in U$ exists, which attains the infimum
of $f^\e(\cd)$. Next, it is clear that $f^\e(\cd)$ decreases as $\e$
decreases, and
$$\bar f\le\lim_{\e\to0}\bar f^\e\equiv\bar f^0.$$
Now, for any $\d>0$, there exists a $u_\d\in U$ such that
$$f(u_\d)<\bar f+\d.$$
On the other hand,
$$\bar f^\e=f^\e(u^\e)\le f^\e(u_\d)=f(u_\d)+\e|u_\d|^2.$$
Hence, letting $\e\to0$, we get
$$\bar f^0\le\bar f+\d.$$
Since $\d>0$ is arbitrary, we must have $\bar f=\bar f^0$. Finally,
by
$$\bar f\le f(u^\e)\le f(u^\e)+\e|u^\e|^2=f^\e(u^\e)=\bar f^\e\to\bar f,$$
we see that
$$\lim_{\e\to0}\e|u^\e|^2=0,$$
proving our conclusion. \endpf

\ms

Now, for the case $\cD(\psi)\ne\cD(H)$, we may introduce
\bel{}\ba{ll}
\ns\ds\dbH^\e(\t,t,x,u,p,P)=\dbH(\t,t,x,u,p,P)+\e|u|^2,\\
\ns\ds\qq\qq\qq\qq\qq\qq\qq\qq(\t,t,x,u,p,P) \in
D[0,T]\times\dbR^n\times U\times\dbR^n\times\cS^n.\ea\ee
Then by Lemma 3.2, one can find a $\psi^\e:\cD(H)\to U$ such that
\bel{}\ba{ll}
\ns\ds\dbH^\e(\t,t,x,\psi^\e(\t,t,x,p,P),p,P)=\inf_{u\in
U}\dbH^\e(\t,t,x,u,p,P)\\
\ns\ds\qq\qq\qq\qq\qq\qq\q~\equiv
H^\e(\t,t,x,p,P),\qq\qq(\t,t,x,p,P)\in\cD(H).\ea\ee
Moreover, $\e\mapsto H^\e(\t,t,x,p,P)$ is decreasing as $\e\da0$,
and
$$\left\{\ba{ll}
\ns\ds\lim_{\e\to0}H^\e(\t,t,x,p,P)=H(\t,t,x,p,P),\\
\ns\ds\lim_{\e\to0}\e|\psi^\e(\t,t,x,p,P)|^2=0,\ea\right.\qq
(\t,t,x,p,P)\in\cD(H).$$
In general, we do not expect the convergence of
$\psi^\e(\t,t,x,p,P)$ as $\e\to0$.

\ms

We point out that in general the map $\psi:\cD(\psi)\subseteq
D[0,T]\times\dbR^n\times\dbR^n\times\cS^n\to U$ is not necessarily
continuous. Here is a simple example.

\ms

\bf Example 3.3. \rm Let $n=m=d=1$, $U=[0,1]$, and
$$b(t,x,u)=u,\q\si(t,x,u)=\si(t,x),\q g(\t,t,x,u)=R(\t,t)u,$$
with $\si(\cd\,,\cd)$ and $R(\cd\,,\cd)$ being continuous, and
$R(\t,t)>0$ for all $(\t,t)\in D[0,T]$. Then
$$\dbH(\t,t,x,u,p,P)=pu+{1\over2}\si(t,x)^2P+R(\t,t)u={1\over2}\si(t,x)^2P
+\big[p+R(\t,t)\big]u.$$
Consequently,
$$H(\t,t,x,p,P)={1\over2}\si(t,x)^2P+\min\{p+R(\t,t),0\},\qq(\t,t,x,p,P)\in D[0,T]\times
\dbR^n\times\dbR^n\times\cS^n,$$
which is continuous, and
$$\psi(\t,t,x,p,P)=I_{(p+R(\t,t)<0)},\qq\forall(\t,t,x,p,P)\in D[0,T]\times
\dbR^3,$$
which has discontinuity along $p+R(\t,t)=0$.

\ms

The following example shows that sometimes, $\psi$ can also be
continuous.

\ms

\bf Example 3.4. \rm Let $n=m=d=1$, $U=[-1,1]$, and
$$b(t,x,u)=u,\q\si(t,x,u)=\si(t,x),\q g(\t,t,x,u)={1\over2}R(\t,t)u^2,$$
with $\si(\cd\,,\cd)$ and $R(\cd\,,\cd)$ being continuous, and
$R(\t,t)>0$ for all $(\t,t)\in D[0,T]$. Then
$$\ba{ll}
\ns\ds\dbH(\t,t,x,u,p,P)=pu+{1\over2}\si(t,x)^2P+{1\over2}R(\t,t)u^2\\
\ns\ds\qq\qq\qq\q~={1\over2}\si(t,x)^2P+{1\over2}R(\t,t) \(u+{p\over
R(\t,t)}\)^2-{p^2\over2R(\t,t)}.\ea$$
Consequently,
$$H(\t,t,x,p,P)={1\over2}\si(t,x)^2P+{1\over2}R(\t,t)
\({|p|\land[R(\t,t)]\over R(\t,t,x)}-{|p|\over
R(\t,t)}\)^2-{p^2\over2R(\t,t)},$$
and
$$\psi(\t,t,x,p,P)=-\big[\sgn p\big]\({|p|\over R(\t,t)}\land1\),\qq\forall(\t,t,x,p,P)\in D[0,T]\times
\dbR^3.$$
Clearly, both $H$ and $\psi$ are continuous. Also, we see that even
if all the coefficients are very smooth, we cannot guarantee that
$H$ and $\psi$ are as smooth as the coefficients, in
general.

\ms

Here is another example which shows that $H$ and $\psi$ could be as
smooth as the coefficients.

\ms

\bf Example 3.5. \rm Let $n=m=d=1$, $U=(-1,1)$, and
$$b(t,x,u)=u,\q\si(t,x,u)=\si(t,x),\q g(\t,t,x,u)=-R(\t,t,x)\ln(1-u^2),$$
with $R(\t,t,x)$ being positive-valued and bounded. Then
$$\dbH(\t,t,x,u,p,P)=pu+{1\over2}\si(t,x)^2P-R(\t,t,x)\ln(1-u^2).$$
A direct computation shows that
$$\ba{ll}
\ns\ds H(\t,t,x,p,P)=p\psi(\t,t,x,p,P)+{1\over2}\si(t,x)^2P-
R(\t,t,x)\ln\big[1-\psi(\t,t,x,p,P)^2\big]\\
\ns\ds={1\over2}\si(t,x)^2P-{p^2\over
R(\t,t,x)+\sqrt{R(\t,t,x)^2+p^2\,}}
+R(\t,t,x)\ln{R(\t,t,x)+\sqrt{R(\t,t,x)^2+p^2}\over2R(\t,t,x)},\ea$$
with
$$\psi(\t,t,x,p,P)={-p\over R(\t,t,x)+\sqrt{R(\t,t,x)^2+p^2}\,}.$$
Therefore, both $H$ and $\psi$ are as smooth as the coefficients.

\ms

From the above discussion, we see that the situation concerning the
map $\psi(\cd)$ is very complicated. For the simplicity of
presentation below, we adopt the following assumption.

\ms

{\bf(H3)} The map $\psi:D[0,T]\times\dbR^n\times\dbR^n\times\cS^n\to
U$ is well-defined and has needed regularity.

\ms

We will address more general situations concerning $\psi(\cd)$ in
our future publications.

\ms

Now, let us recall a standard verification theorem for Problem (C)
stated in Section 1, which will be used below. For our later
purposes, it suffices to consider Problem (C) with the discount rate
$\d=0$. We will denote
$$C^{1,2}([0,T]\times\dbR^n)=\Big\{v(\cd\,,\cd)\in C([0,T]\times\dbR^n)
\bigm|v_t(\cd\,,\cd),v_x(\cd\,,\cd),v_{xx}(\cd\,,\cd)\in
C([0,T]\times\dbR^n)\Big\}.$$
The proof of the following result can be found in
\cite{Fleming-Soner 2006}.

\ms

\bf Proposition 3.6. \sl Suppose $V^0(\cd\,,\cd)\in
C^{1,2}([0,T]\times\dbR^n)$ is a classical solution to the following
Hamilton-Jacobi-Bellman equation:
\bel{HJB0}\left\{\ba{ll}
\ns\ds V^0_t(t,x)+\inf_{u\in
U}\dbH^0\big(t,x,u,V^0_x(t,x),V^0_{xx}(t,x)\big)=0,\qq(t,x)\in[0,T]\times\dbR^n,\\
\ns\ds V^0(T,x)=h^0(x),\qq x\in\dbR^n,\ea\right.\ee
where
\bel{}\dbH^0(t,x,u,p,P)=\lan
b(t,x,u),p\ran+\tr\big[a(t,x,u)P\big]+g^0(t,x,u).\ee
Then
\bel{}V^0(t,x)\le
J^0(t,x;u(\cd)),\qq\forall(t,x)\in[0,T]\times\dbR^n,~
u(\cd)\in\cU[t,T].\ee
If $(t,x)\in[0,T)\times\dbR^n$ is given and $(\bar X(\cd),\bar
u(\cd))$ is a state-control pair with the initial pair $(t,x)$ such
that
\bel{}\bar u(s)\in\arg\min\dbH^0\big(s,\bar X(s),\cd\,,V^0_x(s,\bar
X(s)),V^0_{xx}(s,\bar X(s))\big),\q s\in[t,T],\ee
then
\bel{}V^0(t,x)=J^0(t,x;\bar u(\cd))\ee
and $(\bar X(\cd),\bar u(\cd))$ is an optimal pair of Problem
{\rm(C)} for the initial pair $(t,x)$.

\ms

\rm

We make some remarks on the above verification theorem. First of
all, to guarantee (\ref{HJB0}) to have a classical solution
$V^0(\cd\,,\cd)$, one may pose different conditions. A typical one
is the uniform ellipticity condition:
\bel{}a(t,x,u)\equiv{1\over2}\si(t,x,u)\si(t,x,u)^T\ge\d
I_n,\qq\qq\forall(t,x,u)\in[0,T]\times\dbR^n\times U,\ee
for some $\d>0$. This condition implies that $n\le d$ and
$\si(t,x,u)$ stays full rank for all $(t,x,u)$. Thus, it does not
include the case that $(x,u)\mapsto\si(t,x,u)$ is linear, which is
the case for LQ problems. On the other hand, for a standard LQ
problem with deterministic coefficients, when the Riccati equation
admits a solution $P(\cd)$, the function $V^0(t,x)=\lan P(t)x,x\ran$
is a classical solution to the corresponding HJB equation for which
the uniform ellipticity condition fails. The similar situation
happens for the classical Merton's portfolio problem. This
observation shows that there are quite different conditions under
which the corresponding HJB equation admits a classical solution. In
the following sections, from time to time, we will simply say that
the relevant HJB equation has a classical solution without getting
into detailed sufficient conditions for that. Likewise, suppose
there exists a map
$\psi^0:[0,T]\times\dbR^n\times\dbR^n\times\cS^n\to U$ such that
$$\dbH^0\big(s,x,\psi^0(s,x,p,P),p,P\big)=\inf_{u\in U}\dbH^0(s,x,u,p,P),\qq
\forall(s,x,p,P)\in[0,T]\times\dbR^n\times\dbR^n\times\cS^n.$$
Then the pair $(\bar X(\cd),\bar u(\cd))$ appears in Proposition 3.6
is just a state-control pair satisfying the following
\bel{}\left\{\ba{ll}
\ns\ds d\bar X(s)=b\big(s,\bar X(s),\psi^0(s,\bar X(s),V^0_x(s,\bar
X(s)),V^0_{xx}(s,\bar X(s))\big)ds\\
\ns\ds\qq\qq+\si\big(s,\bar X(s),\psi^0(s,\bar X(s),V^0_x(s,\bar
X(s)),V^0_{xx}(s,\bar X(s))\big)sW(s),\qq s\in[t,T],\\
\ns\ds\bar X(t)=x,\ea\right.\ee
and we do not have to have the Lipschitz continuity of the map
$$x\mapsto\(b\big(s,x,\psi^0(s,x,V^0_x(s,x),V^0_{xx}(s,x))\big),
\si\big(s,x,\psi^0(s,x,V^0_x(s,x),V^0_{xx}(s,x))\big)\).$$
In the following section, from time to time, we will just say some
process is a solution to the relevant closed-loop system without
mentioning if the drift and diffusion are Lipschitz continuous, etc.

\rm

\section{Time-Consistent Equilibria via Multi-Person Differential Games}

In this section, we are going to search time-consistent solution to
Problem (N). Inspired by \cite{Yong 2011, Yong 2012}, we will take
an approach of multi-person differential games.

\ms

To begin with, let us first introduce some necessary notions. Let
$\cP[0,T]$ be the set of all partitions $\Pi=\{t_k\bigm|0\le k\le
N\}$ of $[0,T]$ with $0=t_0<t_1<t_2<\cds<t_{N-1}<t_N=T$, and with
the mesh size $\|\Pi\|$ defined by the following:
$$\|\Pi\|=\max_{1\le k\le N}(t_k-t_{k-1}).$$
We introduce the following important definition.

\ms

\bf Definition 4.1. \rm A continuous map $\Psi:[0,T]\times\dbR^n\to
U$ is called a {\it time-consistent equilibrium strategy} of Problem
(N) if for any $(t,x)\in[0,T)\times\dbR^n$,
\bel{}\left\{\ba{ll}
\ns\ds d\bar X(s)=b\big(s,\bar X(s),\Psi(s,\bar
X(s))\big)ds+\si\big(s,\bar X(s),\Psi(s,\bar X(s))\big)dW(s),\qq
s\in[t,T],\\
\ns\ds\bar X(t)=x,\ea\right.\ee
admits a unique solution $\bar X(\cd)\equiv\bar
X(\cd\,;t,x,\Psi(\cd))$ such that for any $t\in[0,T)$, and $\e>0$
with $t+\e\le T$,
\bel{}\ba{ll}
\ns\ds J\(t,\bar
X\big(t;0,x,\Psi(\cd)\big);\Psi(\cd)\big|_{[t,T]}\)\\
\ns\ds\le J\(t,\bar
X\big(t;0,x,\Psi(\cd)\big);u(\cd)\oplus\Psi(\cd)\big|_{[t+\e,T]}\)+R(\e),\qq\forall
u(\cd)\in\cU[t,t+\e],\ea\ee
where $R(\e)$ represents a generic remainder term satisfying
$R(\e)\to0$ as $\e\to0$, and
\bel{uPsi0}\big[u(\cd)\oplus\Psi(\cd)\big|_{[t+\e,T]}\big](s)=\left\{\ba{ll}
\ns\ds u(s),\qq\qq t\in[t,t+\e),\\
\ns\ds\Psi\big(s,X^\e(s)\big),\qq s\in[t+\e,T],\ea\right.\ee
with
\bel{4.4}\left\{\ba{ll}
\ns\ds
dX^\e(s)=b\big(s,X^\e(s),u(s)\big)ds+\si\big(s,X^\e(s),u(s)\big)dW(s),\qq
s\in[t,t+\e),\\
\ns\ds dX^\e(s)=b\big(s,X^\e(s),\Psi(s,X^\e(s))\big)ds+\si\big(s,
X^\e(s),\Psi(s,X^\e(s))\big)dW(s),\qq
s\in[t+\e,T],\\
\ns\ds X^\e(t)=\bar X\big(t;0,x,\Psi(\cd)\big).\ea\right.\ee
In this case, we call $\bar X(\cd)\equiv\bar X(\cd\,;0,x,\Psi(\cd))$
a {\it time-consistent equilibrium state process}, call the
corresponding $\bar u(\cd)\equiv\Psi\big(\cd\,;\bar X(\cd)\big)$ a
{\it time-consistent equilibrium control} for the initial state $x$,
and call $(\bar X(\cd),\bar u(\cd))$ a {\it time-consistent
equilibrium pair}. Further, function $V:[0,T]\times\dbR^n\to\dbR$ is
called the {\it equilibrium value function} of Problem (N) if
\bel{J=V}V\big(t,\bar X(t;0,x,\Psi(\cd))\big)=J\big(t,\bar
X(t;0,x,\Psi(\cd));\Psi(\cd)\big|_{[t,T]}\big),\qq t\in[0,T].\ee

\ms

Note that our definition of time-consistent equilibrium strategy is
of closed-loop type. Thus, it is comparable with that given in
\cite{Ekeland-Lazrak 2010, Bjork-Murgoci 2008,Bjork-Murgoci-Zhou},
and it is different from that in \cite{Hu-Jin-Zhou 2011}.

\ms

In what follows, the following definition which is equivalent to the
above is more convenient to use.

\ms

\bf Definition 4.1$'$. \rm A continuous map
$\Psi:[0,T]\times\dbR^n\to U$ is called a {\it time-consistent
equilibrium strategy} of Problem (N) if for any
$(t,x)\in[0,T)\times\dbR^n$,
\bel{}\left\{\ba{ll}
\ns\ds d\bar X(s)=b\big(s,\bar X(s),\Psi(s,\bar
X(s))\big)ds+\si\big(s,\bar X(s),\Psi(s,\bar X(s))\big)dW(s),\qq
s\in[t,T],\\
\ns\ds\bar X(t)=x,\ea\right.\ee
admits a unique solution $\bar X(\cd)\equiv\bar
X(\cd\,;t,x,\Psi(\cd))$ such that for any $\Pi\equiv\{t_k\bigm|1\le
k\le N\}\in\cP[0,T]$,
\bel{}\ba{ll}
\ns\ds J\(t_{k-1},\bar
X\big(t_{k-1};0,x,\Psi(\cd)\big);\Psi(\cd)\big|_{[t_{k-1},T]}\)\\
\ns\ds\le J\(t_{k-1},\bar
X\big(t_{k-1};0,x,\Psi(\cd)\big);u^k(\cd)\oplus\Psi(\cd)\big|_{[t_k,T]}\)+R\big(\|\Pi\|\big),\\
\ns\ds\qq\qq\qq\qq\qq\qq\qq\qq\qq\qq\qq\forall
u^k(\cd)\in\cU[t_{k-1},t_k],\ea\ee
where $R(r)$ represents a generic remainder term satisfying
$R(r)\to0$ as $r\to0$, and
\bel{uPsi}\big[u^k(\cd)\oplus\Psi(\cd)\big|_{[t_k,T]}\big](t)=\left\{\ba{ll}
\ns\ds u^k(t),\qq\qq t\in[t_{k-1},t_k),\\
\ns\ds\Psi\big(t,X^k(t)\big),\qq t\in[t_k,T],\ea\right.\ee
with
\bel{4.4}\left\{\ba{ll}
\ns\ds
dX^k(s)=b\big(s,X^k(s),u^k(s)\big)ds+\si\big(s,X^k(s),u^k(s)\big)dW(s),\qq
s\in[t_{k-1},t_k),\\
\ns\ds dX^k(s)=b\big(s,X^k(s),\Psi(s,X^k(s))\big)ds+\si\big(s,
X^k(s),\Psi(s,X^k(s))\big)dW(s),\qq
s\in[t_k,T],\\
\ns\ds X^k(t_{k-1})=\bar
X\big(t_{k-1};0,x,\Psi(\cd)\big).\ea\right.\ee

\subsection{Multi-Person Differential Games.}

We now consider an $N$-person differential game, called Problem
(G$^\Pi$), as briefly described in Section 1. Throughout this
section, we assume that (H1)--(H3) hold. Let us start with Player
$N$ who controls the system on $[t_{N-1}t_N)$. More precisely, for
each $(t,x)\in[t_{N-1},t_N]\times\dbR^n$, consider the following
controlled SDE:
\bel{1.1N}\left\{\ba{ll}
\ns\ds dX^N(s)=b\big(s,X^N(s),u^N(s)\big)ds+\si\big(s,X^N(s),u^N(s)\big)dW(s),\qq s\in[t,t_N],\\
\ns\ds X^N(t)=x,\ea\right.\ee
with cost functional
\bel{1.2N}J^N(t,x;u^N(\cd))=\dbE_t\[\int_t^{t_N}g\big(t_{N-1},s,X^N(s),
u^N(s)\big)ds+h\big(t_{N-1},X^N(t_N)\big)\].\ee
Note that
\bel{J^N=J}J^N(t_{N-1},x;u^N(\cd))=J(t_{N-1},x;u^N(\cd)),\qq(x,u^N(\cd))\in\dbR^n\times
\cU[t_{N-1},t_N].\ee
We pose
the following problem for Player $N$:

\ms

\bf Problem (C$_N$). \rm For any
$(t,x)\in[t_{N-1},t_N)\times\dbR^n$, find a $\bar u^N(\cd)\equiv\bar
u^N(\cd\,;t,x)\in\cU[t,t_N]$ such that
\bel{}J^N(t,x;\bar
u^N(\cd))=\inf_{u^N(\cd)\in\cU[t,t_N]}J^N(t,x;u^N(\cd))\equiv
V^\Pi(t,x).\ee

The above defines the {\it value function} $V^\Pi(\cd\,,\cd)$ on
$[t_{N-1},t_N]\times\dbR^n$, and in the case $\bar u^N(\cd)$ exists,
by (\ref{J^N=J}), we have
\bel{V=J^N}J(t_{N-1},x;\bar u^N(\cd))=V^\Pi(t_{N-1},x),\qq\forall
x\in\dbR^n.\ee
Under proper conditions, $V^\Pi(\cd\,,\cd)$ is the classical
solution to the following HJB equation:
\bel{HJB-N}\left\{\ba{ll}
\ns\ds V^\Pi_t(t,x)+\inf_{u\in
U}\dbH\big(t_{N-1},t,x,u,V^\Pi_x(t,x),
V^\Pi_{xx}(t,x)\big)=0,\qq(t,x)\in[t_{N-1},t_N]\times\dbR^n,\\
\ns\ds V^\Pi(t_N,x)=h(t_{N-1},x),\qq x\in\dbR^n.\ea\right.\ee
By the definition of
$\psi:D[0,T]\times\dbR^n\times\dbR^n\times\cS^n\to U$ (see
(\ref{psi1})--(\ref{psi2})), we may also write (\ref{HJB-N}) as
follows
\bel{HJB-N*}\left\{\ba{ll}
\ns\ds V^\Pi_t(t,x)+\dbH\big(t_{N-1},t,x,\psi(t_{N-1},t,x,V^\Pi_x(t,x),V^\Pi_{xx}(t,x)),V^\Pi_x(t,x),V^\Pi_{xx}(t,x)\big)=0,\\
\ns\ds\qq\qq\qq\qq\qq\qq\qq(t,x)\in[t_{N-1},t_N]\times\dbR^n,\\
\ns\ds V^\Pi(t_N,x)=h(t_{N-1},x),\qq x\in\dbR^n.\ea\right.\ee
Clearly, $V^\Pi(\cd\,,\cd)$, well-defined on
$[t_{N-1},t_N]\times\dbR^n$, depends on $t_{N-1}$ and $t_N$. With
such a solution $V^\Pi(\cd\,,\cd)$ of (\ref{HJB-N}) (or
(\ref{HJB-N*})), let us assume that the following closed-loop system
admits a unique solution $\bar X^N(\cd)\equiv\bar
X^N(\cd\,;t_{N-1},x)$: (suppressing the dependence of $\bar
X^N(\cd)$ on $t_N$ through $V^\Pi(\cd\,,\cd)$)
\bel{closed-loop-N}\left\{\ba{ll}
\ns\ds d\bar X^N(s)=b\big(s,\bar X^N(s),\psi(t_{N-1},s,\bar
X^N(s),V^\Pi_x(s,\bar X^N(s)),V^\Pi_{xx}(s,\bar X^N(s)))\big)ds\\
\ns\ds\qq\qq\q+\si\big(s,\bar X^N(s),\psi(t_{N-1},s,\bar
X^N(s),V^\Pi_x(s,\bar X^N(s)),V^\Pi_{xx}(s,\bar X^N(s)))\big)dW(s),\\
\ns\ds\qq\qq\qq\qq\qq\qq\qq\qq\qq\qq s\in[t_{N-1},t_N],\\
\ns\ds\bar X^N(t_{N-1})=x.\ea\right.\ee
Then under (H3) and Proposition 3.6, an optimal control $\bar
u^N(\cd)$ of Problem (C$_N$) for the initial pair $(t_{N-1},x)$
admits the following feedback representation:
\bel{bar u-N}\ba{ll}
\ns\ds\bar u^N(s)\equiv\bar u^N(s;t_{N-1},x)=\psi\big(t_{N-1},s,\bar
X^N(s),V^\Pi_x(s,\bar
X^N(s)),V^\Pi_{xx}(s,\bar X^N(s))\big)\\
\ns\ds\qq~\equiv\psi\big(t_{N-1},s,\bar
X^N(s;t_{N-1},x),V^\Pi_x(s,\bar
X^N(s;t_{N-1},x)),V^\Pi_{xx}(s,\bar X^N(s;t_{N-1},x))\big)\\
\ns\ds\qq\qq\qq\qq\qq\qq\qq\qq\qq\qq\qq\qq s\in[t_{N-1},t_N],\ea\ee
and $\bar X^N(\cd)\equiv\bar X^N(\cd\,;t_{N-1},x)$ is the
corresponding optimal state process.

\ms

Next, we consider an optimal control problem for Player $(N-1)$ on
$[t_{N-2},t_{N-1})$. Naturally, for any initial pair
$(t,x)\in[t_{N-2},t_{N-1}]\times\dbR^n$, the state equation should
be
\bel{state-(N-1)}\left\{\3n\ba{ll}
\ns\ds
dX^{N-1}(s)=b\big(s,X^{N-1}(s),u^{N-1}(s)\big)ds+\si\big(s,X^{N-1}(s),u^{N-1}(s)\big)dW(s),\q
s\in[t,t_{N-1}),\\
\ns\ds X^{N-1}(t)=x.\ea\right.\ee
To determine the suitable cost functional, we note that Player
$(N-1)$ can only control the system on $[t_{N-2},t_{N-1})$ and
Player $N$ will take over at $t_{N-1}$ to control the system
thereafter. Moreover, Player $(N-1)$ knows that Player $N$ will play
optimally based on the initial pair $(t_{N-1},X^{N-1}(t_{N-1}))$ of
Player $N$, which is the {\it terminal pair} of Player $(N-1)$.
Hence, the {\it sophisticated cost functional} of Player $(N-1)$
should be
\bel{Natural(N-1)}\ba{ll}
\ns\ds
J^{N-1}(t,x;u^{N-1}(\cd))=\dbE_t\[\int_t^{t_{N-1}}g(t_{N-2},s,X^{N-1}(s),u^{N-1}(s))
ds\\
\ns\ds\qq\qq\qq\qq\qq+\2n\int_{t_{N-1}}^{t_N}\3n\2n
g\big(t_{N-2},s,\bar X^N(s;t_{N-1},X^{N-1}(t_{N-1})),\bar
u^N(s;t_{N-1},X^{N-1}(t_{N-1}))\big)ds\\
\ns\ds\qq\qq\qq\qq\qq+h\big(t_{N-2},\bar
X^N(t_N;t_{N-1},X^{N-1}(t_{N-1}))\big)\].\ea\ee
Note that although Player $(N-1)$ knows that Player $N$ will control
the system on $[t_{N-1},t_N]$, he/she still ``discounts'' the future
costs in his/her own way (see $t_{N-2}$ appearing in the running
cost on $[t_{N-1},t_N]$ and in the terminal cost at $t_N$). Now, if
we denote
$$h^{N-1}(x)=\dbE_{t_{N-1}}\[\int_{t_{N-1}}^{t_N}g\big(t_{N-2},s,\bar X^N(s;t_{N-1},x),
\bar u^N(s;t_{N-1},x)\big)ds+h\big(t_{N-2},\bar
X^N(t_N;t_{N-1},x)\big)\],$$
then the cost functional (\ref{Natural(N-1)}) can be written as
\bel{S(N-1)}
J^{N-1}(t,x;u^{N-1}(\cd))=\dbE_t\[\int_t^{t_{N-1}}g\big(t_{N-2},s,X^{N-1}(s),u^{N-1}(s))
ds+h^{N-1}(X^{N-1}(t_{N-1})\big)\].\ee
We see that the optimal control problem associated with the state
equation (\ref{state-(N-1)}) and the cost functional (\ref{S(N-1)})
looks like a standard one. But, the map $x\mapsto h^{N-1}(x)$ seems
to be a little too implicit, which is difficult for us to pass to
the limits later on. We now would like to make it more explicit in
some sense. Inspired by the idea of Four Step Scheme introduced in
\cite{MPY 1994,MY 1999} for FBSDEs with deterministic coefficients,
we proceed as follows. For the optimal state process $\bar
X^N(\cd)\equiv\bar X^N(\cd\,;t_{N-1},x)$ determined by
(\ref{closed-loop-N}) on $[t_{N-1},t_N]$, we introduce the following
BSDE:
\bel{BSDE-N}\left\{\ba{ll}
\ns\ds dY^N(s)=-g\big(t_{N-2},s,\bar X^N(s),\psi(t_{N-1},s,\bar
X^N(s),V^\Pi_x(s,\bar X^N(s)),V^\Pi_{xx}(s,\bar X^N(s)))\big)ds\\
\ns\ds\qq\qq\qq\qq\qq\qq+Z^N(s)dW(s),\qq\qq\qq s\in[t_{N-1},t_N],\\
\ns\ds Y^N(t_N)=h(t_{N-2},\bar X^N(t_N)),\ea\right.\ee
which is equivalent to the following:
\bel{BSDE-Ne}\left\{\ba{ll}
\ns\ds dY^N(s)=-g\big(t_{N-2},s,\bar X^N(s),\bar
u^N(s)\big)ds+Z^N(s)dW(s),\qq s\in[t_{N-1},t_N],\\
\ns\ds Y^N(t_N)=h(t_{N-2},\bar X^N(t_N)),\ea\right.\ee
Note that $t_{N-2}$ appears in the drift of BSDE and in the terminal
condition. This BSDE admits a unique adapted solution
$(Y^N(\cd),Z^N(\cd))\equiv\big(Y^N(\cd\,;x),Z^N(\cd\,;x)\big)$
(\cite{MY 1999,Yong-Zhou 1999}), uniquely depending on $x\in\dbR^n$.
Further, one has
$$\ba{ll}
\ns\ds
Y^N(t_{N-1})=\dbE_{t_{N-1}}\[\int_{t_{N-1}}^{t_N}g\big(t_{N-2},s,
\bar X^N(s),\bar u^N(s)\big)ds+h(t_{N-2},\bar
X^N(t_N))\]=h^{N-1}(x).\ea$$
It is seen that (\ref{closed-loop-N}) and (\ref{BSDE-Ne}) form an
FBSDE. By \cite{MPY 1994} (see also \cite{MY 1999,Yong-Zhou 1999}),
we have the following representation for $Y^N(\cd)$
\bel{rep-N}Y^N(s)=\Th^N(s,\bar X^N(s)),\qq s\in[t_{N-1},t_N],\ee
as long as $\Th^N(\cd\,,\cd)$ is a classical solution to the
following PDE:
\bel{Parabolic-N}\left\{\3n\ba{ll}
\ns\ds\Th^N_s(s,x)+\lan\Th^N_x(s,x),
b\big(s,x,\psi(t_{N-1},s,x,V^\Pi_x(s,x),
V^\Pi_{xx}(s,x))\big)\ran\\
\ns\ds\qq+\tr\big[a\big(s,x,\psi(t_{N-1},s,x,V^\Pi_x(s,x),
V^N_{xx}(s,x))\big)\Th^N_{xx}(s,x)\big]\\
\ns\ds\qq+g\big(t_{N-2},s,x,\psi(t_{N-1},s,x,V^\Pi_x(s,x),V^\Pi_{xx}(s,x))
\big)=0,\qq(s,x)\in[t_{N-1},t_N]\times\dbR^n,\\
\ns\ds\Th^N(t_N,x)=h(t_{N-2},x),\qq x\in\dbR^n,\ea\right.\ee
or equivalently,
\bel{parabolic-N*}\left\{\3n\ba{ll}
\ns\ds\Th^N_s(s,x)+\dbH\big(t_{N-2},s,x,\psi\big(t_{N-1},s,x,V_x^\Pi(s,x),
V^\Pi_{xx}(s,x)\big),\Th^N_x(s,x),
\Th^N_{xx}(s,x)\big)=0,\\
\ns\ds\qq\qq\qq\qq\qq\qq\qq\qq(s,x)\in[t_{N-1},t_N]\times\dbR^n,\\
\ns\ds\Th^N(t_N,x)=h(t_{N-2},x),\qq x\in\dbR^n.\ea\right.\ee
Note that $\Th^N(\cd\,,\cd)$ depends on $(t_{N-2},t_{N-1},t_N)$. We
point out that in general,
\bel{Th-N ne V-N}\ba{ll}
\ns\ds\Th^N(t_{N-1},x)\1n=\1n
Y^N(t_{N-1})\1n=\1n\dbE_{t_{N-1}}\1n\[\int_{t_{N-1}}^{t_N}\3n
g\big(t_{N-2},s,\bar
X^N(s),\bar u^N(s)\big)ds\1n+\1n h(t_{N-2},\bar X^N(T))\]\1n=\1n h^{N-1}(x)\1n\\
\ns\ds\ne\dbE_{t_{N-1}}\[\int_{t_{N-1}}^{t_N}g\big(t_{N-1},s,\bar
X^N(s),\bar u^N(s)\big)ds+h(t_{N-1},\bar
X^N(T))\]=V^\Pi(t_{N-1},x).\ea\ee
With the above representation $\Th^N(\cd\,,\cd)$ of $Y^N(\cd)$, we
can rewrite the cost functional (\ref{S(N-1)}) as follows:
\bel{cost-(N-1)}
J^{N-1}(t,x;u^{N-1}(\cd))\1n=\1n\dbE_t\[\int_t^{t_{N-1}}\3n\3n
g(t_{N-2},s,X^{N-1}(s),u^{N-1}(s))ds
\1n+\1n\Th^N(t_{N-1},X^{N-1}(t_{N-1}))\].\ee
 We now
pose the following problem:

\ms

\bf Problem (C$_{N-1}$). \rm For any
$(t,x)\in[t_{N-2},t_{N-1})\times\dbR^n$, find a $\bar
u^{N-1}(\cd)\equiv\bar u^{N-1}(\cd\,;t,x)\in\cU[t_{N-2},t_{N-1}]$
such that
\bel{}J^{N-1}(t,x;\bar
u^{N-1}(\cd))=\inf_{u^{N-1}(\cd)\in\cU[t,t_{N-1}]}J^{N-1}(t,x;u^{N-1}(\cd))\equiv
V^\Pi(t,x).\ee

\ms

The above defines the value function $V^\Pi(\cd\,,\cd)$ on
$[t_{N-2},t_{N-1})\times\dbR^n$. Under proper conditions,
$V^\Pi(\cd\,,\cd)$ is the classical solution to the following HJB
equation:
\bel{HJB-(N-1)}\left\{\ba{ll}
\ns\ds V^\Pi_t(t,x)+\inf_{u\in
U}\dbH\big(t_{N-2},t,x,u,V^\Pi_x(t,x),V^\Pi_{xx}(t,x)\big)=0,\qq(t,x)
\in[t_{N-2},t_{N-1})\times\dbR^n,\\
\ns\ds V^\Pi(t_{N-1}-0,x)=\Th^N(t_{N-1},x),\qq
x\in\dbR^n.\ea\right.\ee
By the definition of the map $\psi(\cd)$ again (see
(\ref{psi1})--(\ref{psi2})), we may also write the above as
\bel{HJB-(N-1)*}\left\{\ba{ll}
\ns\ds
V^\Pi_t(t,x)+\dbH\big(t_{N-2},t,x,\psi\big(t_{N-2},t,x,V^\Pi_x(t,x),
V^\Pi_{xx}(t,x)\big),V^\Pi_x(t,x),V^\Pi_{xx}(t,x)\big)=0,\\
\ns\ds\qq\qq\qq\qq\qq\qq\qq\qq\qq\qq\qq(t,x)\in[t_{N-2},t_{N-1})\times\dbR^n,\\
\ns\ds V^\Pi(t_{N-1}-0,x)=\Th^N(t_{N-1},x),\qq
x\in\dbR^n.\ea\right.\ee
From (\ref{Th-N ne V-N}), we see that in general,
$$V^\Pi(t_{N-1}-0,x)=\Th^N(t_{N-1},x)\ne V^\Pi(t_{N-1},x).$$
Thus, $V^\Pi(\cd\,,\cd)$, which is now defined on
$[t_{N-2},t_N]\times\dbR^n$, may have a jump at
$\{t_{N-1}\}\times\dbR^n$. For any $x\in\dbR^n$, suppose the
following admits a unique solution $\bar X^{N-1}(\cd)$:
\bel{closed-loop-(N-1)}\left\{\3n\ba{ll}
\ns\ds d\bar X^{N-1}(s)=b\(s,\bar X^{N-1}(s),\psi\big(t_{N-2},s,\bar
X^{N-1}(s),V^\Pi_x(s,\bar X^{N-1}(s)),V^\Pi_{xx}(s,\bar X^{N-1}(s))\big)\)ds\\
\ns\ds\qq\qq\qq+\si\1n\(s,\1n\bar
X^{N-1}\1n(s),\1n\psi\big(t_{N-2},s,\1n\bar X^{N-1}\1n(s),\1n
V^\Pi_x(s,\1n\bar X^{N-1}\1n(s)),
\1n V^\Pi_{xx}(s,\1n\bar X^{N-1}\1n(s))\big)\)dW(s),\\
\ns\ds\qq\qq\qq\qq\qq\qq\qq\qq\qq\qq s\in[t_{N-2},t_{N-1}),\\
\ns\ds\bar X^{N-1}(t_{N-2})=x.\ea\right.\ee
Then we define
$$\ba{ll}
\ns\ds\bar u^{N-1}(s;t_{N-2},x)=\psi\big(t_{N-2},s,\bar
X^{N-1}(s),V^\Pi_x(s,\bar X^{N-1}(s)),
V^\Pi_{xx}(s,\bar X^{N-1}(s))\big),\\
\ns\ds\q\equiv\psi\big(t_{N-2},s,\bar
X^{N-1}(s;t_{N-2},x),V^\Pi_x(s,\bar X^{N-1}(s;t_{N-2},x)),
V^\Pi_{xx}(s,\bar X^{N-1}(s;t_{N-2},x))\big),\\
\ns\ds\qq\qq\qq\qq\qq\qq\qq\qq\qq\qq\qq\qq
s\in[t_{N-2},t_{N-1}),\ea$$
which, again by Proposition 3.6, is an optimal control of Problem
(C$_{N-1}$) with the initial pair $(t_{N-2},x)$. Now, for the
optimal pair
$$\big(\bar X^{N-1}(\cd),\bar u^{N-1}(\cd)\big)
=\big(\bar X^{N-1}(\cd\,;t_{N-2},x),\bar
u^{N-1}(\cd\,;t_{N-2},x)\big)$$
of Problem (C$_{N-1}$) (on $[t_{N-2},t_{N-1}]$), we make a natural
extension on $[t_{N-1},t_N]$ as follows:
\bel{}\left\{\ba{ll}
\ns\ds\bar X^{N-1}(s)=\bar X^N(s;t_{N-1},\bar X^{N-1}(t_{N-1})),\\
\ns\ds\bar u^{N-1}(s)=\bar u^N(s;t_{N-1},\bar
X^{N-1}(t_{N-1})),\ea\right.\qq s\in[t_{N-1},t_N].\ee
Thus, the extended $(\bar X^{N-1}(\cd),\bar u^{N-1}(\cd))$ satisfies
\bel{3.44}\left\{\3n\ba{ll}
\ns\ds d\bar X^{N-1}(s)=b\(s,\bar X^{N-1}(s),\psi\big(t_{N-2},s,\bar
X^{N-1}(s),V^\Pi_x(s,\bar X^{N-1}(s)),V^\Pi_{xx}(s,\bar X^{N-1}(s))\big)\)ds\\
\ns\ds\qq\qq\q+\si\(s,\bar X^{N-1}(s),\psi\big(t_{N-2},s,\bar
X^{N-1}(s),V^\Pi_x(s,\bar X^{N-1}(s)),V^\Pi_{xx}(s,\bar X^{N-1}(s))\big)\)dW(s),\\
\ns\ds\qq\qq\qq\qq\qq\qq\qq\qq\qq\qq s\in[t_{N-2},t_{N-1}),\\
\ns\ds d\bar X^{N-1}(s)=b\(s,\bar X^{N-1}(s),\psi\big(t_{N-1},s,\bar
X^{N-1}(s),V^\Pi_x(s,\bar X^{N-1}(s)),V^\Pi_{xx}(s,\bar X^{N-1}(s))\big)\)ds\\
\ns\ds\qq\qq\q+\si\(s,\bar X^{N-1}(s),\psi\big(t_{N-1},s,\bar
X^{N-1}(s),V^\Pi_x(s,\bar X^{N-1}(s)),V^\Pi_{xx}(s,\bar X^{N-1}(s))\big)\)dW(s),\\
\ns\ds\qq\qq\qq\qq\qq\qq\qq\qq\qq\qq s\in[t_{N-1},t_N],\\
\ns\ds\bar X^{N-1}(t_{N-2})=x.\ea\right.\ee
We refer to such a pair $\big(\bar X^{N-1}(\cd),\bar
u^{N-1}(\cd)\big)$ as a {\it sophisticated equilibrium pair} on
$[t_{N-2},t_N]$. Let
\bel{ell}\ell^\Pi(s)=\sum_{k=1}^Nt_{k-1}I_{[t_{k-1},t_k)}(s),\qq
s\in[0,T].\ee
It is easy to see that
\bel{}0\le s-\ell^\Pi(s)\le\|\Pi\|,\qq s\in[0,T].\ee
Then (\ref{3.44}) can be written compactly as
\bel{3.45}\left\{\3n\ba{ll}
\ns\ds d\bar X^{N-1}(s)=b\(s,\bar
X^{N-1}(s),\psi\big(\ell^\Pi(s),s,\bar
X^{N-1}(s),V^\Pi_x(s,\bar X^{N-1}(s)),V^\Pi_{xx}(s,\bar X^{N-1}(s))\big)\)ds\\
\ns\ds\qq\qq\qq+\si\(s,\1n\bar
X^{N-1}(s),\psi\big(\ell^\Pi(s),s,\1n\bar
X^{N-1}(s),V^\Pi_x(s,\1n\bar X^{N-1}(s)),V^\Pi_{xx}(s,\1n\bar X^{N-1}(s))\big)\)dW(s),\\
\ns\ds\qq\qq\qq\qq\qq\qq\qq\qq\qq\qq\qq\qq\qq s\in[t_{N-2},t_N],\\
\ns\ds\bar X^{N-1}(t_{N-2})=x.\ea\right.\ee
Also, one has
\bel{J^{N-1}=V}\ba{ll}
\ns\ds J(t_{N-2},x;\bar
u^{N-1}(\cd))=\dbE_{t_{N-2}}\[\int_{t_{N-1}}^{t_{N-1}}g(t_{N-2},s,\bar
X^{N-1}(s),\bar u^{N-1}(s))ds\\
\ns\ds\qq\qq\qq\qq\qq\qq+\int_{t_{N-1}}^{t_N}g(t_{N-2},s,\bar
X^{N-1}(s),\bar u^{N-1}(s))ds+h(t_{N-2},\bar X^{N-1}(T))\]\\
\ns\ds\qq\qq\qq\qq~=J^{N-1}(t_{N-2},x;\bar
u^{N-1}(\cd))=V^\Pi(t_{N-2},x),\qq x\in\dbR^n.\ea\ee
We point out that in general it may happen that,
\bel{J>J}J\big(t_{N-2},x;\bar
u^{N-1}(\cd)\big)>\inf_{u(\cd)\in\cU[t_{N-2},t_N]}J\big(t_{N-2},x;u(\cd)\big),\ee
which means that the the sophisticated equilibrium pair might not be
an optimal pair (for the given initial pair).

\ms

Similar to the above, in order to state an optimal control problem
for Player $(N-2)$ on $[t_{N-3},t_{N-2}]$, we introduce the
following BSDE on $[t_{N-2},t_N]$:
\bel{BSDE-(N-1)}\left\{\3n\ba{ll}
\ns\ds dY^{N-1}(s)=-g\big(t_{N-3},s,\bar X^{N-1}(s),\bar u^{N-1}(s)\big)ds+Z^{N-1}(s)dW(s),\q s\in[t_{N-2},t_N],\\
\ns\ds Y^{N-1}(t_N)=h(t_{N-3},\bar X^{N-1}(t_N)),\ea\right.\ee
where $\big(\bar X^{N-1}(\cd),\bar u^{N-1}(\cd)\big)$ is the
sophisticated equilibrium pair determined by (\ref{3.45}) on
$[t_{N-2},t_N]$, uniquely depending on $x\in\dbR^n$. Let
$\big(Y^{N-1}(\cd),Z^{N-1}(\cd)\big)\equiv\big(Y^{N-1}(\cd\,;x),Z^{N-1}(\cd\,;x)\big)$
be the adapted solution of this BSDE. Then (\ref{3.45}) and
(\ref{BSDE-(N-1)}) form an FBSDE. Similar to the above, we have
\bel{rep-(N-1)}Y^{N-1}(s)=\Th^{N-1}(s,\bar X^{N-1}(s)),\qq
s\in[t_{N-2},t_N],\ee
as long as $\Th^{N-1}(\cd\,,\cd)$ is the solution to the following
PDE:
\bel{}\left\{\3n\ba{ll}
\ns\ds\Th^{N-1}_s(s,x)\1n+\1n
\dbH\big(t_{N-3},s,x,\psi\big(\ell^\Pi(s),s,x,V^\Pi_x(s,x),V^\Pi_{xx}(s,x)\big),
\Th^{N-1}_x(s,x),\Th^{N-1}_{xx}(s,x)\big)\1n=\1n0,\\
\ns\ds\qq\qq\qq\qq\qq\qq(s,x)\in[t_{N-2},t_N]\times\dbR^n,\\
\ns\ds\Th^{N-1}(t_N,x)=h(t_{N-3},x),\qq x\in\dbR^n.\ea\right.\ee
Having the above preparation, we now consider, for any
$(t,x)\in[t_{N-3},t_{N-2})\times\dbR^n$, the state equation
\bel{}\left\{\ba{ll}
\ns\ds
dX^{N-2}(s)=b(s,X^{N-2}(s),u^{N-2}(s))ds+\si(s,X^{N-2}(s),u^{N-2}(s))dW(s),\qq
s\in[t,t_{N-2}],\\
\ns\ds X^{N-2}(t)=x,\ea\right.\ee
and the (sophisticated) cost functional
\bel{}J^{N-2}(t,x;u^{N-2}(\cd))=\dbE_t\[\int_t^{t_{N-2}}g(t_{N-3},s,X^{N-2}(s),u^{N-2}(s))ds
+\Th^{N-1}(t_{N-2},X^{N-2}(t_{N-2}))\].\ee
We pose the following problem for Player $(N-2)$:

\ms

\bf Problem (C$_{N-2}$). \rm For any
$(t,x)\in[t_{N-3},t_{N-2})\times\dbR^n$, find a $\bar
u^{N-2}(\cd)\equiv\bar u^{N-2}(\cd\,;t,x)\in\cU[t_{N-3},t_{N-2}]$
such that
\bel{}J^{N-2}(t,x;\bar
u^{N-2}(\cd))=\inf_{u^{N-2}(\cd)\in\cU[t_{N-3},t_{N-2}]}J^{N-2}(t,x;u^{N-2}(\cd))=V^\Pi(t,x).\ee

\ms

Together with the previous definition, we see that
$V^\Pi(\cd\,,\cd)$ is now well-defined on
$[t_{N-3},t_N]\times\dbR^n$. Under proper conditions,
$V^\Pi(\cd\,,\cd)$ is a classical solution to the following HJB
equation:
\bel{HJB-(N-2)}\left\{\ba{ll}
\ns\ds V^\Pi_t(t,x)+\inf_{u\in
U}\dbH\big(\ell^\Pi(t),t,x,u,V^\Pi_x(t,x),V^\Pi_{xx}(t,x)\big)=0,\qq(t,x)\in[t_{N-3},
t_{N-2})\times\dbR^n,\\
\ns\ds V^\Pi(t_{N-2}-0,x)=\Th^{N-1}(t_{N-2},x),\qq
x\in\dbR^n.\ea\right.\ee
Further, by the definition of the map $\psi(\cd)$, we may also write
the above as
\bel{HJB-(N-2)*}\left\{\ba{ll}
\ns\ds
V^\Pi_t(t,x)+\dbH\big(\ell^\Pi(t),t,x,\psi\big(\ell^\Pi(t),t,x,V^\Pi_x(t,x),
V^\Pi_{xx}(t,x)\big),V^\Pi_x(t,x),
V^\Pi_{xx}(t,x)\big)=0,\\
\ns\ds\qq\qq\qq\qq\qq\qq\qq\qq\qq\qq\qq(t,x)\in[t_{N-3},t_{N-2})\times\dbR^n,\\
\ns\ds V^\Pi(t_{N-2}-0,x)=\Th^{N-1}(t_{N-2},x),\qq
x\in\dbR^n.\ea\right.\ee
Also, similar to (\ref{Th-N ne V-N}), we have that in general,
\bel{}V^\Pi(t_{N-2}-0,x)\ne V^\Pi(t_{N-2},x).\ee

The above procedure can be continued recursively. By induction, we
can construct sophisticated cost functional $J^k(t,x;u^k(\cd))$ for
Player $k$, and
\bel{}V^\Pi(t,x)=\inf_{u^k(\cd)\in\cU[t,t_k]}J^k(t,x;u^k(\cd)),
\qq(t,x)\in[t_{k-1},t_k)\times\dbR^n,~1\le
k\le N,\ee
with the value function $V^\Pi(\cd\,,\cd)$ satisfying the following
HJB equations on the time intervals associated with the partition
$\Pi$:
\bel{}\left\{\ba{ll}
\ns\ds
V^\Pi_t(t,x)+\dbH\big(\ell^\Pi(t),t,x,\psi\big(\ell^\Pi(t),t,x,V^\Pi_x(t,x),
V^\Pi_{xx}(t,x)\big),V^\Pi_x(t,x),V^\Pi_{xx}(t,x)\big)=0,\\
\ns\ds\qq\qq\qq\qq\qq\qq\qq\qq\qq\qq(t,x)\in[t_{N-1},t_N]\times\dbR^n,\\
\ns\ds V^\Pi(t_N,x)=h(t_{N-1},x),\qq x\in\dbR^n,\ea\right.\ee
and for $k=1,2,\cds,N-1$,
\bel{}\left\{\ba{ll}
\ns\ds
V^\Pi_t(t,x)+\dbH\big(\ell^\Pi(t),t,x,\psi\big(\ell^\Pi(t),t,x,V^\Pi_x(t,x),
V^\Pi_{xx}(t,x)\big),V^\Pi_x(t,x),
V^\Pi_{xx}(t,x)\big)=0,\\
\ns\ds\qq\qq\qq\qq\qq\qq\qq\qq\qq\qq(t,x)\in[t_{k-1},t_k)\times\dbR^n,\\
\ns\ds V^\Pi(t_k-0,x)=\Th^{k+1}(t_k,x),\qq x\in\dbR^n,\ea\right.\ee
where, for $k=1,2,\cds,N-1$, $\Th^{k+1}(\cd\,,\cd)$ is the solution
to the following (linear) PDE:
\bel{Th(k+1)}\left\{\ba{ll}
\ns\ds\Th^{k+1}_t(t,x)+\dbH\big(t_{k-1},t,x,\psi\big(\ell^\Pi(t),t,x,
V^\Pi_x(t,x),V^\Pi_{xx}(t,x)\big),\Th^{k+1}_x(t,x),\Th^{k+1}_{xx}(t,x)\big)=0,\\
\ns\ds\qq\qq\qq\qq\qq\qq\qq\qq\qq(t,x)\in[t_k,t_N]\times\dbR^n,\\
\ns\ds\Th^{k+1}(t_N,x)=h(t_{k-1},x),\qq\qq x\in\dbR^n.\ea\right.\ee
Now, we define
\bel{}\Psi^\Pi(t,x)=\psi\big(\ell^\Pi(t),t,x,V^\Pi_x(t,x),V^\Pi_{xx}(t,x)\big),\qq
(t,x)\in[0,T]\times\dbR^n.\ee
Then for any given $x\in\dbR^n$, let $\bar X^\Pi(\cd)$ be the
solution to the following closed-loop system:
\bel{}\left\{\ba{ll}
\ns\ds d\bar X^\Pi(s)=b\big(s,\bar X^\Pi(s),\Psi^\Pi(s,\bar
X^\Pi(s))\big)ds+\si\big(s,\bar X^\Pi(s),\Psi^\Pi(s,\bar
X^\Pi(s))\big)dW(s),\qq s\in[0,T],\\
\ns\ds\bar X^\Pi(0)=x,\ea\right.\ee
and denote
\bel{}\bar u^\Pi(s)=\Psi^\Pi\big(s,\bar X^\Pi(s)\big),\qq
s\in[0,T].\ee
According to  our construction, we have
\bel{Nash}\ba{ll}
\ns\ds J\big(t_{k-1},\bar
X^\Pi(t_{k-1});\Psi^\Pi(\cd)\big|_{[t_{k-1},T]}\big)=J\big(t_{k-1},\bar
X^\Pi(t_{k-1});\bar
u^\Pi(\cd)\big|_{[t_{k-1},t_k]}\big)=V^\Pi\big(t_{k-1},\bar
X^\Pi(t_{k-1})\big)\\
\ns\ds=J^k\big(t_{k-1},\bar X^\Pi(t_{k-1});\bar
u^\Pi(\cd)\big|_{[t_{k-1},t_k]}\big)\\
\ns\ds=\inf_{u^k(\cd)\in\cU[t_{k-1},t_k]}J^k\big(t_{k-1},\bar
X^\Pi(t_{k-1});u^k(\cd)\big)\le J^k\big(t_{k-1},\bar
X^\Pi(t_{k-1});u^k(\cd)\big)\\
\ns\ds=J\big(t_{k-1},\bar
X^\Pi(t_{k-1});u^k(\cd)\oplus\Psi^\Pi(\cd)\big|_{[t_k,T]}\big),
\qq\forall u^k(\cd)\in\cU[t_{k-1},t_k],\q1\le k\le N,\ea\ee
where $u^k(\cd)\oplus\Psi^\Pi(\cd)\big|_{[t_k,T]}$ is defined the
same way as (\ref{uPsi})--(\ref{4.4}). Similar to (\ref{J>J}), for
$k=1,2,\cds,N-1$, we have in general that
\bel{}J\big(t_{k-1},\bar X^\Pi(t_{k-1});\bar
u^\Pi(\cd)\big)>\inf_{u(\cd)\in\cU[t_{k-1},T]}J(t_{k-1},\bar
X^\Pi(t_{k-1});u(\cd)\big).\ee
Since the involved $N$ players in Problem (G$^\Pi$) interact through
the initial/terminal pairs $(t_k,X(t_k))$, $k=1,2,\cds,N-1$, one
should actually denote
\bel{wt J}J^k\big(t_{k-1},X(t_{k-1});u^k(\cd)\big)\equiv\wt
J^k\big(x;u^1(\cd),\cds,u^N(\cd)\big),\qq1\le k\le N.\ee
Hence, (\ref{Nash}) means that if we let
$$\bar u^k(\cd)=\bar u^\Pi(\cd)\big|_{[t_{k-1},t_k]},\qq1\le k\le
N,$$
then $(\bar u^1(\cd),\cds,\bar u^N(\cd))$ is a {\it Nash
equilibrium} of the $N$-person non-cooperative differential game
associated with $\wt J^k\big(x;u^1(\cd),\cds,u^N(\cd)\big)$, $1\le
k\le N$ (defined in (\ref{wt J})).

\subsection{The formal limits.}

We now would like to look at the situation when $\|\Pi\|\to0$.
Suppose we have the following:
\bel{limit-V}\lim_{\|\Pi\|\to0}\(|V^\Pi(t,x)-V(t,x)|+|V^\Pi_x(t,x)-V_x(t,x)|
+|V^\Pi_{xx}(t,x)-V_{xx}(t,x)|\)=0,\ee
uniformly for $(t,x)$ in any compact sets, for some $V(\cd\,,\cd)$.
Under (H3), we also have
\bel{}\lim_{\|\Pi\|\to0}|\Psi^\Pi(t,x)-\Psi(t,x)|=0,\ee
uniformly for $(t,x)$ in any compact sets, for
\bel{}\Psi(t,x)=\psi(t,t,x,V_x(t,x),V_{xx}(t,x)\big),\qq(t,x)\in[0,T]\times\dbR^n.\ee
Then the following limit exist:
$$\lim_{\|\Pi\|\to0}\|\bar X^\Pi(\cd)-\bar X(\cd)\|_{L^2_\dbF(\O;C([0,T];\dbR^n))}=0,$$
for $\bar X(\cd)$ solving the following SDE:
\bel{closed-loop}\left\{\ba{ll}
\ns\ds d\bar X(s)=b\big(s,\bar X(s),\bar u(s)\big)ds
+\si\big(s,\bar X(s),\bar u(s)\big)dW(s),\q s\in[0,T].\\
\ns\ds\bar X(0)=x,\ea\right.\ee
where
\bel{bar u}\bar u(s)=\Psi\big(s,\bar X(s)\big),\qq s\in[0,T],\ee
and
$$\ba{ll}
\ns\ds
L^2_\dbF(\O,C([0,T];\dbR^n))=\Big\{X:[0,T]\times\O\to\dbR^n\bigm|\hb{$X(\cd)$
has continuous paths},\\
\ns\ds\qq\qq\qq\qq\qq\qq\qq\qq\qq\qq\qq\dbE\big[\sup_{t\in[0,T]}|X(t)|^2\big]<\infty\Big\}.\ea$$
Clearly,
\bel{}\lim_{\|\Pi\|\to0}\|\bar u^\Pi(\cd)-\bar
u(\cd)\|_{\cU^2[0,T]}=0.\ee
By (\ref{Nash}), we have
$$J\big(\ell^\Pi(t),\bar X^\Pi(\ell^\Pi(t));\bar u^\Pi(\cd)\big)=V^\Pi(\ell^\Pi(t),\bar X^\Pi(t)),
\qq t\in[0,T].$$
Thus, passing to the limits, we have (\ref{J=V}). Also, we have the
following:
$$\ba{ll}
\ns\ds J\big(t_{k-1},\bar
X(t_{k-1});\Psi(\cd)\big|_{[t_{k-1},T]}\big)\equiv
J\big(t_{k-1},\bar
X(t_{k-1});\bar u(\cd)\big|_{[t_{k-1},T]}\big)\\
\ns\ds=\dbE_{t_{k-1}}\[\int_{t_{k-1}}^Tg(t_{k-1},s,\bar X(s),\bar
u(s))ds+h(t_{k-1},\bar X(T))\]\\
\ns\ds\le\dbE_{t_{k-1}}\[\int_{t_{k-1}}^Tg(t_{k-1},s,\bar
X^\Pi(s),\bar u^\Pi(s))ds+h(t_{k-1},\bar
X^\Pi(T))\]+R\big(\|\Pi\|\big)\\
\ns\ds=V^\Pi\big(t_{k-1},\bar
X^\Pi(t_{k-1})\big)+R\big(\|\Pi\|\big)=J^k\big(t_{k-1},\bar
X^\Pi(t_{k-1});\bar
u^\Pi(\cd)\big|_{[t_{k-1},T]}\big)+R\big(\|\Pi\|\big)\\
\ns\ds\le J^k\big(t_{k-1},\bar
X^\Pi(t_{k-1}),u^k(\cd)\big)+R\big(\|\Pi\|\big)\\
\ns\ds=J\big(t_{k-1},\bar
X(t_{k-1});u^k(\cd)\oplus\Psi(\cd)\big|_{[t_k,T]}\big)+R\big(\|\Pi\|\big),\qq\forall
u^k(\cd)\in\cU[t_{k-1},t_k],\ea$$
for some $R(r)$ with $R(r)\to0$ as $r\to0$. Hence, by Definition
4.1, $\Psi(\cd\,,\cd)$ is a time-consistent equilibrium strategy,
and $V(\cd\,,\cd)$ is a time-consistent equilibrium value function
of Problem (N).

\ms

In the rest of this subsection, we will formally pass to the limits
to find the equations that can be used to characterize the
equilibrium value function $V(\cd\,,\cd)$. To this end, let us first
write the equations for $\Th^{k+1}(\cd\,,\cd)$ in the integral
forms: For $k=1,2,\cds,N-1$, one has
\bel{}\ba{ll}
\ns\ds\Th^{k+1}\1n(t,x)\1n=\1n h(t_{k-1},\1n
x)\1n+\3n\int_t^T\3n\dbH\big(t_{k-1},\1n
s,x,\1n\psi\big(\ell^\Pi(s),\1n s,x,\1n V^\Pi_x(s,x),\1n
V^\Pi_{xx}(s,x)\big),\1n
\Th^{k+1}_x\1n(s,x),\1n\Th^{k+1}_{xx}\1n(s,x)\big)ds\\
\ns\ds\qq\qq=h(t_{k-1},x)+\int_t^T\(\lan
\Th^{k+1}_x(s,x),b\big(s,x,\psi(\ell^\Pi(s),s,x,
V^\Pi_x(s,x),V^\Pi_{xx}(s,x))\big)\ran\\
\ns\ds\qq\qq\qq\qq\qq\qq\qq+\tr\big[a\big(s,x,\psi(\ell^\Pi(s),s,x,
V^\Pi_x(s,x),V^\Pi_{xx}(s,x))\big)\Th^{k+1}_{xx}(s,x)\big]\\
\ns\ds\qq\qq\qq\qq\qq\qq\qq+g\big(t_{k-1},s,x,\psi(\ell^\Pi(s),s,x,
V^\Pi_x(s,x),V^\Pi_{xx}(s,x))\big)\)ds,\\
\ns\ds\qq\qq\qq\qq\qq\qq\qq\qq\qq\qq\qq\qq\qq\qq\qq(t,x)\in[t_k,T]\times\dbR^n.\ea\ee
Let us define
\bel{ThPi}\Th^\Pi(\t,t,x)=\sum_{k=1}^{N-1}\Th^{k+1}(t,x)I_{[t_{k-1},t_k)}(\t),\qq
(\t,t,x)\in D[0,T]\times\dbR^n.\ee
Then
\bel{4.48}\ba{ll}
\ns\ds\Th^\Pi(\t,t,x)=
h^\Pi(\t,x)+\int_t^T\Big\{\lan\Th^\Pi_x(\t,s,x),
b\big(s,x,\psi(\ell^\Pi(s),s,x,
V^\Pi_x(s,x),V^\Pi_{xx}(s,x))\big)\ran\\
\ns\ds\qq\qq\qq\qq\qq\qq+\tr\big[a\big(s,x,\psi(\ell^\Pi(s),s,x,V^\Pi_x(s,x),
V^\Pi_{xx}(s,x))\big)\Th^\Pi_{xx}(\t,s,x)\big]\\
\ns\ds\qq\qq\qq\qq\qq\qq+g^\Pi\big(\t,s,x,\psi(\ell^\Pi(s),s,x,
V^\Pi_x(s,x),V^\Pi_{xx}(s,x))\big)\Big\}ds,\ea\ee
where
$$\left\{\ba{ll}
\ns\ds
h^\Pi(\t,x)=\sum_{k=1}^{N-1}h(t_{k-1},x)I_{[t_{k-1},t_k)}(\t),\qq(\t,x)\in[0,T]\times\dbR^n,\\
\ns\ds
g^\Pi(\t,s,x,u)=\sum_{k=1}^{N-1}g(t_{k-1},s,x,u)I_{[t_{k-1},t_k)}(\t),\qq(\t,s,x,u)\in
D[0,T]\times\dbR^n\times U.\ea\right.$$
Let us look at $V^\Pi(\cd\,,\cd)$. For
$(t,x)\in[t_{N-1},t_N]\times\dbR^n$, we have
\bel{}\ba{ll}
\ns\ds V^\Pi(t,x)\1n=\1n
h(t_{N-1},x)\1n+\2n\int_t^T\3n\dbH\big(t_{N-1},s,x,\psi(\ell^\Pi(s),s,x,
V^\Pi_x(s,x),V^\Pi_{xx}(s,x)),V^\Pi_x(s,x),V^\Pi_{xx}(s,x)\big)ds\\
\ns\ds\qq\q~=\1n
h(\ell^\Pi(t),x)\1n+\2n\int_t^T\3n\dbH\big(\ell^\Pi(t),s,x,\psi(\ell^\Pi(s),s,x,
V^\Pi_x(s,x),V^\Pi_{xx}(s,x)),V^\Pi_x(s,x),V^\Pi_{xx}(s,x)\big)ds,\ea\ee
and for $(t,x)\in[t_{k-1},t_k)\times\dbR^n$, $k=1,2,\cds,N-1$,
$$\ba{ll}
\ns\ds
V^\Pi(t,x)=\Th^{k+1}(t_k,x)+\int_t^{t_k}\dbH\(\ell^\Pi(s),s,x,\psi\big(\ell^\Pi(s),s,x,
V^\Pi_x(s,x),V^\Pi_{xx}(s,x)\big),V^\Pi_x(s,x),V^\Pi_{xx}(s,x)\)ds\\
\ns\ds\qq\qq=\Th^\Pi(t_{k-1},t_k,x)+\int_t^{t_k}\Big\{\lan
V^\Pi_x(s,x),b\big(s,x,\psi(\ell^\Pi(s),s,x,
V^\Pi_x(s,x),V^\Pi_{xx}(s,x))\big)\ran\\
\ns\ds\qq\qq\qq\qq\qq\qq\qq\qq+\tr\big[a\big(s,x,\psi(\ell^\Pi(s),s,x,V^\Pi_x(s,x),
V^\Pi_{xx}(s,x))\big)V^\Pi_{xx}(s,x)\big]\\
\ns\ds\qq\qq\qq\qq\qq\qq\qq\qq+g^\Pi\big(t_{k-1},s,x,\psi(\ell^\Pi(s),s,x,
V^\Pi_x(s,x),V^\Pi_{xx}(s,x))\big)\Big\}ds.\ea$$
On the other hand, for $(t,x)\in[t_{k-1},t_k)\times\dbR^n$, we have
$$\ba{ll}
\ns\ds
\Th^\Pi(t_{k-1},t,x)=\Th^\Pi(t_{k-1},t_k,x)+\int_t^{t_k}\Big\{\lan
\Th^\Pi_x(t_{k-1},s,x),b\big(s,x,\psi(\ell^\Pi(s),s,x,
V^\Pi_x(s,x),V^\Pi_{xx}(s,x))\big)\ran\\
\ns\ds\qq\qq\qq\qq\qq\qq\qq\qq+
\tr\big[a\big(s,x,\psi(\ell^\Pi(s),s,x,V^\Pi_x(s,x),
V^\Pi_{xx}(s,x))\big)\Th^\Pi_{xx}(t_{k-1},s,x)\big]\\
\ns\ds\qq\qq\qq\qq\qq\qq\qq\qq+g^\Pi\big(t_{k-1},s,x,\psi(\ell^\Pi(s),s,x,
V^\Pi_x(s,x),V^\Pi_{xx}(s,x))\big)\Big\}ds.\ea$$
Therefore, in the case that (\ref{4.48}) is well-posed, we have
\bel{}V^\Pi(t,x)=\Th^\Pi(t_{k-1},t,x),\qq(t,x)\in[t_{k-1},t_k)\times\dbR^n,\ee
or equivalently,
\bel{}V^\Pi(t,x)=\Th^\Pi(\ell^\Pi(t),t,x),\qq(t,x)\in[0,t_{N-1})\times\dbR^n.\ee
Now, let us assume (\ref{limit-V}) holds for some $V(\cd\,,\cd)$ and
\bel{limit1}\ba{ll}
\ns\ds\lim_{\|\Pi\|\to0}\(|\Th^\Pi(\t,t,x)-\Th(\t,t,x)|+|\Th_x(\t,t,x)-\Th_x(\t,t,x)|
+|\Th_{xx}(\t,t,x)-\Th_{xx}(\t,t,x)|\)=0,\\
\ns\ds\qq\qq\qq\qq\qq\qq\qq\qq\qq\qq\qq\qq\q(\t,t,x)\in
D[0,T]\times\dbR^n,\ea\ee
for some $\Th(\cd\,,\cd\,,\cd)$. Let us assume that
\bel{}|h_\t(\t,x)|+|g_\t(\t,t,x,u)|\le K,\qq\forall(\t,t,x,u)\in
D[0,T]\times\dbR^n\times U.\ee
Then we have
$$\lim_{\|\Pi\|\to0}h^\Pi(\t,x)=h(\t,x),\qq(\t,x)\in[0,T]\times\dbR^n,$$
and
$$\ba{ll}
\ns\ds\lim_{\|\Pi\|\to0}g^\Pi\big(\t,s,x,\psi(\ell^\Pi(s),s,x,
V^\Pi_x(s,x),V^\Pi_{xx}(s,x))\big)=g\big(\t,s,x,\psi(s,s,x,V_x(s,x),V_{xx}(s,x))\big),\\
\ns\ds\qq\qq\qq\qq\qq\qq\qq\qq\qq\qq(\t,s,x)\in
D[0,T]\times\dbR^n.\ea$$
Consequently, we obtain the following integro-partial differential
equation for $\Th(\cd\,,\cd\,,\cd)$:
\bel{4.53}\ba{ll}
\ns\ds\Th(\t,t,x)=h(\t,x)+\int_t^T\dbH\big(\t,s,x,\psi(s,s,x,\Th_x(s,s,x),\Th_{xx}(s,s,x)),
\Th_x(\t,s,x),\Th_{xx}(\t,s,x)\big)ds,\\
\ns\ds\qq\qq\qq\qq\qq\qq\qq\qq\qq\qq\qq(\t,t,x)\in
D[0,T]\times\dbR^n,\ea\ee
and relation
\bel{V=Th}V(t,x)=\Th(t,t,x),\qq(t,x)\in[0,T]\times\dbR^n.\ee
It is clear that (\ref{4.53}) is an integral form of the following
differential equation:
\bel{HJBe}\left\{\ba{ll}
\ns\ds\Th_t(\t,t,x)+\dbH\big(\t,t,x,\psi\big(t,t,x,\Th_x(t,t,x),\Th_{xx}
(t,t,x)\big),
\Th_x(\t,t,x),\Th_{xx}(\t,t,x)\big)=0,\\
\ns\ds\qq\qq\qq\qq\qq\qq\qq\qq\qq\qq\qq(\t,t,x)\in D[0,T]\times\dbR^n,\\
\ns\ds\Th(\t,T,x)=h(\t,x),\qq(\t,x)\in[0,T]\times\dbR^n,\ea\right.\ee
where
\bel{H(t)}\left\{\ba{ll}
\ns\ds\dbH(\t,t,x,u,p,P)=\tr\big[a(t,x,u)P\big]+\lan
b(t,x,u),p\ran+g(\t,t,x,u),\\
\ns\ds\psi(\t,t,x,p,P)\in\arg\min\dbH(\t,t,x,\cd\,,p,P).\ea\right.\ee
Therefore, we may equivalently write (\ref{HJBe}) as follows:
\bel{HJBe1}\left\{\ba{ll}
\ns\ds\Th_t(\t,t,x)+\lan\Th_x(\t,t,x),
b\big(t,x,\psi(t,t,x,\Th_x(t,t,x),\Th_{xx}(t,t,x))
\big)\ran\\
\ns\ds\qq+\tr\big[a\big(t,x,\psi(t,t,x,\Th_x(t,t,x),\Th_{xx}(t,t,x))
\big)\Th_{xx}(\t,t,x)\big]\\
\ns\ds\qq+g\big(\t,t,x,\psi(t,t,x,\Th_x(t,t,x),\Th_{xx}(t,t,x))
\big)=0,\qq\qq(\t,t,x)\in D[0,T]\times\dbR^n,\\
\ns\ds\Th(\t,T,x)=h(\t,x),\qq\qq(\t,x)\in[0,T]\times\dbR^n.\ea\right.\ee
We call the above (\ref{HJBe1}) the {\it equilibrium
Hamilton-Jacobi-Bellamn equation} (equilibrium HJB equation, for
short) of Problem (N). If one can find $\Th(\cd\,,\cd\,,\cd)$ from
the above, then the equilibrium value function $V(\cd\,,\cd)$ can be
determined by (\ref{V=Th}), and the (time-consistent) equilibrium
pair can be determined by (\ref{closed-loop}) and (\ref{bar u}), in
principle.

\ms

Let us make some remarks on (\ref{HJBe1}).

\ms

(i) It is an interesting feature of (\ref{HJBe1}) that both
$\Th(\t,t,x)$ and $\Th(t,t,x)$ appear in the equation where the
later is the restriction of the former on $\t=t$. On one hand,
although the equation is fully nonlinear, due to the fact that
$\Th(t,t,x)$ is different from $\Th(\t,t,x)$, the existing theory
for fully nonlinear parabolic equations cannot apply directly. On
the other hand, it is seen that if $\Th(t,t,x)$ is obtained from an
independent way, then (\ref{HJBe1}) is actually a linear equation
for $\Th(\t,t,x)$ with $\t$ can be purely regarded as a parameter.

\ms

(ii) In the case that $\cD(\psi)$ is not equal to
$D[0,T]\times\dbR^n\times\dbR^n\times\cS^n$, the condition
\bel{}(\t,t,x,\Th_x(t,t,x),\Th_{xx}(t,t,x))\in\cD(\psi)\ee
has to be regarded as a part of the solution. We will see that for
some interesting special cases, the above condition can come
automatically. More generally, we may also write (\ref{HJBe}) as
\bel{}\left\{\ba{ll}
\ns\ds\Th_t(\t,t,x)+\dbH\big(\t,t,x,\arg\min\dbH\big(t,t,x,\Th_x(t,t,x),\Th_{xx}
(t,t,x)\big),
\Th_x(\t,t,x),\Th_{xx}(\t,t,x)\big)\ni0,\\
\ns\ds\qq\qq\qq\qq\qq\qq\qq\qq(\t,t,x)\in D[0,T]\times\dbR^n,\\
\ns\ds\Th(\t,T,x)=h(\t,x),\qq(\t,x)\in[0,T]\times\dbR^n,\ea\right.\ee
since the set $\arg\min\dbH\big(t,t,x,\Th_x(t,t,x),\Th_{xx}
(t,t,x)\big)$ might contain more than one point.

\ms

(iii) In the case $\cD(\psi)\ne\cD(H)$, according to Lemma 3.2, we
can define
$$\dbH^\e(\t,t,x,u,p,P)=\dbH(\t,t,x,u,p,P)+\e|u|^2,$$
and there exists a
$\psi^\e:D[0,T]\times\dbR^n\times\dbR^n\times\cS^n\to U$ such that
$$\dbH^\e(\t,t,x,\psi^\e(\t,t,x,p,P),p,P)=\inf_{u\in
U}\dbH^\e(\t,t,x,u,p,P)\equiv H^\e(\t,t,x,p,P).$$
Further,
$$\lim_{\e\to0}H^\e(\t,t,x,p,P)=H(\t,t,x,p,P),\qq
\lim_{\e\to0}\e|\psi^\e(\t,t,x,p,P)|^2=0.$$
It is not hard to see that the above actually amounts to defining
$$g^\e(\t,t,x,u)=g(\t,t,x,u)+\e|u|^2.$$
We may refer to the corresponding problem as a {\it regularized}
problem. If the corresponding equilibrium value function is denoted
by $V^\e(\cd\,,\cd)$, then it is expected that
\bel{}\lim_{\e\da0}V^\e(t,x)=V(t,x),\qq(t,x)\in[0,T]\times\dbR^n.\ee
However, in general, if $(\bar X^\e(\cd),\bar u^\e(\cd))$ is an
equilibrium pair for the regularized problem, we might not have the
limit of $\bar u^\e(\cd)$ as $\e\to0$. In this case, we should be
satisfied by the above characterization of the equilibrium value
function $V(\cd\,,\cd)$, and $\bar u^\e(\cd)$ can be regarded as
some kind of ``near equilibrium control''.

\ms

(iv) For the case
\bel{}\si(t,x,u)=\si(t,x),\qq\forall(t,x,u)\in[0,T]\times\dbR^n\times
U,\ee
i.e., the control does not enter the diffusion of the state
equation, $\psi(\cd)$ is independent of $P$ and
\bel{}\psi(t,x,p)\in\arg\min\big[\lan
p,b(t,x,\cd)\ran+g(\t,t,x,\cd)\big].\ee
Then the equilibrium HJB equation can be written as
\bel{4.57}\left\{\ba{ll}
\ns\ds\Th_t(\t,t,x)+\tr\big[a(t,x)\Th_{xx}(\t,t,x)\big]+\lan
b\big(t,x,\psi(t,t,x,\Th_x(t,t,x))\big),\Th_x(\t,t,x)\ran\\
\ns\ds\qq\qq+g\big(\t,t,x,\psi(t,t,x,\Th_x(t,t,x))\big)=0,\qq\qq
(\t,t,x)\in D[0,T]\times\dbR^n,\\
\ns\ds\Th(\t,T,x)=h(\t,x),\qq(\t,x)\in[0,T]\times\dbR^n.\ea\right.\ee
We will carefully discuss this case in the next section. Note that
for a deterministic problem, namely Problem (N) for an ordinary
differential equation system, we may take
$$\si(t,x,u)=\e I,\qq(t,x,u)\in[0,T]\times\dbR^n\times U,$$
for $\e>0$ to regularize the problem. Then the corresponding
equilibrium HJB equation reads
\bel{}\left\{\ba{ll}
\ns\ds\Th^\e_t(\t,t,x)+{1\over2}\D\Th^\e_{xx}(\t,t,x)+\lan
b\big(t,x,\psi(t,t,x,\Th^\e_x(t,t,x))\big),\Th^\e_x(\t,t,x)\ran\\
\ns\ds\qq\qq+g\big(\t,t,x,\psi(t,t,x,\Th^\e_x(t,t,x))\big)=0,\qq\qq
(\t,t,x)\in D[0,T]\times\dbR^n,\\
\ns\ds\Th^\e(\t,T,x)=h(\t,x),\qq(\t,x)\in[0,T]\times\dbR^n.\ea\right.\ee
It is expected that $\Th^\e(\cd\,,\cd\,,\cd)\to\Th(\cd\,,\cd\,,\cd)$
in some sense, as $\e\to0$, with
\bel{}\left\{\ba{ll}
\ns\ds\Th_t(\t,t,x)+\lan
b\big(t,x,\psi(t,t,x,\Th_x(t,t,x))\big),\Th_x(\t,t,x)\ran+g\big(\t,t,x,\psi(t,t,x,\Th_x(t,t,x))\big)=0,\\
\ns\ds\qq\qq\qq\qq\qq\qq\qq\qq\qq\qq
(\t,t,x)\in D[0,T]\times\dbR^n,\\
\ns\ds\Th(\t,T,x)=h(\t,x),\qq(\t,x)\in[0,T]\times\dbR^n.\ea\right.\ee
At the moment, it is not clear to us how one can define viscosity
solution to the above equation.

\section{Well-Posedness of the Equilibrium HJB Equation}

In this section, we discuss the well-posedness for the equilibrium
HJB equation (\ref{HJBe1}). Let us first intuitively describe our
idea. For any smooth function $v(\cd\,,\cd)$, denote
\bel{}\ba{ll}
\ns\ds\big[\cL\big(t,v(t,\cd)\big)\f(\cd)\big](x)
=\tr\big[a(t,x,\psi(t,t,x,v_x(t,x),v_{xx}(t,x))
\f_{xx}(x)\big]\\
\ns\ds\qq\qq\qq\qq\qq\qq+\lan
b\big(t,x,\psi(t,t,x,v_x(t,x))\big),\f_x(x)\ran,\qq(t,x)\in[0,T]\times\dbR^n,\ea\ee
and
\bel{}\cG(\t,t,v(t,\cd))(x)=g\big(\t,t,x,\psi(t,t,x,v_x(t,x))\big),\qq(\t,t,x)
\in D[0,T]\times\dbR^n.\ee
Consider the following linear abstract backward evolution equation:
\bel{6.5}\left\{\ba{ll}
\ns\ds\Th_t(\t,t)+\cL(t,v(t))\Th(\t,t)+\cG(\t,t,v(t))=0,\qq
t\in[\t,T],\\
\ns\ds\Th(\t,T)=h(\t).\ea\right.\ee
Under some mild conditions, the above is well-posed, and we have the
following variation of constant formula:
\bel{}\Th(\t,t)=\cE(T,t;v(\cd))h(\t)+\int_t^T\cE(s,t;v(\cd))\cG(\t,s,v(s))ds,\qq
t\in[\t,T],\ee
where $\cE(\cd\,,\cd\,;v(\cd))$ is called the {\it backward
evolution operator} generated by $\cL(\cd\,,v(\cd))$. Consequently,
the (time-consistent) equilibrium value function
$V(t,\cd)=\Th(t,t,\cd)$ should be the solution to the following
nonlinear functional integral equation:
\bel{V(t)}V(t)=\cE(T,t;V(\cd))h(t)+\int_t^T\cE(s,t;V(\cd))\cG(t,s,V(s))ds,\qq
t\in[0,T].\ee
We call (\ref{V(t)}) the {\it equilibrium HJB integral equation} for
Problem (N). Once a solution $V(\cd\,,\cd)$ of (\ref{V(t)}) is
found, we can, in principle, construct a (time-consistent)
equilibrium control and an equilibrium pair for Problem (N). Of
course, if we like, we may also solve the equilibrium HJB equation
(\ref{HJBe1}), which actually is not necessary as far as the
construction of a time-consistent equilibrium pair is concerned.

\ms

The well-posedness of (\ref{V(t)}) seems to be difficult for the
general case. At the moment, we do not have a complete solution for
that and hopefully, we can present some satisfactory results for the
equilibrium HJB integral equation (\ref{V(t)}) in our future
publications. On the other hand, in the rest of this section, we are
going to present a well-posedness result for an interesting special
case of (\ref{V(t)}), from which one can get some taste of the
problem. The main hypothesis that we will assume below is the
following:
\bel{si=si}\si(t,x,u)=\si(t,x),\qq(t,x,u)\in[0,T]\times\dbR^n\times
U,\ee
namely, the control does not enter the diffusion of the state
equation. As we discussed in Section 4, in this case, our
equilibrium HJB equation reads
\bel{6.8}\left\{\ba{ll}
\ns\ds\Th_t(\t,t,x)+{1\over2}\tr\big[\si(t,x)\si(t,x)^T\Th_{xx}(\t,t,x)\big]
+\lan
b\big(t,x,\psi(t,t,x,\Th_x(t,t,x))\big),\Th_x(\t,t,x)\ran\\
\ns\ds\qq\qq+g\big(\t,t,x,\psi(t,t,x,\Th_x(t,t,x))\big)=0,\qq\qq
(\t,t,x)\in D[0,T]\times\dbR^n,\\
\ns\ds\Th(\t,T,x)=h(\t,x),\qq(\t,x)\in[0,T]\times\dbR^n.\ea\right.\ee
The essential feature of (\ref{6.8}) is that $\Th_{xx}(t,t,x)$ does
not appear in the equation (although $\Th_x(t,t,x)$ still appears
there). This leads to the well-posedness problem much more
accessible. Further, from Example 3.5, we see that there are cases
for which $\psi$ is as smooth as the coefficients and
$b(t,x,\psi(t,t,x,p))$ is bounded. Therefore, the case that we are
going to consider below, although very special, includes a big class
of problems.

\ms

To avoid heavy notations, let us consider the following equation
\bel{5.10}\left\{\ba{ll}
\ns\ds\Th_t(\t,t,x)+\tr\big[a(t,x)\Th_{xx}(\t,t,x)\big]+\lan
b\big(t,x,\Th_x(t,t,x)\big),\Th_x(\t,t,x)\ran+g\big(\t,t,x,\Th_x(t,t,x)\big)=0,\\
\ns\ds\qq\qq\qq\qq\qq\qq\qq\qq\qq
(\t,t,x)\in D[0,T]\times\dbR^n,\\
\ns\ds\Th(\t,T,x)=h(\t,x),\qq(\t,x)\in[0,T]\times\dbR^n,\ea\right.\ee
with
\bel{notation}\left\{\ba{ll}
\ns\ds a(t,x)={1\over2}\si(t,x)\si(t,x)^T,\qq(t,x)\in[0,T]\times\dbR^n,\\
\ns\ds b(t,x,p)=b(t,x,\psi(t,t,x,p)),\qq(t,x,p)\in[0,T]\times\dbR^n\times\dbR^n,\\
\ns\ds g(\t,t,x,p)=g(\t,t,x,\psi(t,t,x,p)),\qq(\t,t,x,p)\in
D[0,T]\times\dbR^n\times\dbR^n.\ea\right.\ee
To investigate the well-posedness of (\ref{5.10}) above, let us make
some preparations. Let $C^\a(\dbR^n)$ be the space of all continuous
functions $\f:\dbR^n\to\dbR$ such that
$$\|\f\|_\a\equiv\|\f\|_0+[\f]_\a<\infty,$$
where
$$\|\f\|_0=\sup_{x\in\dbR^n}|\f(x)|,\q[\f]_\a=\sup_{x,y\in\dbR^n,x\ne
y}{|\f(x)-\f(y)|\over|x-y|^\a}.$$
Further, let $C^{1+\a}(\dbR^n)$ and $C^{2+\a}(\dbR^n)$ be the spaces
of all functions $\f:\dbR^n\to\dbR$ such that
$$\|\f\|_{1+\a}\equiv\|\f\|_0+\|\f_x\|_0+[\f_x]_\a<\infty,$$
and
$$\|\f\|_{2+\a}\equiv\|\f\|_0+\|\f_x\|_0+\|\f_{xx}\|_0+[\f_{xx}]_\a<\infty,$$
respectively. Next, let $B([0,T];C^\a(\dbR^n))$ be the set of all
measurable functions $f:[0,T]\times\dbR^n\to\dbR$ such that for each
$t\in[0,T]$, $f(t,\cd)\in C^\a(\dbR^n)$ and
$$\|f(\cd\,,\cd)\|_{B([0,T];C^\a(\dbR^n))}=\sup_{t\in[0,T]}\|f(t,\cd)\|_\a<\infty.$$
Also, we let $C([0,T];C^\a(\dbR^n))$ be the set of all continuous
functions that are also in $B([0,T];C^\a(\dbR^n))$. Thus,
$$C([0,T];C^\a(\dbR^n))\subseteq B([0,T];C^\a(\dbR^n)).$$
Similarly, we define $B([0,T];C^{k+\a}(\dbR^n))$ and
$C([0,T];C^{k+\a}(\dbR^n))$, respectively, for $k=1,2$.

\ms

We introduce the following hypotheses for the above equation
(\ref{5.10}).

\ms

{\bf(P)} The maps $a:[0,T]\times\dbR^n\to\cS^n$,
$b:[0,T]\times\dbR^n\times\dbR^n\to\dbR^n$,
$g:D[0,T]\times\dbR^n\times\dbR^n\to\dbR$ and
$h:[0,T]\times\dbR^n\to\dbR$ are continuous and bounded. Moreover,
there exists a constant $L>0$ such that
\bel{}\ba{ll}
\ns\ds|a_x(t,x)|+|b_x(t,x,p)|+|g_x(\t,t,x,p)|+|b_p(t,x,p)|+|g_p(\t,t,x,p)|+|h_x(\t,x)|\le
L,\\
\ns\ds\qq\qq\qq\qq\qq\qq\qq\qq(\t,t,x,p)\in
D[0,T]\times\dbR^n\times\dbR^n.\ea\ee
Further, $a(t,x)^{-1}$ exists for all $(t,x)\in[0,T]\times\dbR^n$
and there exist constants $\l_0,\l_1>0$ such that
\bel{5.12}\l_0 I\le
a(t,x)^{-1}\le\l_1I,\qq\forall(t,x)\in[0,T]\times\dbR^n.\ee

\ms

We point out here that some of the conditions assumed in (P) can be
substantially relaxed. However, we prefer not to get into those
generalities for the sake of simplicity in our presentation. Note
also that typically, the ellipticity condition of $a(t,x)$ looks
like
$$a(t,x)\ge\d I,\qq\forall(t,x)\in[0,T]\times\dbR^n,$$
for some $\d>0$. It is clear that when $a(\cd\,,\cd)$ is assumed to
be bounded, then the above is equivalent to (\ref{5.12}). The number
$\l_0$ in (\ref{5.12}) will be used below.

\ms

For any $v(\cd\,,\cd)\in C([0,T];C^{1+\a}(\dbR^n))$, we consider the
following linear PDE, parameterized by $\t\in[0,T)$:
\bel{5.13}\left\{\ba{ll}
\ns\ds\Th_t(\t,t,x)+\cL[t,v(\cd)]\Th(\t,t,x)+g\big(\t,t,x,v_x(t,x)\big)=0,\qq
(t,x)\in[0,T)\times\dbR^n,\\
\ns\ds\Th(\t,T,x)=h(\t,x),\qq\qq\qq x\in\dbR^n,\ea\right.\ee
where the differential operator $\cL[t,v(\cd)]$ is defined by the
following:
\bel{cL}\cL[t,v(\cd)]\f(x)=\tr\big[a(t,x)\f_{xx}(x)\big]+\lan
b(t,x,v_x(t,x)),\f_x(x)\ran,\qq\forall\f(\cd)\in C^2(\dbR^n).\ee
In what follows, we let
\bel{C_l}C_\l(\dbR^n)=\Big\{\f(\cd)\in
C(\dbR^n)\bigm|\sup_{x\in\dbR^n}e^{-\l|x|^2}
|\f(x)|<\infty\Big\}.\ee
We have the following result whose proof follows a relevant one
found in \cite{Friedman 1964}, with some minor modifications.

\ms

\bf Proposition 5.1. \sl Let {\rm(P)} hold and fix a $\t\in[0,T)$.
Then for any $v(\cd\,,\cd)\in C([0,T];C^{1+\a}(\dbR^n))$ and any
$h(\t,\cd)\in C_\l(\dbR^n)$ with $\l<{\l_0\over4T}$, Cauchy problem
$(\ref{5.13})$ admits a unique classical solution
$\Th(\t,\cd\,,\cd)$ and the following representation holds:
\bel{variation of constant}\ba{ll}
\ns\ds\Th(\t,t,x)=\int_{\dbR^n}\G^{v(\cd)}(t,x;T,y)h(\t,y)dy+\int_t^T
\int_{\dbR^n}\G^{v(\cd)}(t,x;s,y)g(\t,s,y,v_x(s,y))dyds,\\
\ns\ds\qq\qq\qq\qq\qq\qq\qq\qq\qq\qq\qq(t,x)\in[0,T]\times\dbR^n.\ea\ee
Here, $\G^{v(\cd)}(t,x;s,y)$ is defined on
$[0,T]\times\dbR^n\times[0,T]\times\dbR^n$ with $t<s$ having the
following properties:
 \ms

{\rm(i)} For any fixed $(s,y)\in(0,T]\times\dbR^n$,
$$\G^{v(\cd)}_x(t,x;s,y),~\G^{v(\cd)}_{xx}(t,x;s,y),~\G^{v(\cd)}_t(t,x;s,y)$$
are continuous in $(t,x),(s,y)\in[0,T]\times\dbR^n$ with $t<s$, and
for any fixed $(s,y)\in(0,T]\times\dbR^n$,
\bel{LG=0}\cL[t,v(\cd)]\G^{v(\cd)}(t,x;s,y)=0,\qq(t,x)\in[0,s)\times\dbR^n,\ee

{\rm(ii)} For any $\f(\cd)\in C_\l(\dbR^n)$ with $\l<{\l_0\over4
T}$,
\bel{terminal}\lim_{t\ua
s}\int_{\dbR^n}\G^{v(\cd)}(t,x;s,y)\f(y)dy=\f(x),\qq x\in\dbR^n,\ee

\rm

The map $\G^{v(\cd)}(t,x;s,y)$ in the above proposition is called
the {\it fundamental solution} to the problem (\ref{5.13}).

\ms

Our main result of this section is the following.

\ms

\bf Theorem 5.2. \sl Let {\rm(P)} hold. Then $(\ref{5.10})$ admits a
unique solution $\Th(\cd\,,\cd\,,\cd)$.

\rm

\ms

\it Proof. \rm Let
\bel{cL_0}\cL_0(t)\f(x)=\tr\big[a(t,x)\f_{xx}(x)\big],\qq\forall\f(\cd)\in
C^2(\dbR^n),\ee
which is independent of $v(\cd)$. Then
\bel{}\cL[t,v(\cd)]\f(x)=\cL_0(t)\f(s)+\lan
b(t,x,v_x(t,x)),\f_x(x)\ran,\qq\forall \f(\cd)\in C^2(\dbR^n),\ee
and (\ref{5.13}) can be written as (fix $\t\in[0,T)$)
\bel{5.22}\left\{\ba{ll}
\ns\ds\Th_t(\t,t,x)+\cL_0(t)\Th(\t,t,x)+\lan
b(t,x,v_x(t,x)),\Th_x(\t,t,x)\ran
+g\big(\t,t,x,v_x(t,x)\big)=0,\\
\ns\ds\qq\qq\qq\qq\qq\qq\qq\qq\qq(t,x)\in[0,T)\times\dbR^n,\\
\ns\ds\Th(\t,T,x)=h(\t,x),\qq\qq\qq x\in\dbR^n,\ea\right.\ee
Applying Proposition 5.1, we have
\bel{}\ba{ll}
\ns\ds\Th(\t,t,x)=\int_{\dbR^n}\G^0(t,x;T,y)h(\t,y)dy+\int_t^T\3n\int_{\dbR^n}
\G^0(t,x;s,y)\big[\lan b(s,y,v_x(s,y)),\Th_x(\t,s,y)\ran\\
\ns\ds\qq\qq\qq\qq\qq\qq\qq\qq\qq\qq+g(\t,s,y,v_x(s,y))\big]dyds,\qq(t,x)\in[\t,T]\times
\dbR^n,\ea\ee
where $\G^0(t,x;s,y)$ is the fundamental solution of $\cL_0(\cd)$,
given by the following explicitly:
\bel{G0}\ba{ll}
\ns\ds
\G^0(t,x;s,y)={1\over\big(4\pi(s-t)\big)^{n\over2}\big(\det[a(s,y)]
\big)^{1\over2}} \,e^{-{\lan
a(s,y)^{-1}(x-y),(x-y)\ran\over4(s-t)}},\\
\ns\ds\qq\qq\qq\qq\qq\qq\qq\qq\qq\qq(t,x),(s,y)\in[0,T]\times
\dbR^n,~t<s.\ea\ee
Direct computations show that (\cite{Friedman 1964})
\bel{5.25}\left\{\ba{ll}
%
%
\ns\ds|\G^0(t,x;s,y)|\le {K\over
(s-t)^{n\over2}}\,e^{-{\l|x-y|^2\over4(s-t)}},\\
\ns\ds|\G^0_x(t,x;s,y)|\le{K\over
(s-t)^{n+2\over2}}\,e^{-{\l|x-y|^2\over4(s-t)}},\ea\right.\qq\qq\l<\l_0.\ee
Moreover,
\bel{}\ba{ll}
\ns\ds\G^0_y(t,x;s,y)=-\G^0_x(t,x;s,y)-\G^0(t,x;s,y)
\[\big(\det\big[a(s,y)\big]\big)_y+\lan\big[a(s,y)^{-1}\big]_y(x-y),
x-y\ran\]\\
\ns\ds\qq\qq\q~\equiv-\G^0_x(t,x;s,y)-\G^0(t,x;s,y)\rho(t,x,s,y),\ea\ee
where
$$\lan\big[a(s,y)^{-1}\big]_y(x-y),
x-y\ran=\pmatrix{\lan[a(s,y)^{-1}]_{y_1}(x-y),x-y\ran\cr\lan[a(s,y)^{-1}]_{y_2}(x-y),x-y\ran
\cr\vdots\cr\lan[a(s,y)^{-1}]_{y_n}(x-y),x-y\ran},$$
and
$$\rho(t,x,s,y)={1\over2}\big(\det\big[a(s,y)\big]\big)_y+{\lan\big[a(s,y)^{-1}\big]_y(x-y),
x-y\ran\over4(s-t)}.$$
Under (P), we see that
\bel{}|\rho(s,x,y)|\le K\(1+{|x-y|^2\over
s-t}\),\qq\forall(s,x,y)\in[0,T]\times\dbR^n\times\dbR^n.\ee
Then
$$\ba{ll}
\ns\ds\Th_x(\t,t,x)=\int_{\dbR^n}\G^0_x(t,x;T,y)h(\t,y)dy+\int_t^T\int_{\dbR^n}
\G_x^0(t,x;s,y)\big[\lan b(s,y,v_x(s,y)),\Th_x(\t,s,y)\ran\\
\ns\ds\qq\qq\qq\qq\qq\qq\qq\qq\qq+g(\t,s,y,v_x(s,y))\big]dyds\\
\ns\ds=-\int_{\dbR^n}\G^0_y(t,x;T,y)h(\t,y)dy-\int_{\dbR^n}\G^0(t,x;T,y)\rho(t,x,T,y)
h(\t,y)dy\\
\ns\ds\qq\qq+\int_t^T\int_{\dbR^n}
\G_x^0(t,x;s,y)\big[\lan b(s,y,v_x(s,y)),\Th_x(\t,s,y)\ran+g(\t,s,y,v_x(s,y))\big]dyds\\
\ns\ds=\int_{\dbR^n}\G^0(t,x;T,y)h_y(\t,y)dy-\int_{\dbR^n}\G^0(t,x;s,y)\rho(t,x,T,y)
h(\t,y)dy\\
\ns\ds\qq\qq+\int_t^T\int_{\dbR^n} \G_x^0(t,x;s,y)\big[\lan
b(s,y,v_x(s,y)),\Th_x(\t,s,y)\ran+g(\t,s,y,v_x(s,y))\big]dyds. \ea$$
Hence,
\bel{5.28}\ba{ll}
\ns\ds|\Th_x(\t,t,x)|\le\int_{\dbR^n}{Ke^{-{\l|x-y|^2\over4(T-t)}}
\over(T-t)^{n\over2}}\[|h_y(\t,y)|+\(1+{|x-y|^2\over T-t}\)|h(\t,y)|\]dy\\
\ns\ds\qq\qq\qq\qq\qq+\int_t^T\3n\int_{\dbR^n}{Ke^{-{\l|x-y|^2\over4(s-t)}}\over(s-t)^{n+1\over2}}
\big(|\Th_x(\t,s,y)|+1\big)dyds\\
\ns\ds\le\int_{\dbR^n}{Ke^{-{\l|x-y|^2\over4(T-t)}}
\over(T-t)^{n\over2}}\[\(1+{|x-y|^2\over
T-t}\)|h(\t,y)|+|h_y(\t,y)|\]dy+\int_t^T\int_{\dbR^n}
{Ke^{-{\l|x-y|^2\over4(s-t)}}\over(s-t)^{n+1\over2}}
|\Th_x(\t,s,y)|dyds\\
\ns\ds\le
K\big(1+\|h(\t,\cd)\|_{C^1(\dbR^n)}\big)+\int_t^T\int_{\dbR^n}
{Ke^{-{\l|x-y|^2\over4(s-t)}}\over(s-t)^{n+1\over2}}
|\Th_x(\t,s,y)|dyds.\ea\ee
The above can be iterated as follows:
\bel{5.29}\ba{ll}
\ns\ds|\Th_x(\t,t,x)|\le
K\big(1+\|h(\t,\cd)\|_{C^1(\dbR^n)}\big)+\int_t^T\int_{\dbR^n}
{Ke^{-{\l|x-y|^2\over4(s-t)}}\over(s-t)^{n+1\over2}}
|\Th_x(\t,s,y)|dyds\\
\ns\ds\le
K\big(1+\|h(\t,\cd)\|_{C^1(\dbR^n)}\big)+\int_t^T\int_{\dbR^n}
{Ke^{-{\l|x-y|^2\over4(s-t)}}\over(s-t)^{n+1\over2}}
K\big(1+\|h(\t,\cd)\|_{C^1(\dbR^n)}\big)dyds\\
\ns\ds\qq\qq+\int_t^T\int_{\dbR^n}{Ke^{-{\l|x-y|^2\over4(s-t)}}\over(s-t)^{n+1\over2}}
\int_s^T\int_{\dbR^n}{Ke^{-{\l|y-z|^2\over4(r-s)}}\over(r-s)^{n+1\over2}}
|\Th_x(\t,r,z)|dzdrdyds\\
\ns\ds\le
K\big(1+\|h(\t,\cd)\|_{C^1(\dbR^n)}\big)\\
\ns\ds\qq\qq+\int_t^T\int_{\dbR^n}\(\int_t^r\int_{\dbR^n}
{Ke^{-{\l|x-y|^2\over4(s-t)}}\over(s-t)^{n+1\over2}}
{Ke^{-{\l|y-z|^2\over4(r-s)}}\over(r-s)^{n+1\over2}}dyds\)
|\Th_x(\t,r,z)|dzdr\\
\ns\ds\le
K\big(1+\|h(\t,\cd)\|_{C^1(\dbR^n)}\big)+\int_t^T\int_{\dbR^n}{Ke^{-{\l|x-z|^2\over4(r-t)}}
\over(r-t)^{n\over2}}|\Th_x(\t,r,z)|dzdr.\ea\ee
In the above, we have used Lemma 3 of Chapter 1 in \cite{Friedman
1964}. We can repeat the above procedure $2n$ times and then use
Gronwall's inequality to obtain
\bel{5.30}|\Th_x(\t,t,x)|\le
K\big(1+\|h(\t,\cd)\|_{C^1(\dbR^n)}\big),\qq\forall(t,x)\in[\t,T]\times\dbR^n.\ee
Now, we let $v^i(\cd\,,\cd)\in C([0,T];C^1(\dbR^n))$, $i=0,1$. Let
$\Th^i(\t,\cd\,,\cd)$ be the corresponding solutions of
(\ref{5.22}). Then
\bel{5.31}\ba{ll}
\ns\ds\Th^1(\t,t,x)-\Th^0(\t,t,x)=\int_t^T\3n\int_{\dbR^n}
\G^0(t,x;s,y)\big[\lan b(s,y,v^1_x(s,y)),\Th^1_x(\t,s,y)-\Th^0_x(\t,s,y)\ran\\
\ns\ds\qq\qq\qq\qq\qq\qq\q+\lan
b(s,y,v_x^1(s,y))-b(s,y,v_x^0(s,y)),\Th^0_x(\t,s,y)\ran\\
\ns\ds\qq\qq\qq\qq\qq\qq\q+g(\t,s,y,v^1_x(s,y))-g(\t,s,y,v^0(s,y))\big]dyds,\q(t,x)\in[\t,T]\times
\dbR^n,\ea\ee
and
\bel{}\ba{ll}
\ns\ds\Th^1_x(\t,t,x)-\Th^0_x(\t,t,x)=\int_t^T\3n\int_{\dbR^n}
\G^0_x(t,x;s,y)\big[\lan\Th^1_x(\t,s,y)-\Th^0_x(\t,s,y),b(s,y,v^1_x(s,y))\ran\\
\ns\ds\qq\qq\qq\qq\qq\qq\q+\lan\Th^0_x(\t,s,y),b(s,y,v_x^1(s,y))-b(s,y,v_x^0(s,y))
\ran\\
\ns\ds\qq\qq\qq\qq\qq\qq\q+g(\t,s,y,v^1_x(s,y))-g(\t,s,y,v^0(s,y))\big]dyds,\q(t,x)\in[\t,T]\times
\dbR^n,\ea\ee
Hence,
\bel{5.33}\ba{ll}
\ns\ds|\Th^1_x(\t,t,x)-\Th^0_x(\t,t,x)|\le\int_t^T\3n\int_{\dbR^n}
{Ke^{-{\l|x-y|^2\over4(s-t)}}\over(s-t)^{n+1\over2}}\(|\Th^1_x(\t,s,y)
-\Th^0_x(\t,s,y)|\\
\ns\ds\qq\qq\qq\qq\qq\qq\qq\qq+\big(|\Th^0_x(\t,s,y)|+1\big)\,|v_x^1(s,y)-v_x^0(s,y)|\)dyds\\
\ns\ds\qq\qq\qq\le\int_t^T\3n\int_{\dbR^n}
{Ke^{-{\l|x-y|^2\over4(s-t)}}\over(s-t)^{n+1\over2}}|\Th^1_x(\t,s,y)
-\Th^0_x(\t,s,y)|dyds\\
\ns\ds\qq\qq\qq\qq+K(T-t)^{1\over2}
\big(1+\|h(\t,\cd)\|_{C^1(\dbR^n)}\big)\|v^1_x(\cd\,,\cd)-v^0_x(\cd\,,\cd)\|_{
C([\t,T]\times\dbR^n)}.\ea\ee
Then, similar to (\ref{5.30}) obtained from (\ref{5.28}) via
(\ref{5.29}), we can obtain
\bel{}|\Th^1_x(\t,t,x)-\Th^0_x(\t,t,x)|\le K(T-t)^{1\over2}
\big(1+\|h(\t,\cd)\|_{C^1(\dbR^n)}\big)\|v^1_x(\cd\,,\cd)-v^0_x(\cd\,,\cd)\|_{
C([\t,T]\times\dbR^n)}.\ee
On the other hand, from (\ref{5.31}), we have
\bel{}\ba{ll}
\ns\ds|\Th^1(\t,t,x)-\Th^0(\t,t,x)|\le\int_t^T\3n\int_{\dbR^n}
{Ke^{-{\l|x-y|^2\over4(s-t)}}\over(s-t)^{n\over2}}\(|\Th^1_x(\t,s,y)-
\Th^0_x(\t,s,y)|\\
\ns\ds\qq\qq\qq\qq\qq\qq\qq\qq\q+\big(|\Th^0_x(\t,s,y)|+1\big)|v^1_x(s,y)-v^0_x(s,y)|\)dyds.
\ea\ee
Comparing the above with (\ref{5.33}), we see that the following
must be true:
\bel{}|\Th^1(\t,t,x)-\Th^0(\t,t,x)|\le K(T-t)
\big(1+\|h(\t,\cd)\|_{C^1(\dbR^n)}\big)\|v^1_x(\cd\,,\cd)-v^0_x(\cd\,,\cd)\|_{
C([\t,T]\times\dbR^n)}.\ee
Hence, we obtain
\bel{}\ba{ll}
\ns\ds\|\Th^1(\t,\cd\,,\cd)-\Th^0(\t,\cd\,,\cd)\|_{C([\t,T];C^1(\dbR^n))}\\
\ns\ds\le
K(T-t)^{1\over2}\big(1+\|h(\t,\cd)\|_{C^1(\dbR^n)}\big)\|v^1(\cd\,,\cd)-
v^0(\cd\,,\cd)\|_{C([\t,T];C^1(\dbR^n))}.\ea\ee
From the above procedure, we see that $K>0$ in the above is an
absolute constant, independent of $(\t,t)\in D[0,T]$. Hence, in
particular, we have (denoting $V^i(t,x)=\Th^i(t,t,x)$)
\bel{}\ba{ll}
\ns\ds\|V^1(\cd\,,\cd)-V^0(\cd\,,\cd)\|_{C([T-\d,T];C^1(\dbR^n))}\\
\ns\ds\le
K\d^{1\over2}\big(1+\|h(\cd,\cd)\|_{B([0,T];C^1(\dbR^n))}\big)\|v^1(\cd\,,\cd)-
v^0(\cd\,,\cd)\|_{C([T-\d,T];C^1(\dbR^n))}.\ea\ee
Clearly, by choosing $\d>0$ small, we get a contraction mapping
$v(\cd\,,\cd)\mapsto V(\cd\,,\cd)$ on $C([T-\d,T];C^1(\dbR^n))$.
Therefore, this map admits a unique fixed point. Since we may obtain
similar estimates on $[T-2\d,T-\d])$, etc., one sees that the fixed
point will exists on the whole space $C([0,T;C^1(\dbR^n))$ for the
map $v(\cd\,,\cd)\mapsto V(\cd\,,\cd)$. Then we obtain the
well-posedness of the following:
\bel{}\ba{ll}
\ns\ds\Th(\t,t,x)=\int_{\dbR^n}\G^0(t,x;T,y)h(\t,y)dy+\int_t^T\3n\int_{\dbR^n}
\G^0(t,x;s,y)\big[\lan b(s,y,\Th_x(s,s,y)),\Th_x(\t,s,y)\ran\\
\ns\ds\qq\qq\qq\qq\qq\qq\qq\qq+g(\t,s,y,\Th_x(s,s,y))\big]dyds,\qq(\t,t,x)\in
D[0,T]\times\dbR^n,\ea\ee
Finally, by the regularity of the above expression, we know that
$\Th(\t,t,x)$ is $C^{2+\a}$ in $x$, $C^{1+{\a\over2}}$ in $t$ for
some $\a\in(0,1)$, and the PDE (\ref{5.10}) is satisfied. \endpf

\ms

Recall that in Section 4, by assuming the convergence of
$\Th^\Pi(\cd\,,\cd\,,\cd)$ and $V^\Pi(\cd\,,\cd)$, we get the
equilibrium HJB equation for $\Th(\cd\,,\cd\,,\cd)$ and then
$V(\cd\,,\cd)$ is characterized by an equilibrium HJB integral
equation. We now want to show that under conditions ensuring (P), we
do have the expected convergence. This will make our whole procedure
satisfactorily complete for certain cases, at least. For the sake of
simplicity, we assume that all the involved functions are bounded
and continuously differentiable up to a needed order with bounded
derivatives.

\ms

When (\ref{si=si}) holds, for $k=0,1,\cds, N-1$, we have
\bel{Th(k+1)s}\left\{\ba{ll}
\ns\ds\Th^\Pi_t(\t,t,x)+\tr\big[a(t,x)\Th^\Pi_{xx}(\t,t,x)\big]
+\lan b\big(t,x,\psi(\ell^\Pi(t),t,x,V^\Pi_x(t,x))\big),\Th^\Pi_x(\t,t,x)\ran\\
\ns\ds\qq\qq+g^\Pi(\t,t,x,\psi(\ell^\Pi(t),t,x,V^\Pi_x(t,x))\big)=0,
\qq(t,x)\in[t_k,T]\times\dbR^n,\\
\ns\ds\Th^\Pi(\t,T,x)=h^\Pi(\t,x),\qq\qq x\in\dbR^n,\ea\right.\ee
with
$$\left\{\ba{ll}
\ns\ds
h^\Pi(\t,x)=\sum_{k=1}^{N-1}h(t_{k-1},x)I_{[t_{k-1},t_k)}(\t),\qq(\t,x)\in[0,T]\times\dbR^n,\\
\ns\ds
g^\Pi(\t,s,x,u)=\sum_{k=1}^{N-1}g(t_{k-1},s,x,u)I_{[t_{k-1},t_k)}(\t),\qq(\t,s,x,u)\in
D[0,T]\times\dbR^n\times U.\ea\right.$$
Since (\ref{Th(k+1)s}) is well-posed, we have
\bel{4.49}V^\Pi(t,x)=\Th^\Pi(\ell^\Pi(t),t,x),\qq(t,x)\in[0,t_{N-1})\times\dbR^n.\ee
Also, by the assumed uniform Lipschitz continuity of
$\t\mapsto(h(\t,x),h_y(\t,y),g(\t,t,x,u))$, we have
\bel{}\ba{ll}
\ns\ds|h^\Pi(\t,x)-h(\t,x)|+|h^\Pi_x(\t,x)-h_x(\t,x)|+|g^\Pi(\t,s,x,u)-g(\t,s,x,u)|\le
K\|\Pi\|,\\
\ns\ds\qq\qq\qq\qq\qq\qq\qq\forall(\t,s,x,u)\in
D[0,T]\times\dbR^n\times U.\ea\ee
Next, by Proposition 5.1, we have
\bel{}\ba{ll}
\ns\ds\Th^\Pi(\t,t,x)-\Th(\t,t,x)=\int_{\dbR^n}\G^0(t,x;T,y)\big[h^\Pi(\t,y)-h(\t,y)\big]dy\\
\ns\ds\qq+\int_t^T\3n\int_{\dbR^n} \G^0(t,x;s,y)\big[\lan
b\big(s,y,\psi(\ell^\Pi(s),s,y,V^\Pi_x(s,y))\big),\Th^\Pi_x(\t,s,y)\ran\\
\ns\ds\qq\qq\qq-\lan b\big(s,y,\psi(s,s,y,V_x(s,y))\big),\Th_x(\t,s,y)\ran\\
\ns\ds\qq\qq\qq+g^\Pi\big(\t,s,y,\psi(\ell^\Pi(s),x,V^\Pi_x(s,y))\big)
-g\big(\t,s,y,\psi(s,x,V_x(s,y))\big)
\big]dyds\\
\ns\ds=\int_{\dbR^n}\G^0(t,x;T,y)\big[h^\Pi(\t,y)-h(\t,y)\big]dy\\
\ns\ds\qq+\int_t^T\3n\int_{\dbR^n} \G^0(t,x;s,y)\big[\lan
b\big(s,y,\psi(\ell^\Pi(s),s,y,V^\Pi_x(s,y))\big),\Th^\Pi_x(\t,s,y)-\Th_x(\t,s,y)
\ran\\
\ns\ds\qq\qq\qq+\lan b\big(s,y,\psi(\ell^\Pi(s),s,y,V^\Pi_x(s,y))\big)
-b\big(s,y,\psi(s,s,y,V_x(s,y))\big),\Th_x(\t,s,y)\ran\\
\ns\ds\qq\qq\qq+g^\Pi\big(\t,s,y,\psi(\ell^\Pi(s),x,V^\Pi_x(s,y))\big)
-g\big(\t,s,y,\psi(s,x,V_x(s,y))\big) \big]dyds.\ea\ee
Consequently,
\bel{}\ba{ll}
\ns\ds\Th_x^\Pi(\t,t,x)-\Th_x(\t,t,x)=\int_{\dbR^n}\G^0(t,x;T,y)\big[h^\Pi_y(\t,y)-h_y(\t,y)\big]dy\\
\ns\ds\qq-\int_{\dbR^n}\G^0(t,x;T,y)\rho(t,x,T,y)\big[h^\Pi(\t,y)-h(\t,y)\big]dy\\
\ns\ds\qq+\int_t^T\3n\int_{\dbR^n} \G^0_x(t,x;s,y)\big[\lan
b\big(s,y,\psi(\ell^\Pi(s),s,y,V^\Pi_x(s,y))\big),\Th^\Pi_x(\t,s,y)-\Th_x(\t,s,y)
\ran\\
\ns\ds\qq\qq\qq+\lan b\big(s,y,\psi(\ell^\Pi(s),s,y,V^\Pi_x(s,y))\big)
-b\big(s,y,\psi(s,s,y,V_x(s,y))\big),\Th_x(\t,s,y)\ran\\
\ns\ds\qq\qq\qq+g^\Pi\big(\t,s,y,\psi(\ell^\Pi(s),s,x,V^\Pi_x(s,y))\big)
-g\big(\t,s,y,\psi(s,x,V_x(s,y))\big) \big]dyds.\ea\ee
Thus,
\bel{}\ba{ll}
\ns\ds|\Th_x^\Pi(\t,t,x)-\Th_x(\t,t,x)|\le\int_{\dbR^n}{Ke^{-{\l|x-y|^2\over4(T-t)}}\over
(T-t)^{n\over2}}|h^\Pi_y(\t,y)-h_y(\t,y)|dy\\
\ns\ds\qq+\int_{\dbR^n}{Ke^{-{\l|x-y|^2\over4(T-t)}}\over
(T-t)^{n\over2}}\(1+{|x-y|^2\over(T-t)}\)|h^\Pi(\t,y)-h(\t,y)|dy\\
\ns\ds\qq+\int_t^T\3n\int_{\dbR^n}
{Ke^{-{\l|x-y|^2\over4(T-t)}}\over(T-t)^{n+1\over2}}\(|\Th^\Pi_x(\t,s,y)-\Th_x(\t,s,y)|
+|\ell^\Pi(s)-s|+|V^\Pi_x(s,y)-V_x(s,y)|\\
\ns\ds\qq\qq\qq+|g^\Pi\big(\t,s,y,\psi(\t,s,x,V_x(s,y))\big)
-g\big(\t,s,y,\psi(s,x,V_x(s,y))\big)|\)dyds\\
\ns\ds\le\int_{\dbR^n}{Ke^{-{\l|x-y|^2\over4(T-t)}}\over
(T-t)^{n\over2}}\|\Pi\|dy+\int_{\dbR^n}{Ke^{-{\l|x-y|^2\over4(T-t)}}\over
(T-t)^{n\over2}}\(1+{|x-y|^2\over(T-t)}\)\|\Pi\|dy\\
\ns\ds\qq+\int_t^T\3n\int_{\dbR^n}
{Ke^{-{\l|x-y|^2\over4(T-t)}}\over(T-t)^{n+1\over2}}\(|\Th^\Pi_x(\t,s,y)-\Th_x(\t,s,y)|
+\|\Pi\|+|\Th^\Pi_x(\ell^\Pi(s),s,y)-\Th_x(s,s,y)|\)dyds.\ea\ee
Therefore,
\bel{}\ba{ll}
\ns\ds\sup_{\t\in[0,t]}|\Th_x^\Pi(\t,t,x)-\Th_x(\t,t,x)|\le
K\|\Pi\|\1n+\2n\int_t^T\3n\int_{\dbR^n}
{Ke^{-{\l|x-y|^2\over4(T-t)}}\over(T-t)^{n+1\over2}}|\Th_x(\ell^\Pi(s),s,x)-\Th_x(s,s,y)|
dyds\\
\ns\ds\qq\qq\qq\qq\qq\qq\qq+\int_t^T\3n\int_{\dbR^n}
{Ke^{-{\l|x-y|^2\over4(T-t)}}\over(T-t)^{n+1\over2}}\sup_{\t\in[0,s]}
|\Th^\Pi_x(\t,s,y)-\Th_x(\t,s,y)|dyds.\ea\ee
This yields that
\bel{}\sup_{\t\in[0,t]}|\Th_x^\Pi(\t,t,x)-\Th_x(\t,t,x)|\le
K\(\|\Pi\|+\2n\int_t^T\3n\int_{\dbR^n}
{Ke^{-{\l|x-y|^2\over4(T-t)}}\over(T-t)^{n+1\over2}}|\Th_x(\ell^\Pi(s),s,x)-\Th_x(s,s,y)|
dyds\).\ee
Likewise, we can have a similar estimate for
$|\Th^\Pi(\t,t,x)-\Th(\t,t,x)|$. Hence, we obtain
\bel{}\ba{ll}
\ns\ds\sup_{(\t,t)\in
D[0,T]}\|\Th^\Pi(\t,t,\cd)-\Th(\t,t,\cd)\|_{C^1(\dbR^n)}\\
\ns\ds\le K\(\|\Pi\|+\2n\int_t^T\3n\int_{\dbR^n}
{Ke^{-{\l|x-y|^2\over4(T-t)}}\over(T-t)^{n+1\over2}}|\Th_x(\ell^\Pi(s),s,x)-\Th_x(s,s,y)|
dyds\).\ea\ee
From this, our expected convergence follows.

\section{Some Special Cases}

In this section, we are going to look at several important special
cases. We will mainly look at the corresponding forms of our
equilibrium HJB equations.

\subsection{A Linear-Quadratic Problem}

Let us look at the LQ problem. For any initial pair
$(t,x)\in[0,T)\times\dbR^n$, the state equation is
\bel{}\left\{\ba{ll}
\ns\ds
dX(s)=\big[A(s)X(s)+B(s)u(s)\big]ds+\big[A_1(s)X(s)+B_1(s)u(s)\big]dW(s),\qq
s\in[t,T],\\
\ns\ds X(t)=x,\ea\right.\ee
with the cost functional
\bel{}\ba{ll}
\ns\ds J(t,x;u(\cd))={1\over2}\dbE_t\[\int_t^T\big(\lan
Q(t,s)X(s),X(s)\ran+\lan R(t,s)u(s),u(s)\ran\big)ds+\lan
G(t)X(T),X(T)\ran\].\ea\ee
Then
$$\ba{ll}
\ns\ds\dbH(\t,t,x,u,p,P)=\lan
p,A(t)x+B(t)u\ran+{1\over2}\tr\[\big(A_1(t)x+B_1(t)u\big)
\big(A_1(t)x+B_1(t)u\big)^TP\]\\
\ns\ds\qq\qq\qq\qq\qq+{1\over2}\[\lan Q(\t,t)x,x\ran+\lan R(\t,t)u,u\ran\]\\
\ns\ds=\lan p,A(t)x\ran+{1\over2}\lan
\big[A_1(t)^TPA_1(t)+Q(\t,t)\big]x,x\ran\\
\ns\ds\qq\qq+{1\over2}\lan\big[R(\t,t)+B_1(t)^TPB_1(t)\big]u,u\ran+\lan
u,B(t)^Tp+B_1(t)^TPA_1(t)x\ran.\ea$$
This yields
$$\psi(\t,t,x,p,P)=-\big[R(\t,t)+B_1(t)^TPB_1(t)\big]^{-1}\big[B(t)^Tp
+B_1(t)^TPA_1(t)x\big],$$
and
$$\ba{ll}
\ns\ds\dbH(\t,t,x,\psi(t,t,x,\bar p,\bar P),p,P)\\
\ns\ds=\lan p,A(t)x\ran+{1\over2}\lan
\big[A_1(t)^TPA_1(t)+Q(\t,t)\big]x,x\ran\\
\ns\ds\qq+{1\over2}\lan\big[R(\t,t)+B_1(t)^TPB_1(t)\big]\psi(t,t,x,\bar
p,\bar P),\psi(t,t,x,\bar p,\bar P)\ran\\
\ns\ds\qq+\lan\psi(t,t,x,\bar p,\bar
P),B(t)^Tp+B_1(t)^TPA_1(t)x\ran\\
\ns\ds=\lan p,A(t)x\ran+{1\over2}\lan\big[A_1(t)^TPA_1(t)+Q(\t,t)\big]x,x\ran\\
\ns\ds\q+{1\over2}\lan\big[R(\t,t)+B_1(t)^TPB_1(t)\big]\big[R(t,t)+B_1(t)^T\bar
PB_1(t)\big]^{-1}
\big[B(t)^T\bar p+B_1(t)^T\bar PA_1(t)x\big],\\
\ns\ds\qq\qq\big[R(t,t)+B_1(t)^T\bar
PB_1(t)\big]^{-1}\big[B(t)^T\bar p+B_1(t)^T
\bar PA_1(t)x\big]\ran\\
\ns\ds\q-\lan\big[R(t,t)+B_1(t)^T\bar
PB_1(t)\big]^{-1}\big[B(t)^T\bar p+B_1(t)^T\bar
PA_1(t)x\big],B(t)^Tp+B_1(t)^TPA_1(t)x\ran.\ea$$
Hence, the equilibrium HJB equation takes the following form:
\bel{HJB4}\left\{\ba{ll}
\ns\ds\Th_t(\t,t,x)+\lan\Th_x(\t,t,x),A(t)x\ran+{1\over2}
\lan\big[A_1(t)^T\Th_{xx}(\t,t,x)A_1(t)+Q(\t,t)\big]x,x\ran\\
\ns\ds+{1\over2}\lan\big[R(\t,t)+B_1(t)^T\Th_{xx}(\t,t,x)B_1(t)\big]
\big[R(t,t)+B_1(t)^T\Th_{xx}(t,t,x)B_1(t)\big]^{-1}\cr
\ns\ds\qq\cd
\big[B(t)^T\Th_x(t,t,x)+B_1(t)^T\Th_{xx}(t,t,x)A_1(t)x\big],\\
\ns\ds\qq\big[R(t,t)+B_1(t)^T\Th_{xx}(t,t,x)B_1(t)\big]^{-1}
\big[B(t)^T\Th_x(t,t,x)+B_1(t)^T
\Th_{xx}(t,t,x)A_1(t)x\big]\ran\\
\ns\ds-\lan\big[R(t,t)+B_1(t)^T\Th_{xx}(t,t,x)B_1(t)\big]^{-1}\big[B(t)^T\Th_x(t,t,x)
+B_1(t)^T\Th_{xx}(t,t,x)A_1(t)x\big],\\
\ns\ds\qq B(t)^T\Th_x(\t,t,x)+B_1(t)^T\Th_{xx}(\t,t,x)
A_1(t)x\ran=0,\qq(\t,t,x)\in D[0,T]\times\dbR^n,\\
\ns\ds\Th(\t,T,x)={1\over2}\lan G(\t)x,x\ran,\qq
(\t,x)\in[0,T]\times\dbR^n.\ea\right.\ee
Although having a little bit complicated looking, the above has a
quadratic structure which can help us to study the well-posedness of
it. To see that, let
\bel{}\Th(\t,t,x)={1\over2}\lan P(\t,t)x,x\ran,\qq(\t,t,x)\in
D[0,T]\times\dbR^n,\ee
with some undetermined map $P:D[0,T]\to\cS^n$. Plugging the above
into (\ref{HJB4}), we see that the map $P(\cd\,,\cd)$ should satisfy
the following equation:
$$\ba{ll}
\ns\ds0=P_t(\t,t)+P(\t,t)A(t)+A(t)^TP(\t,t)+A_1(t)^TP(\t,t)A_1(t)+Q(\t,t)\\
\ns\ds\q+[P(t,t)B(t)+A_1(t)^TP(t,t)B_1(t)][R(t,t)+B_1(t)^TP(t,t)B_1(t)]^{-1}\\
\ns\ds\q\cd[R(\t,t)\1n+\1n B_1(t)^T\1n P(\t,t)B_1(t)][R(t,t)\1n+\1n
B_1(t)^T\1n P(t,t)B_1(t)]^{-1}[B(t)^T\1n P(t,t)\1n+\1n B_1(t)^T\1n
P(t,t)A_1(t)]\\
\ns\ds\q-[P(t,t)B(t)+A_1(t)^TP(t,t)B_1(t)][R(t,t)+B_1(t)^TP(t,t)B_1(t)]^{-1}[B(t)^TP(\t,t)+B_1(t)^TP(\t,t)A_1(t)]\\
\ns\ds\q-[P(\t,t)B(t)+A_1(t)^TP(\t,t)B_1(t)][R(t,t)+B_1(t)^TP(t,t)B_1(t)]^{-1}[B(t)^TP(t,t)+B_1(t)^TP(t,t)A_1(t)],\ea$$
with the terminal condition
$$P(\t,T)=G(\t).$$
Note that if $P(t,t)$ and $R(t,t)$ are replaced by $P(\t,t)$ and
$R(\t,t)$, respectively, the above becomes a standard Riccati
equation with a parameter $\t$. The appearance of $P(t,t)$ and
$R(t,t)$ makes the above non-standard. We may rewrite the above as
follows (suppressing $t$ in $P(\t,t)$, etc., for simplicity):
$$\ba{ll}
\ns\ds0=P_t(\t)+P(\t)A+A^TP(\t)+A_1^TP(\t)A_1+Q(\t)
+\big[P(t)B+A_1^TP(t)B_1\big]\big[R(t)+B_1^TP(t)B_1\big]^{-1}\\
\ns\ds\q\cd\big[R(\t)+B_1^TP(\t)B_1\big]\big[R(t)+B_1^TP(t)B_1\big]^{-1}
\big[B^TP(t)+B_1^TP(t)A_1\big]\\
\ns\ds\q-\big[P(t)B+A_1^TP(t)B_1\big]\big[R(t)+B_1^TP(t)B_1\big]^{-1}
\big[B^TP(\t)+B_1^TP(\t)A_1\big]\\
\ns\ds\q-\big[P(\t)B+A_1^TP(\t)B_1\big]\big[R(t)+B_1^TP(t)B_1\big]^{-1}
\big[B^TP(t)+B_1^TP(t)A_1\big]\\
\ns\ds\q=P_t(\t)+P(\t)\(A-B\big[R(t)+B_1^TP(t)B_1\big]^{-1}\big[B^TP(t)
+B_1^TP(t)A_1\big]\)\\
\ns\ds\qq+\(A-B\big[R(t)+B_1^TP(t)B_1\big]^{-1}\big[B^TP(t)+B_1^TP(t)A_1\big]
\)^TP(\t)+Q(\t)\\
\ns\ds\qq+\(A_1-B_1\big[R(t)+B_1^TP(t)B_1\big]^{-1}\big[B^TP(t)+B_1^TP(t)A_1
\big]\)^TP(\t)\\
\ns\ds\qq\cd\(A_1-B_1\big[R(t)+B_1^TP(t)B_1\big]^{-1}
\big[B^TP(t)+B_1^TP(t)A_1\big]\)\\
\ns\ds\qq+\big[P(t)B\1n+\1n A_1^TP(t)B_1\big]\big[R(t)\1n+\1n
B_1^TP(t)B_1\big]^{-1}\1n R(\t)\big[R(t)\1n+\1n
B_1^TP(t)B_1\big]^{-1}\1n\big[B^TP(t)\1n+\1n B_1^TP(t)A_1\big].\ea$$
Denote
\bel{}\left\{\2n\ba{ll}
\ns\ds\G(t)=\big[R(t,t)+B_1(t)^TP(t,t)B_1(t)\big]^{-1}\big[B(t)^TP(t,t)
+B_1(t)^TP(t,t)A_1(t)\big],\\
\ns\ds\h A(t)=A(t)-B(t)\G(t),\qq\h A_1(t)=A_1(t)-B_1(t)\G(t),\\
\ns\ds\h Q(\t,t)=Q(\t,t)+\G(t)^TR(\t,t)\G(t),\ea\right.\q(\t,t)\in
D[0,T].\ee
Then the equation for $P(\cd\,,\cd)$ can be written as follows:
\bel{P(tau,t)}\left\{\ba{ll}
\ns\ds P_t(\t,t)+P(\t,t)\h A(t)+\h A(t)^TP(\t,t)+\h
A_1(t)^TP(\t,t)\h A_1(t)+\h Q(\t,t)=0,\qq(\t,t)\in D[0,T],\\
\ns\ds P(\t,T)=G(\t),\qq\qq\t\in[0,T].\ea\right.\ee
Next, let $\F(\cd\,,\cd)$ be the fundamental matrix of $\big(\h
A(\cd),\h A_1(\cd)\big)$, i.e., the following holds:
\bel{}\left\{\ba{ll}
\ns\ds d\F(s,t)=\h A(s)\F(s,t)ds+\h A_1(s)\F(s,t)dW(s),\qq s\in[t,T],\\
\ns\ds\F(t,t)=I.\ea\right.\ee
Applying It\^o's formula to $s\mapsto\lan
P(\t,s)\F(s,t)x,\F(s,t)x\ran$ on $[t,T]$, we have
$$\ba{ll}
\ns\ds\lan\F(T,t)^TG(\t)\F(T,t)x,x\ran-\lan
P(\t,t)x,x\ran=\int_t^T-\lan\h Q(\t,s)\F(s,t)x,\F(s,t)x\ran
ds\\
\ns\ds\qq\qq\qq\qq+\int_t^T\lan\big[P(\t,s)\h A_1(s)+\h
A_1(s)^TP(\t,s)\big]\F(s,t)x,\F(s,t)x\ran dW(s),\ea$$
which leads to
\bel{P(tau,t)}\ba{ll}
\ns\ds P(\t,t)=\F(T,t)^TG(\t)\F(T,t)+\int_t^T\F(s,t)^T\h
Q(\t,s)\F(s,t)ds\\
\ns\ds\qq\qq\qq-\int_t^T\F(s,t)^T\big[P(\t,s)\h A_1(s)+\h
A_1(s)^TP(\t,s)\big]\F(s,t)dW(s),\q(\t,t)\in D[0,T].\ea\ee
Note that although $P(\t,t)$ is a deterministic function, the above
representation is stochastic. From the above, we have
\bel{P(tau,t)1}P(\t,t)\1n=\1n\dbE_t\[\F(T,t)^T\1n
G(\t)\F(T,t)\1n+\2n\int_t^T\3n\F(s,t)^T\1n\big[
Q(\t,s)\1n+\1n\G(s)^T\1n R(\t,s)\G(s)\big]\F(s,t)ds\],\q(\t,t)\in
D[0,T].\ee
In particular, taking $\t=t$ and denoting $P(t)=P(t,t)$, one has
\bel{}\ba{ll}
\ns\ds
P(t)=\dbE_t\[\F(T,t)^TG(t)\F(T,t)\1n+\2n\int_t^T\2n\F(s,t)^T\big[Q(t,s)\1n
+\1n\G(s)^TR(t,s)\G(s)\big]\F(s,t)ds\],\qq t\in[0,T].\ea\ee
Combining the above, we end up with the following system for the
function $P(\cd)$:
\bel{6.9}\3n\left\{\2n\ba{ll}
\ns\ds
P(t)=\dbE_t\[\F(T,t)^TG(t)\F(T,t)\1n+\2n\int_t^T\2n\F(s,t)^T\big[Q(t,s)\1n+\1n\G(s)^TR(t,s)\G(s)\big]
\F(s,t)ds\],\qq t\in[0,T],\\
\ns\ds\F(s,t)\1n=\1n I\1n+\2n\int_t^s\2n\big[A(r)\1n-\1n B(r)\G(r)\big]\F(r,t)dr\1n
+\2n\int_t^s\2n\big[A_1(r)\1n-\1n B_1(r)\G(r)\big]\F(r,t)dW(r),\q(t,s)\in D[0,T],\\
\ns\ds
\G(t)=\big[R(t,t)+B_1(t)^TP(t)B_1(t)\big]^{-1}\big[B(t)^TP(t)+B_1(t)^T
P(t)A_1(t)\big],\qq t\in[0,T].\ea\right.\ee
We refer to the above as a {\it Riccati-Volterra integral equation
system} for the corresponding (time-inconsistent) LQ problem. Note
that the above is actually a coupled forward-backward stochastic
Volterra integral equation system (FBSVIE, for short). Some relevant
results concerning backward stochastic Volterra integral equations
(BSVIEs, for short) can be found in \cite{Yong 2006, Yong 2008}. If
$\big(\F(\cd\,,\cd),P(\cd)\big)$ is a solution to the above, then
the time-consistent equilibrium control is given by
\bel{6.13}\bar u(t)=-\G(t)\bar X(t),\q t\in[0,T].\ee
In the case that
\bel{A_1=0}A_1(\cd)=0,\qq B_1(\cd)=0,\ee
the above (\ref{6.9}) is reduced to the following Riccati-Volterra
integral equation system (for a deterministic time-inconsistent LQ
problem):
\bel{RV}\left\{\ba{ll}
\ns\ds
P(t)=\F(T,t)^TG(t)\F(T,t)+\int_t^T\F(s,t)^T\big[Q(t,s)+\G(s)^TR(t,s)\G(s)\big]
\F(s,t)ds,\q t\in[0,T],\\
\ns\ds\F(s,t)=I+\int_t^s\big[A(r)-B(r)\G(r)\big]\F(r,t)dr,\qq
(t,s)\in D[0,T],\\
\ns\ds\G(t)=R(t,t)^{-1}B(t)^TP(t),\qq t\in[0,T],\ea\right.\ee
and the time-consistent equilibrium control is given by (\ref{6.13})
with a simpler $\G(\cd)$. This recovers the case presented in
\cite{Yong 2011} where the well-posedness of (\ref{RV}) was
established.

\ms

For (\ref{6.9}), we have the following result.

\ms

\bf Proposition 6.1. \sl Suppose
\bel{6.15}\left\{\ba{ll}
\ns\ds A(\cd),A_1(\cd)\in C([0,T];\dbR^{n\times n}),\q
B(\cd),B_1(\cd)\in C([0,T];\dbR^{n\times m}),\\
\ns\ds Q(\cd)\in C(D[0,T];\cl{\cS^n_+}\,),\q R(\cd)\in
C(D[0,T];\cS^m_+),\q G(\cd)\in C([0,T];\cl{\cS^n_+}\,).\ea\right.\ee
Further, suppose
\bel{6.16}\sup_{P\in\cS^n_+}\Big|\big[R(t,t)+B_1(t)^TPB_1(t)\big]^{-1}
\big[B(t)^TP+B_1(t)^TPA_1(t)\big]\Big|\equiv L<\infty.\ee
Then $(\ref{6.9})$ admits a unique solution.

\ms

\it Proof. \rm Let $\cX[\t,T]\equiv
C\big([\t,T];\cl{\cS^n_+}\,\big)$ which is a complete metric space
with the metric induced by the norm in $C([\t,T];\cS^n)$. For any
$p(\cd)\in\cX[0,T]$, define
$$\G(t;p(\cd))=\big[R(t,t)+B_1(t)^Tp(t)B_1(t)\big]^{-1}\big[B(t)^Tp(t)+B_1(t)^Tp(t)A_1(t)\big],
\qq t\in[0,T].$$
By (\ref{6.16}), we have
$$|\G(t;p(\cd))|\le K,\qq t\in[0,T],$$
with the bound independent of $p(\cd)\in\cX$. Let
$\F(\cd\,,\cd)\equiv\F(\cd\,,\cd\,;p(\cd))$ be the solution to the
following:
$$\ba{ll}
\ns\ds\F(s,t)=I+\int_t^s\big[A(r)-B(r)\G(r;p(\cd))\big]\F(r,t)dr\\
\ns\ds\qq\qq\q+\int_t^s\big[A_1(r)-B_1(r)\G(r;p(\cd))\big]\F(r,t)dW(r),\qq(t,s)\in
D[0,T].\ea$$
We have
$$\dbE\[\sup_{(t,s)\in D[0,T]}|\F(s,t)|^2\]\le K,$$
for some absolute constant $K$. Next, we define
$$\ba{ll}
\ns\ds P(t)\equiv P(t;p(\cd))=\dbE_t\[\F(T,t)^TG(t)\F(T,t)\\
\ns\ds\qq\qq+\int_t^T\F(s,t)^T\big[Q(t,s)+\G(s;p(\cd))^TR(t,s)\G(s;p(\cd))\big]
\F(s,t)ds,\qq t\in[0,T].\ea$$
Clearly, $P(\cd)\in\cX[0,T]$. We want to show that $p(\cd)\mapsto
P(\cd\,;p(\cd))$ admits a unique fixed point. To this end, for any
$p^1(\cd),p^2(\cd)\in\cX[0,T]$, we denote
$$\G^i(\cd)=\G(\cd\,;p^i(\cd)),\q\F^i(\cd\,,\cd)=\F(\cd\,,\cd\,;p^i(\cd)),\q
P^i(\cd)=P(\cd\,;p^i(\cd)),\qq i=1,2.$$
Then (suppressing $t$)
$$\ba{ll}
\ns\ds|\G^1(t)-\G^2(t)|=\big|(R+B_1^Tp_1B_1)^{-1}
(B^Tp_1+B_1^Tp_1A_1)-(R+B_1^Tp_2B_1)^{-1}
(B^Tp_2+B_1^Tp_2A_1)|\\
\ns\ds\le\big|(R+B_1^Tp_1B_1)^{-1}\big|\,\big|B^T(p_1-p_2)+B_1^T(p_1-p_2)A_1\big|\\
\ns\ds\qq+\big|\big[(R+B_1^Tp_1B_1)^{-1}-(R+B_1^Tp_2B_1)^{-1}\big](B^Tp_2+B_1^Tp_2B_1)\big|\\
\ns\ds\le
K|p_1-p_2|+\big|(R+B_1^Tp_1B_1)^{-1}B_1^T(p_1-p_2)B_1(R+B_1^Tp_2B_1)^{-1}(B^Tp_2
+B_1^Tp_2A_1)\big|\\
\ns\ds\le K|p_1(t)-p_2(t)|.\ea$$
Next,
$$\ba{ll}
\ns\ds\F^1(s,t)-\F^2(s,t)=\int_t^s\big\{A(r)\big[\F^1(r,t)-\F^2(r,t)\big]
-B(r)\big[\G^1(r)-\G^2(r)\big]\F^1(r,t)\\
\ns\ds\qq\qq\qq\qq\qq\qq-B(r)\G^2(r)\big[\F^1(r,t)-\F^2(r,t)\big]\big\}dr\\
\ns\ds\qq\qq\qq\qq\qq+\int_t^s\big\{A_1(r)\big[\F^1(r,t)-\F^2(r,t)\big]
-B_1(r)\big[\G^1(r)-\G^2(r)\big]\F^1(r,t)\\
\ns\ds\qq\qq\qq\qq\qq\qq-B_1(r)\G^2(r)\big[\F^1(r,t)-\F^2(r,t)\big]\big\}dW(r).\ea$$
Then
$$\ba{ll}
\ns\ds\dbE\[\sup_{s\in[t,T]}|\F^1(s,t)-\F^2(s,t)|^2\]\le
K\int_t^s|\G^1(r)-\G^2(r)|^2dr\le K\int_t^s|p^1(r)-p^2(r)|^2dr.\ea$$
Consequently,
$$\ba{ll}
\ns\ds|P^1(t)-P^2(t)|\le\dbE\Big\{\big|\F^1(T,t)^TG(t)\F^1(T,t)-\F^2(T,t)G(t)\F^2(T,t)\big|\\
\ns\ds\qq\qq\qq\qq+\int_t^T\big|\F^1(s,t)\big[Q(t,s)+\G^1(s)R(t,s)\G^1(s)\big]\F^1(r,t)\\
\ns\ds\qq\qq\qq\qq\qq-\F^2(s,t)\big[Q(t,s)+\G^2(s)R(t,s)\G^2(s)\big]\F^2(r,t)\big|dr
\Big\}\\
\ns\ds\le
K\dbE\Big\{\(|\F^1(T,t)|+|\F^2(T,t)|\)|\F^1(T,t)-\F^2(T,t)|\\
\ns\ds\qq+\int_t^T\[\(|\F^1(s,t)|+|\F^2(s,t)|\)|\F^1(s,t)-\F^2(s,t)|\\
\ns\ds\qq\qq+\(|\F^1(s,t)\G^1(s)|+|\F^2(s,t)\G^2(s)|\)
|\F^1(s,t)\G^1(s)-\F^2(s,t)\G^2(s)|\]ds\Big\}\\
\ns\ds\le K\Big\{\(\dbE|\F^1(T,t)-\F^2(T,t)|^2\)^{1\over2}+\int_t^T\[\(\dbE|\F^1(s,t)-\F^2(s,t)|^2\)^{1\over2}\\
\ns\ds\qq\qq+\(\dbE|\F^1(s,t)\G^1(s)-\F^2(s,t)\G^2(s)|^2\)^{1\over2}\]ds\Big\}\\
\ns\ds\le K\Big\{\int_t^T|p^1(s)-p^2(s)|^2ds
+\int_t^T\[\(\int_t^s|p^1(r)-p^2(r)|^2dr\)^{1\over2}\\
\ns\ds\qq+\(\dbE|\F^1(s,t)-\F^2(s,t)|^2+|\G^1(s)-\G^2(s)|^2\)^{1\over2}\]ds\Big\}
\le
K\(\int_t^T|p^1(s)-p^2(s)|^2ds\)^{1\over2},\ea$$
with $K>0$ being an absolute constant. Hence, we obtain
$$\max_{t\in[T-\k,T]}|P^1(t)-P^2(t)|\le
K(T-\t)^{1\over2}\max_{t\in[T-\t,T]}|p^1(t)-p^2(t)|.$$
Hence, the map $p(\cd)\mapsto P(\cd)$ is contractive on $\cX[\t,T]$
as long as $T-\t>0$ small. Then a usual argument applies to obtain a
unique fixed point of $p(\cd)\mapsto P(\cd)$ on $\cX[0,T]$. This
proves the well-posedness of (\ref{6.9}). \endpf

\ms

To conclude this section, let us make a remark on the condition
(\ref{6.16}). It is not hard to show that when $m=n$ and
$B_1(\cd)^{-1}$ exists and bounded, then (\ref{6.16}) holds.
Apparently, this is a restrictive condition. We hope that in our
future publications, such a condition can be removed.

\subsection{A generalized Merton's problem}

In this subsection, let us look at the Merton's portfolio problem
with general discounting. The results in this subsection is
comparable with some of the results in \cite{Ekeland-Mbodji-Pirvu
2012} and \cite{Marin-Solano-Navas 2010}. Let us recall
\bel{}\left\{\ba{ll}
\ns\ds dX(s)=\big[rX(s)+(\m-r)u(s)-c(s)\big]ds+\si u(s)dW(s),\qq
s\in[t,T],\\
\ns\ds X(t)=x,\ea\right.\ee
and
\bel{}J(t,x;u(\cd),c(\cd))=\dbE_t\[\int_t^T\n(t,s)c(s)^\b
ds+\rho(t)X(T)^\b\].\ee
Then
$$\ba{ll}
\ns\ds\dbH(t,s,x,u,c,p,P)=p[rx+(\m-r)u-c]+{1\over2}\si^2u^2P+\n(t,s)c^\b\\
\ns\ds\qq\qq\qq\qq=rxp+{\si^2P\over2}\[u^2+2{(\m-r)p\over\si^2P}u\]+\n(t,s)c^\b-pc.\ea$$
The maximum of $(u,c)\mapsto\dbH(t,s,x,u,c,p,P)$ is attained at
$$\bar u=-{(\m-r)p\over\si^2P},\qq\bar
c=\({\b\n(t,s)\over p}\)^{1\over1-\b}.$$
We denote
$$\psi(t,s,x,p,P)\equiv(\bar u,\bar c)=\(-{(\m-r)p\over\si^2P},\[{\b\n(t,s)\over
p}\]^{1\over1-\b}\).$$
Then
$$\psi(t,t,x,\bar p,\bar P)\equiv(\bar u,\bar c)=\(-{(\m-r)\bar p\over\si^2\bar P},
\[{\b\n(t,t)\over \bar p}\]^{1\over1-\b}\).$$
and
$$\ba{ll}
\ns\ds\dbH(\t,t,x,\psi(t,t,x,\bar p,\bar P),p,P)\\
\ns\ds=rxp+{\si^2P\over2}\[{(\m-r)^2\bar p^2\over\si^4\bar
P^2}-2{(\m-r)p\over\si^2P}{(\m-r)\bar p\over\si^2\bar P}\]
+\n(\t,t)\({\b\n(t,t)\over\bar
p}\)^{\b\over1-\b}-p\({\b\n(t,t)\over\bar p}\)^{1\over1-\b}\\
\ns\ds=rxp+{(\m-r)^2\bar pP\over2\si^2\bar P}\({\bar p\over\bar
P}-2{p\over P}\) +{[\b\n(t,t)]^{\b\over1-\b}\over\bar
p^{1\over1-\b}}[\n(\t,t)\bar p-\b\n(t,t)p].\ea$$
Hence, the equilibrium HJB equation reads
\bel{HJB2a}\left\{\ba{ll}
\ns\ds\Th_t(\t,t,x)+rx\Th_x(\t,t,x)+{(\m-r)^2\Th_x(t,t,x)\Th_{xx}(\t,t,x)
\over2\si^2\Th_{xx}(t,t,x)}\({\Th_x(t,t,x)\over
\Th_{xx}(t,t,x)}-2{\Th_x(\t,t,x)\over\Th_{xx}(\t,t,x)}\)\\
\ns\ds+{[\b\n(t,t)]^{\b\over1-\b}\over\Th_x(t,t,x)^{1\over1-\b}}
\big[\n(\t,t)\Th_x(t,t,x)-\b\n(t,t)\Th_x(\t,t,x)\big]=0,\qq(\t,t,x)\in
D[0,T)\times(0,\infty),\\
\ns\ds\Th(\t,t,0)=0,\qq\qq\qq(\t,t)\in D[0,T],\\
\ns\ds\Th(\t,T,x)=\rho(\t)x^\b,\qq\qq(\t,x)\in(0,\infty)\times[0,\infty).\ea\right.\ee
Let
$$\Th(\t,t,x)=\f(\t,t)x^\b,\qq(\t,t,x)\in D[0,T]\times\dbR^n.$$
Then
$$\ba{ll}
\ns\ds0=\f_t(\t,t)+r\b\f(\t,t)+{(\m-r)^2\b\over2\si^2(1-\b)}\f(\t,t)+{[\b\n(t,t)]^{\b\over1-\b}\over[\b\f(t,t)]^{1\over1-\b}}
\big[\n(\t,t)\b\f(t,t)-\b^2\n(t,t)\f(\t,t)\big]\\
\ns\ds\q=\f_t(\t,t)+r\b\f(\t,t)+{(\m-r)^2\b\over2\si^2(1-\b)}\f(\t,t)
+{\n(t,t)^{\b\over1-\b}\over\f(t,t)^{1\over1-\b}}
\big[\n(\t,t)\f(t,t)-\b\n(t,t)\f(\t,t)\big]\\
\ns\ds\q=\f_t(\t,t)+\b\[r+{(\m-r)^2\over2\si^2(1-\b)}-\({\n(t,t)\over
\f(t,t)}\)^{1\over1-\b}\]\f(\t,t)
+\({\n(t,t)\over\f(t,t)}\)^{\b\over1-\b}\n(\t,t)\\
\ns\ds\q\equiv\f_t(\t,t)+\[\l-\b\({\n(t,t)\over
\f(t,t)}\)^{1\over1-\b}\]\f(\t,t)
+\({\n(t,t)\over\f(t,t)}\)^{\b\over1-\b}\n(\t,t),\ea$$
and
$$\f(\t,T)=\rho(\t),$$
with $\l$ given by (\ref{l2}). Thus,
$$\ba{ll}
\ns\ds\f(\t,t)=e^{\l(T-t)
-\b\int_t^T\big({\n(s,s)\over\f(s,s)}\big)^{1\over1-\b}ds}\rho(\t)+\int_t^Te^{\l(s-t)
-\b\int_t^s\big({\n(s',s')\over\f(s',s')}\big)^{1\over1-\b}ds'}
\({\n(s,s)\over\f(s,s)}\)^{\b\over1-\b}\n(\t,s)ds,\\
\ns\ds\qq\qq\qq\qq\qq\qq\qq\qq\qq\qq\qq\qq(\t,t)\in D[0,T].\ea$$
Hence, we obtain the following integral equation for
$t\mapsto\f(t,t)$:
\bel{f(t,t)}\ba{ll}
\ns\ds\f(t,t)=e^{\l(T-t)-\b\int_t^T\big({\n(s,s)\over\f(s,s)}\big)^{1\over1-\b}ds}\rho(t)\\
\ns\ds\qq\qq\qq+\int_t^Te^{\l(s-t)-\b\int_t^s\big({\n(s',s')\over\f(s',s')}\big)^{1\over1-\b}ds'}
\({\n(s,s)\over\f(s,s)}\)^{\b\over1-\b}\n(t,s)ds,\qq
t\in[0,T].\ea\ee
Once such an integral equation admits a unique solution
$t\mapsto\f(t,t)$, we will obtain the time-consistent equilibrium
value function
$$V(t,x)=\Th(t,t,x)=\f(t,t)x^\b,\qq(t,x)\in[0,T]\times[0,\infty),$$
and the time-consistent equilibrium control
\bel{uc}\bar u(t)=-{\m-r\over\si^2(1-\b)}\bar X(t),\qq\bar
c(t)=\({\n(t,t)\over \f(t,t)}\)^{1\over1-\b}\bar X(t),\qq
t\in[0,T].\ee
To establish the well-posedness of (\ref{f(t,t)}), we set
$$z(t)={\f(t,t)\over\n(t,t)},\qq t\in[0,T].$$
Then (\ref{f(t,t)}) is equivalent to the following:
\bel{z}z(t)=e^{\l(T-t)-\b\int_t^Tz(s)^{1\over\b-1}ds}\rho(t)
+\int_t^Te^{\l(\t-t)-\b\int_t^\t
z(s)^{1\over\b-1}ds}z(\t)^{\b\over\b-1}\n(t,\t)d\t,\qq t\in[0,T].\ee
We have the following result.

\ms

\bf Proposition 6.2. \sl Let $\n:D[0,T]\to(0,\b]$ and
$\rho:[0,T]\to(0,\infty)$ be continuous, and
\bel{6.19}\bar\l\equiv\sup_{0\le t<s\le T}{-\ln\n(t,s)\over
s-t}<\infty.\ee
Suppose $z:[0,T]\to(0,\infty)$ is a solution to $(\ref{z})$. Then
\bel{}e^{(\l-\bar\l)(T-t)}\min_{t\in[0,T]}\rho(t)\le z(t)\le
e^{\l(T-t)}\max_{t\in[0,T]}\rho(t),\qq t\in[0,T],\ee

\it Proof. \rm We may write (\ref{z}) as follows:
$$z(t)e^{-\l(T-t)+\b\int_t^Tz(s)^{1\over\b-1}ds}=\rho(t)
+\int_t^T\n(t,\t)z(\t)^{1\over\b-1}\[z(\t)e^{-\l(T-\t)+\b\int_\t^Tz(s)^{1\over\b-1}ds}\]d\t.$$
Denoting
$$\h z(t)=z(t)e^{-\l(T-t)+\b\int_t^Tz(s)^{1\over\b-1}ds},\qq t\in[0,T],$$
we have
\bel{hz}\h z(t)=\rho(t)+\int_t^T\n(t,\t)z(\t)^{1\over\b-1}\h
z(\t)d\t,\qq t\in[0,T].\ee
Note that
\bel{}\rho_0\equiv\min_{t\in[0,T]}\rho(t)\le\rho(t)\le\rho_1\equiv\max_{t\in[0,T]}\rho(t),
\qq0<\n(t,\t)\le\b.\ee
Thus, (\ref{hz}) implies
$$\h z(t)\le\rho_1+\b\int_t^Tz(\t)^{1\over\b-1}\h
z(\t)d\t,\qq t\in[0,T].$$
Then by Gronwall's inequality, we obtain
$$z(t)e^{-\l(T-t)+\b\int_t^Tz(s)^{1\over\b-1}ds}\equiv\h z(t)\le
\rho_1e^{\b\int_t^Tz(s)^{1\over\b-1}ds},\qq t\in[0,T],$$
which leads to
$$z(t)\le\rho_1e^{\l(T-t)},\qq t\in[0,T].$$
Next, (\ref{6.19}) implies
\bel{}\n(t,s)\ge e^{-\bar\l(s-t)},\qq(t,s)\in D[0,T].\ee
Consequently, we obtain from (\ref{hz}) that
\bel{hz2}\h
z(t)\ge\rho_0+\int_t^Te^{-\bar\l(\t-t)}z(\t)^{1\over\b-1}\h
z(\t)d\t,\qq t\in[0,T],\ee
which is equivalent to the following:
$$\h z(t)e^{-\bar\l t}\ge\rho_0e^{-\bar\l t}+\int_t^Tz(\t)^{1\over\b-1}\big[\h
z(\t)e^{-\bar\l\t}\big]d\t\equiv\z(t).$$
Then
$$\z'(t)=-\bar\l\rho_0e^{-\bar\l t}-z(t)^{1\over\b-1}\big[\h z(t)e^{-\bar\l
t}\big]\le-\bar\l\rho_0e^{-\bar\l t}-z(t)^{1\over\b-1}\z(t),$$
which yields
$$\[\z(t)e^{-\int_t^Tz(s)^{1\over\b-1}ds}\]'\le-\bar\l\rho_0e^{-\bar\l t-\int_t^Tz(s)^{1\over\b-1}ds}.$$
Hence,
$$\rho_0e^{-\bar\l
T}-\z(t)e^{-\int_t^Tz(s)^{1\over\b-1}ds}\le-\bar\l\rho_0\int_t^Te^{-\bar\l\t}
e^{-\int_\t^Tz(s)^{1\over\b-1}ds}d\t$$
Then
$$\ba{ll}
\ns\ds z(t)=\h z(t)e^{\l(T-t)}e^{-\int_t^Tz(s)^{1\over\b-1}ds}\ge
e^{\bar\l t}e^{\l(T-t)}\big[\z(t)e^{-\int_t^Tz(s)^{1\over\b-1}ds}\big]\\
\ns\ds\qq\ge e^{\bar\l t}e^{\l(T-t)}\[\rho_0e^{-\bar\l
T}+\bar\l\rho_0\int_t^Te^{-\bar\l\t}e^{-\int_\t^Tz(s)^{1\over\b-1}ds}d\t\]
\ge\rho_0e^{(\l-\bar\l)(T-t)},\qq t\in[0,T].\ea$$
This proves our proposition. \endpf

\ms

Having the above proposition, we then can easily obtain the
well-posedness of (\ref{f(t,t)}) by means of contraction mapping
theorem (giving the local solvability) and a usual continuation
argument. Then time-consistent equilibrium control can be
constructed by (\ref{uc}).

\ms

Note that by multiplying a constant to the payoff functional, if
necessary, we can always make
$$\max_{(\t,t)\in D[0,T]}\n(\t,t)\le\b.$$
On the other hand, in the case that $s\mapsto\n(s,t)$ is
differentiable, we have
$$\lim_{s\da t}{-\ln\n(t,s)\over s-t}=-{\n_s(t,t)\over\n(t,t)}.$$
Thus, condition (\ref{6.19}) is ensured by the boundedness of the
right hand side of the above, which is not too restrictive.

\section{Appendix}

\rm

\ms

In this appendix, we present some detailed calculations.

\ms

\bf Example 2.1. \rm Recall that we are considering the following
one-dimensional controlled linear SDE:
\bel{ex2.1a0}\left\{\ba{ll}
\ns\ds dX(s)=u(s)ds+\si X(s)dW(s),\qq s\in[t,T],\\
\ns\ds X(t)=x,\ea\right.\ee
with cost functional
\bel{ex2.1b0}J(t,x;u(\cd))=\dbE_t\[\int_t^T|u(s)|^2ds+g(t)|X(T)|^2\],\ee
where $\si>0$ is a constant and $g(t)$ is a deterministic
non-constant, continuous and positive function. For such a
linear-quadratic optimal control problem on $[t,T]$ (with
deterministic coefficients) and with $t\in[0,T)$ fixed, the Riccati
equation takes the following form: (note that $t\in[0,T)$ is a
parameter)
$$\left\{\ba{ll}
\ns\ds P_s(s,t)=P(s,t)^2-\si^2P(s,t),\qq s\in[t,T],\\
\ns\ds P(T,t)=g(t).\ea\right.$$
Let us solve the above Riccati equation. By separation of variables,
we have
$$ds={dP\over P^2-\si^2P}={1\over\si^2}\({1\over P-\si^2}-{1\over P}\)dP.$$
Integrating from $s$ to $T$, one has
$$\si^2(T-s)=\ln\({g(t)-\si^2\over g(t)}\)-\ln\({P(s,t)-\si^2\over
P(s,t)}\)=\ln\({g(t)-\si^2\over g(t)}{P(s,t)\over P(s,t)-\si^2}\).$$
Then
$${P(s,t)[g(t)-\si^2]\over[P(s,t)-\si^2]g(t)}=e^{\si^2(T-s)}.$$
Hence,
$$P(s,t)={\si^2g(t)e^{\si^2(T-s)}\over\si^2+g(t)(e^{\si^2(T-s)}-1)},\qq s\in[t,T],$$
and the optimal control is given by
$$\bar u(s)=-P(s,t)\bar X(s),\qq s\in[t,T].$$
Thus, the closed-loop SDE reads
$$\left\{\ba{ll}
\ns\ds d\bar X(s)=-P(s,t)\bar X(s)ds+\si\bar X(s)dW(s),\qq
s\in[t,T],\\
\ns\ds\bar X(t)=x.\ea\right.$$
Consequently, the optimal state process is given by
$$\ba{ll}
\ns\ds\bar
X(s)=xe^{-\int_t^s[P(\th,t)+{\si^2\over2}]d\th+\si[W(s)-W(t)]}\\
\ns\ds\qq~={\si^2+g(t)(e^{\si^2(T-s)}-1)\over\si^2+g(t)
(e^{\si^2(T-t)}-1)}\,e^{-{\si^2\over2}(s-t)+\si[W(s)-W(t)]}x,\qq
s\in[t,T].\ea$$
since
$$\ba{ll}
\ns\ds-\int_t^sP(\th,t)d\th=-\int_t^s{\si^2g(t)e^{\si^2(T-\th)}\over\si^2
+g(t)(e^{\si^2(T-\th)}-1)}\,d\th\\
\ns\ds=\int_t^s{d[\si^2+g(t)(e^{\si^2(T-\th)}-1)]\over\si^2+g(t)(e^{\si^2(T-\th)}-1)}=
\ln\({\si^2+g(t)(e^{\si^2(T-s)}-1)\over\si^2+g(t)(e^{\si^2(T-t)}-1)}\).\ea$$
Consequently, the optimal control also admits the following
open-loop form:
$$\bar u(s)=\bar u(s;t,x)=-P(s,t)\bar X(s)=-{\si^2g(t)e^{\si^2(T-s)}
e^{-{\si^2\over2}(s-t)+\si[W(s)-W(t)]}x\over\si^2+g(t)(e^{\si^2(T-t)}-1)}.$$
Recall that
$$\dbE_t\[e^{k[W(s)-W(t)]}\]=e^{{k^2\over2}(s-t)},\qq s\ge
t\ge0.$$
In fact, for fixed $t$, applying It\^o's formula to $s\mapsto
\f(s)\equiv e^{k[W(s)-W(t)]}$, we have
$$\left\{\ba{ll}
\ns\ds d\f(s)=k\f(s)dW(s)+{k^2\over2}\f(s)ds,\qq s\ge
t,\\
\ns\ds\f(t)=1.\ea\right.$$
Hence,
$$\left\{\ba{ll}
\ns\ds
d\{\dbE_t[\f(s)]\}={k^2\over2}\dbE_t[\f(s)],\qq s\ge t,\\
\ns\ds\dbE_t[\f(t)]=1,\ea\right.$$
which leads to what we want. We can check that
$$\ba{ll}
\ns\ds J(t,x;\bar u(\cd))=\dbE_t\[\int_t^T|\bar u(s)|^2ds+g(t)|\bar
X(T)|^2\]
%
%
%
%
%
%
%
={\si^2g(t)e^{\si^2(T-t)}x^2\over
\si^2+g(t)(e^{\si^2(T-t)}-1)}=P(t,t)x^2.\ea$$

Next, let $\t\in(t,T)$. We consider the corresponding LQ problem
starting from the initial pair $(\t,\bar X(\t))$. Then the optimal
state process, denoted by $\h X(\cd)$, must be given by
$$\ba{ll}
\ns\ds\h X(s)={\si^2+g(\t)(e^{\si^2(T-s)}-1)\over\si^2+g(\t)
(e^{\si^2(T-\t)}-1)}\,e^{-{\si^2\over2}(s-\t)+\si[W(s)-W(\t)]}\bar
X(\t),\qq s\in[\t,T],\ea$$
and the optimal control, denoted by $\h u(\cd)$, should be given by
$$\ba{ll}
\ns\ds\h u(s)=-P(s,\t)\h
X(s)=-{\si^2g(\t)e^{\si^2(T-s)}e^{-{\si^2\over2}(s-\t)+\si[W(s)-W(\t)]}\bar
X(\t)\over\si^2+g(\t)(e^{\si^2(T-\t)}-1)}\,,\qq s\in[\t,T].\ea$$
Note that
$$J(\t,\bar X(\t);\h u(\cd))=P(\t,\t)|\bar
X(\t)|^2={\si^2g(\t)e^{\si^2(T-\t)}\over\si^2+g(\t)(e^{\si^2(T-\t)}-1)}|\bar
X(\t)|^2.$$
On the other hand,
$$\bar X(\t)={\si^2+g(t)(e^{\si^2(T-\t)}-1)\over\si^2+g(t)
(e^{\si^2(T-t)}-1)}\,e^{-{\si^2\over2}(\t-t)+\si[W(\t)-W(t)]}x.$$
Thus, for $s\in[\t,T]$,
$$\ba{ll}
\ns\ds\bar X(s)={\si^2+g(t)(e^{\si^2(T-s)}-1)\over\si^2+g(t)
(e^{\si^2(T-t)}-1)}\,e^{-{\si^2\over2}(s-t)+\si[W(s)-W(t)]}x\\
\ns\ds\qq={\si^2+g(t)(e^{\si^2(T-s)}-1)\over\si^2+g(t)
(e^{\si^2(T-\t)}-1)}\,e^{-{\si^2\over2}(s-\t)+\si[W(s)-W(\t)]}\bar
X(\t),\ea$$
and
$$\ba{ll}
\ns\ds\bar u(s)=-P(s,t)\bar
X(s)\\
\ns\ds=-{\si^2g(t)e^{\si^2(T-s)}\over\si^2+g(t)(e^{\si^2(T-s)}-1)}
\cd{\si^2+g(t)(e^{\si^2(T-s)}-1)\over\si^2+g(t)
(e^{\si^2(T-\t)}-1)}e^{-{\si^2\over2}(s-\t)+\si[W(s)-W(\t)]}\bar X(\t)\\
\ns\ds=-{\si^2g(t)e^{\si^2(T-s)}\over
\si^2+g(t)(e^{\si^2(T-\t)}-1)}\,e^{-{\si^2\over2}(s-\t)+\si[W(s)-W(\t)]}\bar
X(\t).\ea$$
Then
$$\ba{ll}
\ns\ds J(\t,\bar X(\t);\bar u(\cd))=\dbE_\t\[\int_\t^T|\bar
u(s)|^2ds+g(\t)|\bar X(T)|^2\]\\
\ns\ds=\dbE_\t\[\int_\t^T{\si^4g(t)^2e^{2\si^2(T-s)-\si^2(s-\t)
+2\si[W(s)-W(\t)]}|\bar
X(\t)|^2\over[\si^2+g(t)(e^{\si^2(T-\t)}-1)]^2}\,ds\\
\ns\ds\qq\qq+{\si^4g(\t)e^{-\si^2(T-\t)+2\si[W(T)-W(\t)]}|\bar
X(\t)|^2\over[\si^2+g(t)(e^{\si^2(T-\t)}-1)]^2}\]\\
\ns\ds={\si^4g(t)^2|\bar
X(\t)|^2\over[\si^2+g(t)(e^{\si^2(T-\t)}-1)]^2}
\int_\t^Te^{2\si^2(T-s)+\si^2(s-\t)}ds+{\si^4g(\t)e^{\si^2(T-\t)}|\bar
X(\t)|^2\over[\si^2+g(t)(e^{\si^2(T-\t)}-1)]^2}\\
\ns\ds={\si^2g(t)^2|\bar
X(\t)|^2e^{\si^2(T-\t)}(e^{\si^2(T-\t)}-1)+\si^4g(\t)e^{\si^2(T-\t)}|\bar
X(\t)|^2\over[\si^2+g(t)(e^{\si^2(T-\t)}-1)]^2}\\
\ns\ds={[g(t)^2(e^{\si^2(T-\t)}-1)+\si^2g(\t)]\si^2e^{\si^2(T-\t)}|\bar
X(\t)|^2\over[\si^2+g(t)(e^{\si^2(T-\t)}-1)]^2}.\ea$$
Consequently, the following holds:
$$\ba{ll}
\ns\ds J(\t,\bar X(\t);\bar u(\cd))-J(\t,\bar X(\t);\h u(\cd))\\
\ns\ds=\Big\{{g(t)^2(e^{\si^2(T-\t)}-1)+\si^2g(\t)\over
[\si^2+g(t)(e^{\si^2(T-\t)}-1)]^2}-{g(\t)\over\si^2+g(\t)(e^{\si^2(T-\t)}-1)}
\Big\}\si^2e^{\si^2(T-\t)}|\bar X(\t)|^2\\
\ns\ds={\si^4(e^{\si^2(T-\t)}-1)e^{\si^2(T-\t)}[g(t)-g(\t)]^2|\bar
X(\t)|^2\over[\si^2+g(t)(e^{\si^2(T-\t)}-1)]^2[\si^2+g(\t)(e^{\si^2(T-\t)}-1)]}\\
%
%
\ns\ds={\si^4(e^{\si^2(T-\t)}-1)[g(t)-g(\t)]^2e^{\si^2(T-\t)-\si^2(\t-t)+2\si[W(\t)-W(t)]}x^2
\over[\si^2+g(t)(e^{\si^2(T-t)}-1)]^2[\si^2+g(\t)(e^{\si^2(T-\t)}-1)]}
\,,\ea$$
which is strictly positive unless $x=0$, or $g(\t)=g(t)$. This means
that the problem is time-inconsistent.

\ms

\bf Example 2.2. (Generalized Merton's portfolio problem) \rm Recall
the following controlled SDE:
\bel{State1}\left\{\ba{ll}
\ns\ds dX(s)=\big[rX(s)+(\m-r)u(s)-c(s)\big]ds+\si u(s)dW(s),\qq
s\in[t,T],\\
\ns\ds X(t)=x,\ea\right.\ee
with payoff functional
\bel{6.4}J(t,x;u(\cd),c(\cd))=\dbE_t\[\int_t^T\n(t,s)c(s)^\b
ds+\rho(t)X\big(T;t,x,u(\cd),c(\cd)\big)^\b\],\ee
where $\n(\cd\,,\cd)$ and $\rho(\cd)$ are given positive-valued
functions, and $\b\in(0,1)$. As a convention, we define
$$x^\b=-\infty,\qq x<0.$$
The optimal control problem is to find a pair $(\bar u(\cd),\bar
c(\cd))$ such that $J(t,x;u(\cd),c(\cd))$ is maximized. To approach
this problem, we use dynamic programming method. More precisely, for
any $(s,y)\in[t,T)\times[0,\infty)$, let
\bel{}J^t\big(s,y;u(\cd),c(\cd)\big)=\dbE_t\[\int_s^T\n(t,\t)c(\t)^\b
d\t+\rho(t)X\big(T;s,y, u(\cd),c(\cd)\big)^\b\].\ee
Define the value function $V^t(\cd\,,\cd)$ (which is parameterized
by $t\in[0,T]$) by the following:
\bel{}V^t(s,y)=\sup_{(u(\cd),c(\cd))}J^t\big(s,y;u(\cd),c(\cd)\big),\qq(s,y)\in[t,T]\times[0,\infty).\ee
Due to the problem being a maximization problem, the above
convention forces $X(s)$ to stay nonnegative, in particular, the
initial state $x\ge0$ has to be assumed. In another word, $V^t(s,x)$
is only defined on $[0,T]\times[0,\infty)$. By Girsanov's theorem,
we know that
$$\wt W(s)={\m-r\over\si}s+\si W(s),\qq s\ge0$$
is a standard Brownian motion, with the natural filtration
(augmented by all the $\dbP$-null sets) coincides with
$\dbF\equiv\{\cF_t\}_{t\ge0}$. Then (\ref{State1}) is equivalent to
the following:
\bel{6.7}\left\{\ba{ll}
\ns\ds dX(s)=\big[rX(s)-c(s)\big]ds+\si u(s)d\wt W(s),\qq
s\in[t,T],\\
\ns\ds X(t)=x,\ea\right.\ee
Therefore, for any initial pair $(t,x)\in[0,T)\times[0,\infty)$, and
a control $(u(\cd),c(\cd))$, the unique solution $X(\cd)\equiv
X(\cd\,;t,x,u(\cd),c(\cd))$ of the above is given by the following:
\bel{}X(s)=e^{r(s-t)}x-\int_t^se^{r(s-\t)}c(\t)d\t+\si\int_t^se^{r(s-\t)}u(\t)d\wt
W(\t),\qq s\in[t,T].\ee
If the initial wealth $x=0$, to keep the wealth $X(\cd)$
non-negative, we have to take $u(\cd)=c(\cd)=0$, leading to
$$X(s)=0,\qq s\in[t,T].$$
This means that
\bel{}V^t(t,0)=0,\qq t\in[0,T],\ee
which gives the boundary condition for $V^t(\cd\,,\cd)$ on $x=0$.
Now, let us return to (\ref{State1})--(\ref{6.4}). In the case that
$V^t(\cd\,,\cd)$ is differentiable, it satisfies the following HJB
equation:
\bel{}\ba{ll}
\ns\ds0=V^t_s(s,y)+\sup_{(u,c)}\[V^t_y(s,y)\big[ry+(\m-r)u-c\big]+{1\over2}
\si^2u^2
V^t_{yy}(s,y)+\n(t,s)c^\b\]\\
\ns\ds\q=V^t_s(s,y)+ryV^t_y(s,y)+\sup_{u\in\dbR}\[(\m-r)
V^t_y(s,y)u+{1\over2}\si^2V^t_{yy}(s,y)u^2\]+\sup_{c\ge0}\big[\n(t,s)c^\b-cV^t_y(s,y)\big].\ea\ee
Assume, for the time being, that the following holds:
\bel{V_yy<0}V^t_y(s,y)>0,\q
V^t_{yy}(s,y)<0,\qq(s,y)\in[t,T]\times(0,\infty).\ee
Then
\bel{}\sup_{u\in\dbR}\[(\m-r)V^t_y(s,y)u+{1\over2}\si^2V^t_{yy}(s,y)u^2\]=-{(\m-r)^2
V^t_y(s,y)^2\over2\si^2V^t_{yy}(s,y)}>0,\ee
with the maximum attained at
\bel{}\bar u^t
(s,y)=-{(\m-r)V^t_y(s,y)\over\si^2V^t_{yy}(s,y)}>0,\ee
and
\bel{}\sup_{c>0}\big[\n(t,s)c^\b-cV^t_y(s,y)\big]=(1-\b)\b^{\b\over1-\b}
\n(t,s)^{1\over1-\b}V^t_y(s,y)^{\b\over\b-1}>0,\ee
with the maximum attained at
\bel{}\bar c^t(s,y)=\({\b\n(t,s)\over
V^t_y(s,y)}\)^{1\over1-\b}>0.\ee
Consequently, the HJB equation reads
\bel{HJB2}\left\{\ba{ll}
\ns\ds
V^t_s(s,y)+ryV^t_y(s,y)-{(\m-r)^2V^t_y(s,y)^2\over2\si^2V^t_{yy}(s,y)}
+(1-\b)\b^{\b\over1-\b}\n(t,s)^{1\over1-\b}V^t_y(s,y)^{\b\over\b-1}=0,\\
\ns\ds\qq\qq\qq\qq\qq\qq\qq\qq\qq(s,y)\in[t,T)\times(0,\infty),\\
\ns\ds V^t(s,0)=0,\qq\qq\qq s\in[t,T],\\
\ns\ds V^t(T,y)=\rho(t)y^\b,\qq\qq y\in(0,\infty).\ea\right.\ee
We try to find the solution of the following form:
\bel{}V^t(s,y)=\f(s)y^\b,\qq(s,y)\in[t,T]\times[0,\infty).\ee
Clearly, with such a form, (\ref{V_yy<0}) is satisfied. Then we
should have
$$\ba{ll}
\ns\ds0=\f'(s)y^\b+r\b\f(s)y^\b+{(\m-r)^2\b
\f(s)y^\b\over2\si^2(1-\b)}
+(1-\b)\n(t,s)^{1\over1-\b}\f(s)^{\b\over\b-1}y^\b\ea$$
This leads to an ordinary differential equation for $\f(\cd)$:
\bel{ODE}\left\{\ba{ll}
\ns\ds\f'(s)+\[r\b+{(\m-r)^2\b\over2\si^2(1-\b)}\]\f(s)
+(1-\b)\n(t,s)^{1\over1-\b}\f(s)^{\b\over\b-1}=0,\q s\in[t,T],\\
\ns\ds\f(T)=\rho(t).\ea\right.\ee
If we denote
\bel{l7}\l=r\b+{(\m-r)^2\b\over2\si^2(1-\b)}={[2r\si^2(1-\b)+(\m-r)^2]\b\over2\si^2(1-\b)},\ee
then (\ref{ODE}) becomes
\bel{ODE2}\left\{\ba{ll}
\ns\ds\f'(s)+\l\f(s)+(1-\b)\n(t,s)^{1\over1-\b}\f(s)^{\b\over\b-1}=0,\qq s\in[t,T],\\
\ns\ds\f(T)=\rho(t).\ea\right.\ee
This is a Bernoulli equation. To solve it, let
\bel{}\psi(s)=\f(s)^{1\over1-\b},\qq s\in[t,T].\ee
Then
\bel{}\ba{ll}
\ns\ds\psi'(s)={1\over1-\b}\f(s)^{\b\over1-\b}\f'(s)=-{1\over1-\b}
\f(s)^{\b\over1-\b}
\(\l\f(s)+(1-\b)\n(t,s)^{1\over1-\b}\f(s)^{\b\over\b-1}\)\\
\ns\ds\qq~=-{\l\over1-\b}\psi(s)-\n(t,s)^{1\over1-\b},\qq
s\in[t,T],\ea\ee
with
\bel{}\psi(T)=\rho(t)^{1\over1-\b}.\ee
Hence,
\bel{}\psi(s)=e^{{\l\over1-\b}(T-s)}\rho(t)^{1\over1-\b}+\int_s^Te^{{\l\over1-\b}(\t-s)}\n(t,\t)^{1\over1-\b}d\t,
\qq s\in[t,T].\ee
Then
\bel{}\f(s)=\[e^{{\l\over1-\b}(T-s)}\rho(t)^{1\over1-\b}+\int_s^Te^{{\l\over1-\b}(\t-s)}\n(t,\t)^{1\over1-\b}d\t\]^{1-\b},
\qq s\in[t,T].\ee
Consequently,
\bel{}V^t(s,y)=\[e^{{\l\over1-\b}(T-s)}\rho(t)^{1\over1-\b}+\int_s^Te^{{\l\over1-\b}(\t-s)}\n(t,\t)^{1\over1-\b}d\t\]^{1-\b}y^\b,\qq(s,y)\in[t,T]\times[0,\infty).\ee
Therefore, the optimal control $(\bar u^t(\cd),\bar c^t(\cd))$ for
the initial pair $(t,x)$ is given by the following:
\bel{}\bar u^t(s)=-{(\m-r)V^t_y(s,\bar
X^t(s))\over\si^2[V^t_{yy}(s,\bar X^t(s))]}={(\m-r)\bar
X^t(s)\over\si^2 (1-\b)}, \qq s\in[t,T],\ee
which is not depending on the parameter $t$ directly, and
\be{}\ba{ll}
\ns\ds\bar c^t(s)=\b^{1\over1-\b}\n(t,s)^{1\over1-\b}V^t_y(s,\bar
X^t(s))^{1\over\b-1} =\n(t,s)^{1\over1-\b}
\f(s)^{1\over\b-1}\bar X^t(s)\\
\ns\ds\qq\;={\n(t,s)^{1\over1-\b}\bar X^t(s)\over
e^{{\l\over1-\b}(T-s)}\rho(t)^{1\over1-\b}
+\int_s^Te^{{\l\over1-\b}(\t-s)}\n(t,\t)^{1\over1-\b}d\t},\qq\q
s\in[t,T],\ea\ee
which is depending on the parameter $t$ directly. Combining the
above, we have
\bel{}\ba{ll}
\ns\ds\sup_{u(\cd),c(\cd)}J(t,x;u(\cd),c(\cd))=V(t,x)\equiv
V^t(t,x)\\
\ns\ds\qq=\[e^{{\l\over1-\b}(T-t)}\rho(t)^{1\over1-\b}
+\int_t^Te^{{\l\over1-\b}(\t-t)}\n(t,\t)^{1\over1-\b}
d\t\]^{1-\b}x^\b\\
\ns\ds\qq=\dbE_t\[\int_t^T\n(t,\t)\bar c^t(\t)^\b d\t+\rho(t)\bar
X^t(T)^\b\],\qq(t,x)\in[0,T]\times[0,\infty).\ea\ee
Now, let $\bar t\in(t,T)$. Then, applying the above argument, one
has
\bel{bar V}\ba{ll}
\ns\ds V(\bar t,\bar X^t(\bar t))=\[e^{{\l\over1-\b}(T-\bar
t)}\rho(\bar t)^{1\over1-\b}+\int_{\bar t}^Te^{{\l\over1-\b}(\t-\bar
t)}\n(\bar t,\t)^{1\over1-\b} d\t\]^{1-\b}\bar X^t(\bar t)^\b.\ea\ee
We claim that if the following is assumed
\bel{ne}\int_{\bar
t}^T\[{e^{\l\t}\n(t,\t)\over\rho(t)}\]^{1\over1-\b}d\t\ne\int_{\bar
t}^T\[{e^{\l\t}\n(\bar t,\t)\over\rho(\bar t)}\]^{1\over1-\b}d\t,\ee
then
\bel{J<V}J\big(\bar t,\bar X^t(\bar t);\bar u(\cd)\big|_{[\bar
t,T]},\bar c(\cd) \big|_{[\bar t,T]}\big)<V(\bar t,\bar X^t(\bar
t)).\ee
This means that the restriction $\big(\bar u(\cd)\big|_{[\bar
t,T]},\bar c(\cd)\big|_{[\bar t,T]}\big)$ of $(\bar u(\cd),\bar
c(\cd))$ on $[\bar t,T]$ is not optimal for the initial pair $(\bar
t,\bar X^t(\bar t))$. To prove our claim, let us denote
\bel{G7}\ba{ll}
\ns\ds\G(t,s)={\n(t,s)^{1\over1-\b}\over
e^{{\l\over1-\b}(T-s)}\rho(t)^{1\over1-\b}
+\int_s^Te^{{\l\over1-\b}(\t-s)}\n(t,\t)^{1\over1-\b}d\t}\\
\ns\ds\qq\q={e^{-{\l\over1-\b}(T-s)}\n(t,s)^{1\over1-\b}\over
\rho(t)^{1\over1-\b}+\int_s^Te^{-{\l\over1-\b}(T-\t)}\n(t,\t)^{1\over1-\b}d\t}\ea\ee
Then
$$\bar c^t(s)=\G(t,s)\bar X^t(s),\qq\si\bar u^t(s)={\m-r\over\si(1-\b)}\bar X^t(s),\qq s\in[t,T],$$
where $\bar X^t(\cd)$ is the solution to the following closed-loop
system
$$\left\{\ba{ll}
\ns\ds d\bar X^t(s)=\(r+{(\m-r)^2\over\si^2(1-\b)}-\G(t,s)\)\bar
X^t(s)ds+{\m-r\over\si(1-\b)}\bar X^t(s)dW(s),\qq
s\in[t,T],\\
\ns\ds\bar X^t(t)=x.\ea\right.$$
By denoting
$$b=r+{(\m-r)^2\over\si^2(1-\b)}\,,\qq a={\m-r\over\si(1-\b)}\,,$$
we may write the above as
$$\left\{\ba{ll}
\ns\ds d\bar X^t(s)=\big[b-\G(t,s)\big]\bar X^t(s)ds+a\bar
X^t(s)dW(s),\qq s\in[t,T],\\
\ns\ds\bar X^t(t)=x.\ea\right.$$
Hence,
$$\bar X^t(s)=\bar X^t(\bar t)e^{(b-{a^2\over2})(s-\bar t)-\int_{\bar t}^s\G(t,\t)d\t+a[
W(s)-W(\bar t)]},\qq s\in[\bar t,T].$$
Consequently,
$$\ba{ll}
\ns\ds\dbE_{\bar t}\big[\bar X^t(s)^\b\big]=\bar X^t(\bar t)^\b
e^{\b(b-{a^2\over2}) (s-\bar t)-\b\int_{\bar
t}^s\G(t,\t)d\t+{\b^2a^2\over2}(s-\bar t)}\\
\ns\ds\qq\qq\q\,=\bar X^t(\bar t)^\b e^{\l(s-\bar t)-\b\int_{\bar
t}^s\G(t,\t)d\t},\qq s\in[\bar t,T].\ea$$
Here, we note that (recall (\ref{l7}))
$$\b\(b-{a^2\over2}\)+{\b^2a^2\over2}=\b\(b-{a^2(1-\b)\over2}\)=\b\(r+{(\m-r)^2\over2\si^2(1-\b)}\)=\l.$$
Then
$$\ba{ll}
\ns\ds J(\bar t,\bar X^t(\bar t);\bar u(\cd)\big|_{[\bar t,T]},\bar
c(\cd)\big|_{[\bar t,T]})=\dbE_{\bar t}\[\int_{\bar t}^T\n(\bar
t,s)\bar
c^t(s)^\b ds+\rho(\bar t)\bar X^t(T)^\b\]\\
\ns\ds=\int_{\bar t}^T\n(\bar t,s)\G(t,s)^\b\dbE_{\bar t}[\bar
X^t(s)^\b]ds
+\rho(\bar t)\dbE_{\bar t}[\bar X^t(T)^\b]\\
\ns\ds=\[\int_{\bar t}^T\n(\bar t,s)\G(t,s)^\b e^{\l(s-\bar
t)-\b\int_{\bar t}^s\G(t,\t)d\t}ds+\rho(\bar t)e^{\l(T-\bar
t)-\b\int_{\bar t}^T\G(t,\t)d\t}\]\bar X^t(\bar t)^\b\\
\ns\ds=\[\int_{\bar t}^Te^{\l(s-\bar t)}\n(\bar t,s)\(\G(t,s)
e^{-\int_{\bar t}^s\G(t,\t)d\t}\)^\b ds+e^{\l(T-\bar t)}\rho(\bar
t)e^{-\b\int_{\bar t}^T\G(t,\t)d\t}\]\bar X^t(\bar t)^\b.\ea$$
By H\"older's inequality, one has
$$\ba{ll}
\ns\ds\int_{\bar t}^Te^{\l(s-\bar t)}\n(\bar t,s)\(\G(t,s)
e^{-\int_{\bar t}^s\G(t,\t)d\t}\)^\b ds+e^{\l(T-\bar t)}\rho(\bar
t)e^{-\b\int_{\bar t}^T\G(t,\t)d\t}\\
\ns\ds\le\[\int_{\bar t}^Te^{{\l\over1-\b}(s-\bar t)}\n(\bar
t,s)^{1\over1-\b}ds\]^{1-\b}\[\int_{\bar t}^T\G(t,s)e^{-\int_{\bar
t}^s\G(t,\t)d\t}ds\]^\b+e^{\l(T-\bar t)}\rho(\bar
t)e^{-\b\int_{\bar t}^T\G(t,\t)d\t}\\
\ns\ds\le\[\int_{\bar t}^Te^{{\l\over1-\b}(s-\bar t)}\n(\bar
t,s)^{1\over1-\b}ds\]^{1-\b}\[1-e^{-\int_{\bar
t}^T\G(t,\t)d\t}\]^\b+e^{\l(T-\bar t)}\rho(\bar
t)e^{-\b\int_{\bar t}^T\G(t,\t)d\t}\\
\ns\ds\equiv\a_1(\bar t)^{1-\b}\g(t,\bar t)^\b+\a_2(\bar
t)^{1-\b}\big[1-\g(t,\bar t)\big]^\b,\ea$$
where
$$\left\{\ba{ll}
\ns\ds\a_1(\bar t)=\int_{\bar t}^Te^{{\l\over1-\b}(s-\bar t)}\n(\bar
t,s)^{1\over1-\b}ds,\qq\a_2(\bar t)=e^{{\l\over1-\b}(T-\bar t)}\rho(\bar t)^{1\over1-\b}\ge0,\\
\ns\ds\g(t,\bar t)=1-e^{-\int_{\bar
t}^T\G(t,\t)d\t}\in[0,1].\ea\right.$$
To proceed further, we need the following elementary result.

\ms

\bf Lemma 7.1. \rm Let $\a_1,\a_2>0$, $\b\in(0,1)$, and
$$f(\g)=\a_1^{1-\b}\g^\b+\a_2^{1-\b}(1-\g)^\b,\qq\g\in[0,1].$$
Then $\g\mapsto f(\g)$ is strictly concave and
$$\max_{\g\in[0,1]}f(\g)=f\({\a_1\over\a_1+\a_2}\)=(\a_1+\a_2)^{1-\b}.$$

\it Proof. \rm Note that
$$f'(\g)=\b\[\a_1^{1-\b}\g^{\b-1}-\a_2^{1-\b}(1-\g)^{\b-1}\],$$
and
$$f''(\g)=-\b(1-\b)\[\a_1^{1-\b}\g^{\b-2}+\a_2^{1-\b}(1-\g)^{\b-2}\]<0,\qq\forall
\g\in(0,1).$$
Thus, $\g\mapsto f(\g)$ is strictly concave. By setting $f'(\g)=0$,
we have
$${\a_1\over\g}={\a_2\over1-\g}\,,$$
which gives the unique maximum:
$$\g={\a_1\over\a_1+\a_2}\,.$$
Clearly,
$$f\({\a_1\over\a_1+\a_2}\)=\a_1^{1-\b}{\a_1^\b\over(\a_1+\a_2)^\b}
+\a_2^{1-\b}{\a_2^\b\over(\a_1+\a_2)^\b}=(\a_1+\a_2)^{1-\b}.$$
This proves the lemma. \endpf

\ms

By the above lemmas, we obtain (note (\ref{bar V}))
$$\ba{ll}
\ns\ds J\big(\bar t,\bar X^t(\bar t);\bar u(\cd)\big|_{[\bar
t,T]},\bar
c(\cd)\big|_{[\bar t,T]}\big)\\
\ns\ds=\[\int_{\bar t}^Te^{\l(s-\bar t)}\n(\bar t,s)\(\G(t,s)
e^{-\int_{\bar t}^s\G(t,\t)d\t}\)^\b ds+e^{\l(T-\bar t)}\rho(\bar
t)e^{-\b\int_{\bar t}^T\G(t,\t)d\t}\]\bar X^t(\bar t)^\b\\
\ns\ds\le\[\a_1(\bar t)^{1-\b}\g(t,\bar t)^\b+\a_2(\bar
t)^{1-\b}[1-\g(t,\bar t)]^\b\]\bar X^t(\bar t)^\b\le[\a_1(\bar
t)+\a_2(\bar t)]^{1-\b}\bar
X^t(\bar t)^\b\\
\ns\ds=\[\int_{\bar t}^Te^{{\l\over1-\b}(s-\bar t)}\n(\bar
t,s)^{1\over1-\b} ds+e^{{\l\over1-\b}(T-\bar t)}\rho(\bar
t)^{1\over1-\b}\]^{1-\b}\bar X^t(\bar t)^\b=V(\bar t,\bar X^t(\bar
t)).\ea$$
In order to have a strict inequality in the above, it suffices to
have
$$\a_1(\bar t)^{1-\b}\g(t,\bar t)^\b+\a_2(\bar
t)^{1-\b}[1-\g(t,\bar t)]^\b<[\a_1(\bar t)+\a_2(\bar t)]^{1-\b},$$
which is implied by
$$\g(t,\bar t)\ne{\a_1(\bar t)\over\a_1(\bar t)+\a_2(\bar t)}\,.$$
This is equivalent to the following:
$$\ba{ll}
\ns\ds e^{-\int_{\bar t}^T\G(t,\t)d\t}\ne{\a_2(\bar t)\over\a_1(\bar
t)+\a_2(\bar t)}={e^{{\l\over1-\b}(T-\bar t)}\rho(\bar
t)^{1\over1-\b}\over e^{{\l\over1-\b}(T-\bar t)}\rho(\bar
t)^{1\over1-\b}+\int_{\bar
t}^Te^{{\l\over1-\b}(s-\bar t)}\n(\bar t,s)^{1\over1-\b}ds}\\
\ns\ds={\rho(\bar t)^{1\over1-\b}\over\rho(\bar
t)^{1\over1-\b}+\int_{\bar t}^Te^{-{\l\over1-\b}(T-s)}\n(\bar
t,s)^{1\over1-\b}ds}.\ea$$
Note that (recall (\ref{G7}))
$$\ba{ll}
\ns\ds-\int_{\bar t}^T\G(t,s)ds=-\int_{\bar
t}^T{e^{-{\l\over1-\b}(T-s)}\n(t,s)^{1\over1-\b}\over\rho(t)^{1\over1-\b}
+\int_s^Te^{-{\l\over1-\b}(T-\t)}\n(t,\t)^{1\over1-\b}d\t}\,ds\\
\ns\ds=\int_{\bar
t}^T{d[\rho(t)^{1\over1-\b}+\int_s^Te^{-{\l\over1-\b}(T-\t)}\n(t,\t)^{1\over1-\b}
d\t]
\over\rho(t)^{1\over1-\b}+\int_s^Te^{-{\l\over1-\b}(T-\t)}\n(t,\t)^{1\over1-\b}
d\t}\\
\ns\ds=\ln\[\rho(t)^{1\over1-\b}+\int_s^Te^{-{\l\over1-\b}(T-\t)}\n(t,\t)^{1\over1-\b}
d\t\]\Big|_{\bar
t}^T=\ln{\rho(t)^{1\over1-\b}\over\rho(t)^{1\over1-\b}+\int_{\bar
t}^Te^{-{\l\over1-\b}(T-\t)}\n(t,\t)^{1\over1-\b} d\t}\,.\ea$$
Hence, we need
$${\rho(t)^{1\over1-\b}\over\rho(t)^{1\over1-\b}+\int_{\bar
t}^Te^{-{\l\over1-\b}(T-\t)}\n(t,\t)^{1\over1-\b}d\t}\ne {\rho(\bar
t)^{1\over1-\b}\over\rho(\bar t)^{1\over1-\b}+\int_{\bar
t}^Te^{-{\l\over1-\b}(T-\t)}\n(\bar t,\t)^{1\over1-\b}d\t}\,.$$
which is equivalent to (\ref{ne}), proving our claim.

\ms

Note that in the case
\bel{7.34}\n(t,s)=e^{-\d(s-t)},\qq\rho(t)=e^{-\d(T-t)},\qq0\le t\le
s\le T,\ee
we have
$${\n(t,\t)\over\rho(t)}={e^{-\d(\t-t)}\over
e^{-\d(T-t)}}=e^{\d(T-\t)}.$$
Thus, (\ref{ne}) will not be true. When (\ref{7.34}) holds, the
problem is referred to as the (classical) Merton's portfolio
problem. In this case, (with $a={\l-\d\over1-\b}$)
\be{}\ba{ll}
\ns\ds\bar c(s,y)={e^{-{\d\over1-\b}(s-t)}y\over
e^{\k(T-s)}e^{-{\d\over1-\b}(T-t)}+\int_s^Te^{\k(\t-s)}e^{-{\d\over1-\b}(\t-t)}
d\t}\\
\ns\ds\qq\q={y\over e^{\k(T-s)}e^{-{\d
(T-s)\over1-\b}}+\int_s^Te^{\k(\t-s)}e^{-{\d(\t-s)\over1-\b}}d\t}\\
\ns\ds\qq\q={y\over e^{a(T-s)}+\int_s^Te^{a(\t-s)}d\t}={ay\over
ae^{a(T-s)}+e^{a(T-s)}-1},\ea\ee
which is independent of $t$. This recovers the solution to the
classical Merton's portfolio problem.

\ms

\end{document}